\documentclass[11pt,reqno]{amsart}
\usepackage[T1]{fontenc}
\usepackage{lmodern,amsmath,amssymb,mathtools,mathrsfs,stmaryrd}
\usepackage[margin=1.05in]{geometry}
\usepackage{microtype,enumitem}
\usepackage[hidelinks]{hyperref}
\newtheorem{theorem}{Theorem}[section]
\newtheorem{proposition}[theorem]{Proposition}
\newtheorem{lemma}[theorem]{Lemma}

\theoremstyle{definition}

\theoremstyle{remark}

\numberwithin{equation}{section}
\newcommand{\R}{\mathbb R}
\newcommand{\C}{\mathbb C}
\newcommand{\spt}{\operatorname{spt}}
\newcommand{\Sing}{\operatorname{Sing}}
\newcommand{\Graph}{\operatorname{Graph}}
\newcommand{\osc}{\operatorname{osc}}
\newcommand{\tr}{\operatorname{tr}}
\newcommand{\adj}{\operatorname{adj}}
\newcommand{\dist}{\operatorname{dist}}
\newcommand{\A}{\mathcal A}
\newcommand{\Gr}{\operatorname{Gr}}
\newcommand{\res}{\mathbin{\llcorner}}
\newcommand{\Iplane}{\mathcal I}
\newcommand{\Pplane}{\mathcal P}
\DeclareMathOperator{\diag}{diag}
\DeclareMathOperator{\degmap}{deg}
\allowdisplaybreaks[2]
\title[Entire special Lagrangian graphs with critical phase]
{Nonpolynomial entire solutions of the\\
special Lagrangian equation with critical phase}
\author{Caiyan Li}
\address{(C.L.) School of Mathematical Sciences, Xiamen University, 361005, Xiamen, P.~R.~China}
\email{caiyanli@xmu.edu.cn} 
\thanks{C.L. is supported by NSFC (Grant No. 12501274).}
\date{}
\begin{document}
\begin{abstract}
In this paper, we construct smooth nonpolynomial entire solutions of the special Lagrangian equation with critical phase $(n-2)\pi/2$  for every $n\ge4$.
Together with the nonquadratic example in dimension two and Warren's  example in dimension three, our construction shows that the
supercritical condition
$|\theta|>(n-2)\pi/2$ in the Bernstein theorem proved by Yuan is sharp
in every dimension $n\ge2$.
\end{abstract}
\maketitle

\section{Introduction}\label{sec:intro}

In this paper, we consider the special Lagrangian equation
\begin{equation}\label{eq:special-lagrangian-equation}
\Theta(D^2u)=\arctan\lambda_1+\cdots+\arctan\lambda_n
=\theta\qquad\text{in }\R^n,
\end{equation}
where $\lambda_1,\ldots,\lambda_n$ are the eigenvalues of $D^2u$ and
$\theta$ is constant. More generally,
$\Theta(A)\coloneqq\sum_{i=1}^n\arctan\lambda_i(A)$ for a real symmetric
$n\times n$ matrix $A$.

Equation \eqref{eq:special-lagrangian-equation} arises from the theory of calibrated geometries developed by Harvey and Lawson~\cite{HL}. By symmetry of $D^2u$, the standard symplectic form on $\C^n$ vanishes along the gradient graph $(x,Du(x))\subset\R^n\times\R^n=\C^n$, which is therefore Lagrangian. The phase of the gradient graph $(x,Du(x))$ is the argument of the complex number $(1+i\lambda_1)\cdots(1+i\lambda_n)$, taken as $\sum_j\arctan\lambda_j$. The induced metric on $(x,Du(x))$ is $g=I+(D^2u)^2$, and the corresponding volume element is $dv_g=\sqrt{\det(I+(D^2u)^2)}\,dx$. If $u$ satisfies \eqref{eq:special-lagrangian-equation}, then the gradient graph is special Lagrangian with phase $\theta$. By Harvey--Lawson~\cite[Theorem~2.3 and Proposition~2.17]{HL}, the special Lagrangian graph $(x,Du(x))$, with its natural orientation, is calibrated by $\operatorname{Re}(e^{-i\theta}dz_1\wedge\cdots\wedge dz_n)$ and hence is volume minimizing, where $z_1,\ldots,z_n$ are the standard complex coordinates on $\C^n$. The Bernstein problem for \eqref{eq:special-lagrangian-equation} is to determine under what conditions an entire solution $u$ must be a quadratic polynomial, or equivalently, the gradient graph $(x,Du(x))$ must be an affine plane.

For $n\ge2$, Yuan~\cite[Lemma~2.1 and Remark~1]{YuanSupercritical}
showed that the nonempty level sets
$\{\lambda\in\R^n:\sum_i\arctan\lambda_i=\theta\}$ are convex
hypersurfaces exactly when $|\theta|\ge(n-2)\pi/2$. Accordingly, the phase
$\theta$ is called \textit{critical} if $|\theta|=(n-2)\pi/2$,
\textit{subcritical} if $|\theta|<(n-2)\pi/2$, and \textit{supercritical}
if $|\theta|>(n-2)\pi/2$.

Our main result is as follows.
\begin{theorem}\label{thm:main}
Let $n\ge4$. Then there exists a nonpolynomial
$u\in C^\infty(\R^n)$ satisfying
\[
\sum_{i=1}^n\arctan\lambda_i(D^2u)=\frac{n-2}{2}\pi
\quad\text{in }\R^n.
\]
\end{theorem}

In dimension two, Fu~\cite[Theorem 1.3]{Fu} proved that every entire
solution of \eqref{eq:special-lagrangian-equation} with nonzero phase
is quadratic. At the critical phase $\theta=0$, the equation is
$\Delta u=0$. For instance,
\[
u(x_1,x_2)=x_1^2x_2-\tfrac13x_2^3
\]
is a nonquadratic entire solution \cite{LiYuanExplicit}.
In arbitrary dimension, Yuan~\cite[Theorem 1.1]{YuanConvex} proved that
convex entire solutions of \eqref{eq:special-lagrangian-equation} are quadratic.
In dimension three, Yuan~\cite[Theorem 1.3]{YuanConvex} also proved
that an entire solution of \eqref{eq:special-lagrangian-equation} is
quadratic if $D^2u\ge-KI$ in $\R^3$ for some constant $K>0$.
In dimension four, the Bernstein result under $D^2u\ge-KI$ follows by combining the argument of
Yuan~\cite[pp.~124--125]{YuanConvex} with the rigidity theorem of
Nadirashvili--Vl\u{a}du\c{t}~\cite{NV13} for special Lagrangian cones given by
real analytic gradient graphs over $\R^4\setminus\{0\}$
(see Ogden--Yuan~\cite[p.~3]{OgdenYuanRank}).
In general dimensions, Ogden--Yuan~\cite[Theorem 1.2]{OgdenYuanRank}
proved that a smooth entire solution of
\eqref{eq:special-lagrangian-equation} is quadratic if
\[
\inf_{\R^n}\arctan\lambda_{\min}(D^2u)>\frac{\theta-\pi}{n}
\qquad\text{or}\qquad
\sup_{\R^n}\arctan\lambda_{\max}(D^2u)<\frac{\theta+\pi}{n}.
\]

Under the supercritical phase condition $|\theta|>(n-2)\pi/2$,
Yuan~\cite{YuanSupercritical}
proved that all entire solutions of \eqref{eq:special-lagrangian-equation}
are quadratic. 
At the critical phase $\theta=(n-2)\pi/2$,
Yuan~\cite[p.~1356]{YuanSupercritical} conjectured that nonquadratic
entire solutions exist in higher dimensions. For arbitrary subcritical phase $|\theta|<(n-2)\pi/2$,
Li--Yuan~\cite{LiYuanExplicit} construct entire explicit solutions of
\eqref{eq:special-lagrangian-equation}.

In dimension three, equation \eqref{eq:special-lagrangian-equation} with critical phase $\pi/2$ corresponds to
$\sigma_2(D^2u)=1$ with $\sigma_1(D^2u)>0$.
Warren~\cite{Warren} constructed the
nonquadratic entire solution
\[
u(x_1,x_2,x_3)=(x_1^2+x_2^2-1)e^{x_3}+\tfrac14e^{-x_3}
\]
to $\sigma_2(D^2u)=1$. Warren \cite{Warren}
also asked whether nonquadratic entire solutions exist at the critical phase
for $n\ge4$.

Together with the examples in dimensions two and three,
Theorem~\ref{thm:main} shows that the supercritical phase condition
$|\theta|>(n-2)\pi/2$ in \cite[Theorem~1.1]{YuanSupercritical}
is sharp for every $n\ge2$.

The philosophy is to control $Du$ inside a fixed ball by finding one
gradient level that occurs only outside the ball. We rescale $x$ and
$p$ separately and show that the rescaled $x$-coordinates stay away
from zero when the rescaled $|p|=1$. Since the $x$-scale tends to
infinity, we can choose $R_j$ so that $|p|=R_j$ on the limiting current
implies $|x|>R$. Since $p=Du(x)$ on the graph, we have $|Du|<R_j$ along
paths from $0$ in $B_R(0)$ by continuity from $Du(0)=0$. We thereby
keep the graph bounded near finite boundary points and reduce the
extension to the $x$-projection of the singular set.

The rest of the paper is organized as follows.
In Section~\ref{sec:ball}, we take a current limit $S$ of the rescaled
gradient graphs of the solutions to \eqref{eq:dirichlet}. We then show
that the connected component $\Lambda$ of
$\spt S\setminus\Sing S$ through the origin is a gradient graph
$\Graph(Du)$ over an unbounded domain.

Next, in Section~\ref{sec:cones}, we prove that the unique cone at infinity
of $S$ is the plane $\{x=0\}$ with multiplicity $2n-1$. Hence
$|x|/|p|\to0$ on $\spt S$ as $|p|\to\infty$.

In Section~\ref{sec:profile}, we complete the proof of Theorem \ref{thm:main} in three steps.
First, to study the $x$-coordinates of $\spt S$ at large $|p|$, we rescale
$S$ on a smaller scale in $x$ than in $p$ and construct a Dirichlet minimizer
$F$ with $0\notin\spt F(p)$ for
$p\ne0$, and show that $\dist(x,\spt F(p))$ tends uniformly to zero
on the rescaled support with $|p|=1$ as the $p$-scale tends to infinity.
Second, we use \eqref{eq:uniform-support-distance} and $0\notin\spt F(p)$
on $|p|=1$ to control $Du$ near $\partial(\pi_x(\Lambda))$ and extend $u$
smoothly to $\R^n$.
Third, we rule out polynomial $u$ by comparing $|x|/|Du(x)|$ along the
first coordinate axis with $|x|/|p|\to0$ on $\spt S$ as $|p|\to\infty$.


\section{A smooth limit on an unbounded domain}\label{sec:ball}
\label{sec:first}
\label{subsec:degree}
\label{subsec:area}
\label{subsec:inputs}
\label{subsec:boxes}

Fix $n\ge4$.
Let $e_i\coloneqq(0,\ldots,0,\underset{i\text{th}}{1},0,\ldots,0)$ be the $i$th standard
basis vector of $\R^n$, $1\le i\le n$.
Let $G$ be the group generated by coordinate permutations and
reflections across the coordinate hyperplanes in $\R^n$.
For $R>0$, let $B_R\subset\R^n$ and $\mathbf B_R\subset\C^n$ be
the balls centered at zero with radius $R$, and set
$\omega_n\coloneqq|B_1|$.

Define the $G$-invariant harmonic quartic by
\begin{equation}\label{eq:notation}
H_n(x)\coloneqq\sum_{i=1}^n x_i^4-\frac3{n+2}|x|^4.
\end{equation}
By \cite[Corollary 1.1]{BMS}, for each $t>0$ there
exists a unique $v^{(t)}\in C^3(B_1)\cap C^0(\overline B_1)$ satisfying
\begin{equation}\label{eq:dirichlet}
\begin{cases}
\Theta(D^2v^{(t)})=\dfrac{n-2}{2}\pi&\text{in }B_1,\\
v^{(t)}=\dfrac{\cot(\pi/n)}2-tH_n&\text{on }\partial B_1.
\end{cases}
\end{equation}

For each $t>0$, let
\begin{equation}\label{eq:first-radius}
r_t\coloneqq\min\Bigl\{\frac12,\,
\inf\bigl\{r\in(0,\tfrac12]:\osc_{B_r}v^{(t)}\ge\cot(\pi/n)r^2\bigr\}\Bigr\}.
\end{equation}
When $r_t>0$,   denote  
\begin{equation}\label{eq:rescaled-potentials}
u_t(y)\coloneqq r_t^{-2}\bigl(v^{(t)}(r_ty)-v^{(t)}(0)\bigr),
\qquad\forall y\in B_{1/r_t}.
\end{equation}

In this section, we prove the following result.
\begin{proposition}\label{prop:first-limit}
Let $r_t$ be as in \eqref{eq:first-radius}. Then
\begin{enumerate}[label=\textup{(\roman*)}]
\item $r_t>0$ for every $t>0$, and $r_t\to0$ as $t\to\infty$.
\item Let $u_{t_j}$ be as in \eqref{eq:rescaled-potentials} with  $t_j\to\infty$ as $j\to\infty$.
Then there exist  a calibrated current $S$ without boundary,
a connected component $\Lambda$ of $\spt S\setminus\Sing S$ containing $(0,0)$,
and $u\in C^\infty(\pi_x(\Lambda))$ such that, up to a subsequence, 
\[
\begin{aligned}
\Graph(Du_{t_j})&\rightharpoonup S&&\text{locally as currents},\\
u_{t_j}&\longrightarrow u&&\text{in }C^\infty_{\mathrm{loc}}(\pi_x(\Lambda))
\end{aligned}
\qquad\text{as }j\to\infty.
\]
Moreover, $\pi_x(\Lambda)$ is a connected unbounded open subset of $\R^n$
and $\Lambda=\Graph(Du)$.
\end{enumerate}
\end{proposition}

To prove Proposition~\ref{prop:first-limit}, we need
Lemmas~\ref{lem:critical-admissibility}--\ref{lem:bounded-source} below.
To state these results, for a real symmetric $n\times n$ matrix
$A\in\operatorname{Sym}(n)$, we define the $k$th elementary symmetric polynomial $\sigma_k$
and the  G{\aa}rding cone $\Gamma_k$ by
\[
\begin{gathered}
\sigma_0(A)\coloneqq1,\qquad
\sigma_k(A)\coloneqq\sum_{1\le i_1<\cdots<i_k\le n}
\lambda_{i_1}(A)\cdots\lambda_{i_k}(A),\\
\Gamma_k\coloneqq\{A\in\operatorname{Sym}(n):
\sigma_j(A)>0\text{ for all }1\le j\le k\},
\qquad 1\le k\le n.
\end{gathered}
\]
\begin{lemma}\label{lem:critical-admissibility}
Let $n\ge4$. Then
\[
\{A\in\operatorname{Sym}(n):\Theta(A)=\frac{n-2}{2}\pi\}
\subset\Gamma_{n-1}.
\]
\end{lemma}
\begin{proof}
For $A\in\operatorname{Sym}(n)$ satisfying $\Theta(A)=\frac{n-2}{2}\pi$,
let $\lambda_i\coloneqq\lambda_i(A)$ and
$\theta_i\coloneqq\arctan\lambda_i\in(-\pi/2,\pi/2)$,
with $\lambda_1\ge\cdots\ge\lambda_n$.
Subtracting the largest $n-2$ angles from the phase, we obtain
\begin{equation}\label{eq:two-smallest-angles}
\theta_{n-1}+\theta_n
=\frac{n-2}{2}\pi-\sum_{i=1}^{n-2}\theta_i>0.
\end{equation}
Hence $\theta_{n-1}>0$, so $\lambda_i\leq 0$ can hold only for $i=n$.
If $\lambda_n<0$, then
\[
0<\sum_{i=1}^{n-1}(\pi/2-\theta_i)
=\frac{n-1}{2}\pi-\sum_{i=1}^{n-1}\theta_i
=\frac\pi2+\theta_n<\frac\pi2.
\]
Note that, for $\theta',\theta''>0$ with $\theta'+\theta''<\pi/2$,
\begin{equation}\label{eq:tangent-addition-positive}
\tan(\theta'+\theta'')-\tan\theta'-\tan\theta''
=\tan\theta'\,\tan\theta''\,\tan(\theta'+\theta'')>0.
\end{equation}
Applying \eqref{eq:tangent-addition-positive} successively to the angles
$\pi/2-\theta_i$ and using $\tan\theta_n<0$, we obtain
\begin{equation}\label{eq:critical-reciprocal-sign}
\begin{aligned}
1+\lambda_n\sum_{i=1}^{n-1}\lambda_i^{-1}
&=1+\tan\theta_n\sum_{i=1}^{n-1}\cot\theta_i >1+\tan\theta_n\tan\Bigl(\sum_{i=1}^{n-1}(\pi/2-\theta_i)\Bigr)\\
&=1+\tan\theta_n\tan(\pi/2+\theta_n)
=1-\tan\theta_n\cot\theta_n=0.
\end{aligned}
\end{equation}
By \eqref{eq:critical-reciprocal-sign} when $\lambda_n<0$, and directly
when $\lambda_n\ge0$, we have
\[
\sigma_{n-1}(A)
=\Bigl(\prod_{i=1}^{n-1}\lambda_i\Bigr)
\Bigl(1+\lambda_n\sum_{i=1}^{n-1}\lambda_i^{-1}\Bigr)>0.
\]
Since $\lambda_1,\ldots,\lambda_{n-1}>0$,
\begin{equation}\label{eq:last-eigenvalue-shift}
\sigma_{n-1}(\lambda_1,\ldots,\lambda_{n-1},\lambda_n+\tau)
=\sigma_{n-1}(A)+\tau\sigma_{n-2}(\lambda_1,\ldots,\lambda_{n-1})>0
\qquad\forall\tau\ge0.
\end{equation}
By \eqref{eq:last-eigenvalue-shift},
$(\lambda_1,\ldots,\lambda_{n-1},\lambda_n+\tau)
\in\{\lambda\in\R^n:\sigma_{n-1}(\lambda)>0\}$ for $\tau\ge0$.
For $\tau>-\lambda_n$, $(\lambda_1,\ldots,\lambda_{n-1},\lambda_n+\tau)\in(0,\infty)^n$.
Thus varying $\lambda_n+\tau$ for $\tau\ge0$ joins the eigenvalues of $A$
to $(0,\infty)^n$ within $\{\sigma_{n-1}>0\}$,
so $A\in\Gamma_{n-1}$ \cite[p.~582, following (2.5)]{TW}.
\end{proof}

\begin{lemma}\label{lem:gauge}
Let $n\ge4$ and $A\in\Gamma_{n-1}$. Then there exists a unique
$s>0$ satisfying $\Theta(A/s)=\frac{n-2}{2}\pi$.

Moreover, the function $\mathcal F_n$ defined by $\mathcal F_n(A)\coloneqq s$ is smooth,
strictly elliptic, homogeneous of degree one and concave on
$\Gamma_{n-1}$, and satisfies
\begin{equation}\label{eq:gauge-properties}
\operatorname{tr}D_A\mathcal F_n(A)\ge\tan\frac\pi n
\qquad\forall A\in\Gamma_{n-1},
\end{equation}
and $\mathcal F_n$ can be   continuous extended to
$\overline{\Gamma}_{n-1}$ with $\mathcal F_n=0$ on $\partial\Gamma_{n-1}$.
\end{lemma}
\begin{proof}
Let $\lambda_i\coloneqq\lambda_i(A)$, with
$\lambda_1\ge\cdots\ge\lambda_n$.
If $\lambda_{n-1}=0$, then
\[
0<\sigma_{n-1}(A)
=\lambda_n\prod_{i=1}^{n-2}\lambda_i\le0,
\]
a contradiction. Since $\Gamma_{n-1}$ is connected
\cite[p.~582, (2.1) and (2.5)]{TW} and $\lambda_{n-1}(I)=1$,
we have $\lambda_{n-1}>0$.
If $\lambda_n<0$, then
\begin{equation}\label{eq:reciprocal-admissibility}
\sigma_{n-1}(A)
=\Bigl(\prod_{i=1}^{n-1}\lambda_i\Bigr)
\Bigl(1+\lambda_n\sum_{i=1}^{n-1}\lambda_i^{-1}\Bigr)>0.
\end{equation}
Therefore
\[
\lim_{s\downarrow0}
\frac{\Theta(A/s)-\frac{n-2}{2}\pi}{s}
=-\lambda_n^{-1}-\sum_{i=1}^{n-1}\lambda_i^{-1}>0.
\]
When $A\ge0$, we also have
$\lim_{s\downarrow0}\Theta(A/s)\ge(n-1)\pi/2>\frac{n-2}{2}\pi$.
In both cases, $\Theta(A/s)>(n-2)\pi/2$ for small $s>0$ and
$\Theta(A/s)\to0$ as $s\to\infty$. Thus there exists $s>0$ with
$\Theta(A/s)=(n-2)\pi/2$. At a positive solution $s$, the supporting inequality at
$(\lambda_i/s)_{i=1}^n$ for the convex set
$\{\mu\in\R^n:\sum_i\arctan\mu_i\ge(n-2)\pi/2\}$
\cite[Lemma~2.1]{YuanSupercritical} reads
\begin{equation}\label{eq:critical-scale-derivative}
0\le\sum_i\frac{(r-s^{-1})\lambda_i}{1+\lambda_i^2/s^2}
=-s^2(r-s^{-1})\partial_s\Theta(A/s)
\qquad\text{if }\Theta(rA)\ge(n-2)\pi/2.
\end{equation}
For large $r>s^{-1}$, $\Theta(rA)>(n-2)\pi/2$, so the inequality in
\eqref{eq:critical-scale-derivative} is strict and $\partial_s\Theta(A/s)<0$.
For $s'>s$, we have $-s^2(1/s'-1/s)\partial_s\Theta(A/s)<0$.
Thus \eqref{eq:critical-scale-derivative} with $r=1/s'$ implies
$\Theta(A/s')<(n-2)\pi/2$.
Thus $s=\mathcal F_n(A)$ is unique. As $\partial_s\Theta(A/s)<0$ at
$s=\mathcal F_n(A)$, $\mathcal F_n\in C^\infty(\Gamma_{n-1})$ by the
implicit function theorem.
Since, at $s=\mathcal F_n(A)$ and for
$\mathsf E\in\operatorname{Sym}(n)$,
\[
\begin{aligned}
0
 =\left.\frac{d}{dh}\Theta\Bigl(\frac{A+h\mathsf E}
 {\mathcal F_n(A+h\mathsf E)}\Bigr)\right|_{h=0} 
&=\left.\frac{d}{dh}\Theta\Bigl(\frac{A+h\mathsf E}{s}\Bigr)\right|_{h=0}
 +\partial_s\Theta(A/s)\,D_A\mathcal F_n(A)[\mathsf E]\\
&=s\tr\bigl((s^2I+A^2)^{-1}\mathsf E\bigr)
 +\partial_s\Theta(A/s)\,D_A\mathcal F_n(A)[\mathsf E].
\end{aligned}
\]
then 
\begin{equation}\label{eq:gauge-derivative}
D_A\mathcal F_n(A)
=\frac{s(s^2I+A^2)^{-1}}{-\partial_s\Theta(A/s)}.
\end{equation}
Since $-\partial_s\Theta(A/s)>0$ and $s(s^2I+A^2)^{-1}>0$,
\eqref{eq:gauge-derivative} implies $D_A\mathcal F_n(A)>0$.
For $t>0$,
\[
\Theta\Bigl(\frac{tA}{t\mathcal F_n(A)}\Bigr)
=\Theta\Bigl(\frac{A}{\mathcal F_n(A)}\Bigr)
=\frac{n-2}{2}\pi.
\]
By uniqueness of the positive solution of $\Theta(tA/s)=\frac{n-2}{2}\pi$,
\[
\mathcal F_n(tA)=t\mathcal F_n(A)\qquad\forall t>0.
\]

By continuity and uniqueness, $\Theta(A/s)>(n-2)\pi/2$ for
$0<s<\mathcal F_n(A)$. Since $\Theta(A/s)<(n-2)\pi/2$ for
$s>\mathcal F_n(A)$, we have
\[
\{\mathcal F_n\ge s\}
=\{A\in\Gamma_{n-1}:\Theta(A/s)\ge\frac{n-2}{2}\pi\}
\qquad\forall s>0.
\]
By \cite[Lemma 2.1]{YuanSupercritical}, $\{\mathcal F_n\ge s\}$ is convex.
Consequently, by homogeneity and convexity of $\{\mathcal F_n\ge1\}$, we obtain
for $A,B\in\Gamma_{n-1}$
\[
\frac{A+B}{\mathcal F_n(A)+\mathcal F_n(B)}
=\frac{\mathcal F_n(A)}{\mathcal F_n(A)+\mathcal F_n(B)}\frac A{\mathcal F_n(A)}
+\frac{\mathcal F_n(B)}{\mathcal F_n(A)+\mathcal F_n(B)}\frac B{\mathcal F_n(B)}
\in\{\mathcal F_n\ge1\}.
\]
Hence, again by homogeneity,
\[
\mathcal F_n(A+B)\ge\mathcal F_n(A)+\mathcal F_n(B)
\qquad\forall A,B\in\Gamma_{n-1}.
\]
Thus $\mathcal F_n$ is concave. Taking $B=\epsilon I$ with $\epsilon>0$,
we have $\mathcal F_n(A+\epsilon I)-\mathcal F_n(A)\ge\epsilon\mathcal F_n(I)$,
so
\[
\tr D_A\mathcal F_n(A)
=\lim_{\epsilon\downarrow0}
\frac{\mathcal F_n(A+\epsilon I)-\mathcal F_n(A)}{\epsilon}
\ge\mathcal F_n(I)=\tan\frac\pi n.
\]
Finally, let $A_j\in\Gamma_{n-1}$ converge to $A\in\partial\Gamma_{n-1}$.
Since $A_j\le\lambda_1(A_j)I$ and $\lambda_1(A_j)>0$, by ellipticity and
homogeneity, we obtain
\[
0<\mathcal F_n(A_j)
\le\mathcal F_n(\lambda_1(A_j)I)=\lambda_1(A_j)\tan\frac\pi n.
\]
If, up to a subsequence, $\mathcal F_n(A_j)\to s>0$, then
$\Theta(A/s)=\frac{n-2}{2}\pi$, so $A/s\in\Gamma_{n-1}$
by Lemma~\ref{lem:critical-admissibility}, a contradiction.
Therefore $\mathcal F_n(A_j)\to0$ whenever
$A_j\in\Gamma_{n-1}$ and $A_j\to A\in\partial\Gamma_{n-1}$.
Thus $\mathcal F_n\in C^0(\overline{\Gamma}_{n-1})$ with
$\mathcal F_n=0$ on $\partial\Gamma_{n-1}$.

\end{proof}

\begin{lemma}\label{lem:ball-mass}
Let $M_t\coloneqq\Graph(Dv^{(t)})$. Then the following hold.
\begin{enumerate}[label=\textup{(\roman*)}]
\item  
\begin{equation}\label{eq:ball-boundary}
\liminf_{t\to\infty}\frac{\dist(0,\spt\partial M_t)}t>0.
\end{equation}
\item There exists $\delta>0$ such that
\begin{equation}\label{eq:ball-area}
\limsup_{t\to\infty}\sup_{0<R<\delta t/2}
\frac{\|M_t\|(\mathbf B_R)}{\omega_nR^n}\le2n-1.
\end{equation}
\end{enumerate}
\end{lemma}
\begin{proof}
\textup{(i)} Set $U_t\coloneqq t^{-1}v^{(t)}$. 
Applying Lemma~\ref{lem:critical-admissibility} to \eqref{eq:dirichlet},
we have $ D^2U_t=t^{-1}D^2v^{(t)}\in\Gamma_{n-1}$.  
Using Lemma~\ref{lem:gauge} for the equation and
\eqref{eq:dirichlet} for the boundary value, we obtain
\begin{equation}\label{eq:quotient}
\begin{cases}
\mathcal F_n(D^2U_t)
=t^{-1}\mathcal F_n(D^2v^{(t)})=t^{-1},\quad D^2U_t\in\Gamma_{n-1}
&\text{in }B_1,\\
U_t-\dfrac{\cot(\pi/n)}{2t}
=t^{-1}\Bigl(v^{(t)}-\dfrac{\cot(\pi/n)}2\Bigr)=-H_n
&\text{on }\partial B_1.
\end{cases}
\end{equation}
Choose $C(n)\ge\max\{\|D^2H_n\|_{C^0(\overline B_1)},\cot(\pi/n)\}$.
Then $-D^2H_n+2C(n)I\ge C(n)I$, so ellipticity of $\mathcal F_n$ implies
\begin{equation}\label{eq:dirichlet-subsolution}
\begin{aligned}
\mathcal F_n\bigl(D^2[-H_n+C(n)(|x|^2-1)]\bigr)
&=\mathcal F_n(-D^2H_n+2C(n)I)\\
&\ge\mathcal F_n(C(n)I)=C(n)\tan(\pi/n)\ge1\ge t^{-1}
\qquad\forall t\ge1.
\end{aligned}
\end{equation}
Thus, for $t\ge1$, there exists
$W_t\in C^\infty(B_1)\cap C^2(\overline B_1)$ satisfying
\[
\begin{cases}
\mathcal F_n(D^2W_t)=t^{-1},\quad D^2W_t\in\Gamma_{n-1}
&\text{in }B_1,\\
W_t=-H_n&\text{on }\partial B_1,
\end{cases}
\]
by the continuity method with the subsolution
$-H_n+C(n)(|x|^2-1)$ from \eqref{eq:dirichlet-subsolution},
which equals $-H_n$ on $\partial B_1$, and the a priori estimates
of \cite[Theorem 1.1]{ITW} for $\mathcal F_n$ from Lemma~\ref{lem:gauge}.
The comparison principle for \eqref{eq:quotient} identifies $W_t$ with
$U_t-\cot(\pi/n)/(2t)$. Hence
\[
U_t=W_t+\frac{\cot(\pi/n)}{2t}\in C^2(\overline B_1).
\]

Applying the global Hessian estimate of \cite[Theorem 1.1]{ITW} for the concave
operator $\mathcal F_n$ in Lemma~\ref{lem:gauge}, with the subsolution in
\eqref{eq:dirichlet-subsolution}, we have
\[
\|D^2U_t\|_{C^0(\overline B_1)}
\le C\bigl(n,\|t^{-1}\|_{C^2(\overline B_1)}\bigr)\le C(n)
\qquad\forall t\ge1.
\]
By uniqueness in \eqref{eq:dirichlet}, $U_t$ is even in each coordinate,
so $DU_t(0)=0$. Then
\[
|DU_t(x)|=\Bigl|\int_0^1D^2U_t(sx)x\,ds\Bigr|\le C(n)
\qquad\forall x\in\overline B_1,\quad t\ge1.
\]
Integrating $DU_t$ from $e_1$ to $x$ and using the boundary value in
\eqref{eq:quotient}, we obtain
\begin{equation}\label{eq:C11}
\|U_t\|_{C^{1,1}(\overline B_1)}
\le C(n)\bigl(|U_t(e_1)|+\|D^2U_t\|_{C^0(\overline B_1)}\bigr)
\le C(n)\qquad\forall t\ge1.
\end{equation}

Next, for $0<h<1$, let
\[
\phi_h(x)\coloneqq
(1-x_1+h)^{\frac{n-3}{n-2}}
+(n-1)(|x|^2-x_1^2)(1-x_1+h)^{-\frac{n-1}{n-2}}
\qquad\forall x\in\overline B_1.
\]
Note that $\xi_1=0$ for $\xi\in e_1^\perp$. Thus,
\[
\langle D^2\phi_h\xi,\xi\rangle
=2(n-1)(1-x_1+h)^{-\frac{n-1}{n-2}}|\xi|^2>0
\qquad\forall\xi\in e_1^\perp\setminus\{0\}.
\]
As $e_1^\perp$ has dimension $n-1$, $\lambda_{n-1}(D^2\phi_h)>0$.
Moreover, summing the principal minors of order $n-1$, we obtain
\begin{equation}\label{eq:power-barrier}
\sigma_{n-1}(D^2\phi_h)
=\frac{2^{n-2}(n-1)^{n-1}}
{(n-2)^2(1-x_1+h)^{(n-1)^2/(n-2)}}
\Bigl[2(n-2)^2-n+3+
\frac{(n-1)^2(|x|^2-x_1^2)}{(1-x_1+h)^2}\Bigr]>0,
\end{equation}
since $n\ge4$ and $|x|^2-x_1^2\ge0$.
Since $\lambda_{n-1}(D^2\phi_h)>0$,
\eqref{eq:power-barrier} and \eqref{eq:last-eigenvalue-shift} imply
that increasing the smallest eigenvalue connects $D^2\phi_h$
to a positive definite matrix within $\{\sigma_{n-1}>0\}$.
Hence $D^2\phi_h\in\Gamma_{n-1}$.

Let
\[
w_h(x)\coloneqq\phi_h(x)-\phi_h(e_1)-(n+2)(1-x_1).
\]
Then $w_h(e_1)=0$ and $D^2w_h=D^2\phi_h\in\Gamma_{n-1}$ on $\overline B_1$.
For $x\in\mathbb S^{n-1}$ near $e_1$, we expand $w_h$ at $e_1$ using
$Dw_h(e_1)\parallel e_1$ and obtain
\[
\begin{aligned}
w_h(x)&\le\partial_1w_h(e_1)(x_1-1)
  +\tfrac12\|D^2w_h\|_{C^0(\overline B_1)}|x-e_1|^2\\
&\le\bigl(|\partial_1w_h(e_1)|+\|D^2w_h\|_{C^0(\overline B_1)}\bigr)
\Bigl(1-\sum_i x_i^4\Bigr),
\end{aligned}
\]
where $|x-e_1|^2=2(1-x_1)$ and $\sum_i x_i^4\le x_1^2\sum_i x_i^2\le x_1$
since $x_1=\max_i|x_i|$ on $\mathbb S^{n-1}$ near $e_1$. On the other hand,
\[
\phi_h(-e_1)<3,\qquad \phi_h(\pm e_i)\le n+1\quad\forall i=2,\ldots,n.
\]
Since $1-x_1\ge1$ at $x\in\{-e_1,\pm e_2,\ldots,\pm e_n\}$, the estimates
for $\phi_h(x)$ imply $w_h(x)\le n+1-(n+2)<0$.
Thus $w_h(x)/(1-\sum_i x_i^4)\leq 0$
near $-e_1,\pm e_2,\ldots,\pm e_n$ whenever $\sum_i x_i^4<1$.
The quotient is continuous on the compact complement in $\mathbb S^{n-1}$
of the neighborhoods of the coordinate poles. By the Taylor estimate at
$e_1$ and the sign near the other poles, we can choose $C(n,h)>0$ with
\[
C(n,h)\ge\sup_{x\in\mathbb S^{n-1},\ \sum_i x_i^4<1}
\frac{w_h(x)}{1-\sum_i x_i^4}.
\]
Since $w_h(\pm e_i)\le0$ for $1\le i\le n$, we obtain from \eqref{eq:notation}
\begin{equation}\label{eq:pole-boundary-comparison}
-\frac{n-1}{n+2}+\frac{w_h(x)}{C(n,h)}
\le-\frac{n-1}{n+2}+1-\sum_i x_i^4=-H_n(x)
\qquad\forall x\in\mathbb S^{n-1}.
\end{equation}
Fix $0<h<1$ with $\frac{n-3}{n-2}h^{-\frac1{n-2}}>n+2$.
For $t^{-1}\le\min_{\overline B_1}\mathcal F_n(D^2w_h/C(n,h))$, the
comparison principle for \eqref{eq:quotient} and
\eqref{eq:pole-boundary-comparison} yield
\begin{equation}\label{eq:pole-barrier-comparison}
-\frac{n-1}{n+2}+\frac{w_h(x)}{C(n,h)}
\le U_t(x)-\frac{\cot(\pi/n)}{2t}
\qquad\forall x\in\overline B_1.
\end{equation}
Let $\nu(x)\coloneqq x$ be the outward unit normal to $\partial B_1$.
Equality holds in \eqref{eq:pole-barrier-comparison} at $e_1$, so
\begin{equation}\label{eq:axis-normal-sign}
\begin{aligned}
\partial_\nu U_t(e_1)
&\le\lim_{\epsilon\downarrow0}
\frac{w_h(e_1)-w_h((1-\epsilon)e_1)}{C(n,h)\epsilon}\\
&=\frac{\partial_\nu w_h(e_1)}{C(n,h)}=\frac{1}{C(n,h)}
\Bigl(n+2-\frac{n-3}{n-2}h^{-\frac1{n-2}}\Bigr)<0.
\end{aligned}
\end{equation}
By $G$-invariance of $U_t$,
$\partial_\nu U_t(\pm e_i)=\partial_\nu U_t(e_1)$, $\forall i=1,\ldots,n$.

On the other hand, with $A=D^2v^{(t)}$ in \eqref{eq:two-smallest-angles}, we have
$\lambda_{n-1}(D^2v^{(t)})>-\lambda_n(D^2v^{(t)})$.
Thus, for distinct $i,j\in\{1,\ldots,n\}$,
\begin{equation}\label{eq:least-pair}
\Delta_{ij}v^{(t)}
=\partial_{ii}v^{(t)}+\partial_{jj}v^{(t)}
\ge\lambda_{n-1}(D^2v^{(t)})+\lambda_n(D^2v^{(t)})>0
\quad\text{in }B_1.
\end{equation}
Let $x\in\partial B_1$ with $R\coloneqq(x_i^2+x_j^2)^{1/2}>0$.
Fix $y_k=x_k$ for $k\notin\{i,j\}$. Then
\[
\Delta_{ij}\Bigl[y_i^4+y_j^4-\frac34(y_i^2+y_j^2)^2\Bigr]
=\Bigl(12-\frac34\cdot16\Bigr)(y_i^2+y_j^2)=0.
\]
Consequently, by \eqref{eq:least-pair}, the function
\[
v^{(t)}(y)+t\Bigl[y_i^4+y_j^4-\tfrac34(y_i^2+y_j^2)^2\Bigr]
\]
is subharmonic in the disc $y_i^2+y_j^2<R^2$. On the circle
$y_i^2+y_j^2=R^2$, we have $|y|^2=R^2+\sum_{k\notin\{i,j\}}x_k^2=1$.
Substituting the boundary value in \eqref{eq:dirichlet}, we obtain
\begin{equation}\label{eq:disc-boundary-value}
v^{(t)}(y)+t\Bigl[y_i^4+y_j^4-\tfrac34(y_i^2+y_j^2)^2\Bigr]
=\frac{\cot(\pi/n)}2
-t\Bigl[\sum_{k\notin\{i,j\}}x_k^4+\tfrac34R^4-\tfrac3{n+2}\Bigr].
\end{equation}
Taking tangential derivatives in \eqref{eq:dirichlet}, we obtain
\begin{equation}\label{eq:boundary-coordinate-derivative}
\partial_i v^{(t)}=x_i\Bigl(\partial_\nu v^{(t)}+4t\sum_kx_k^4-4tx_i^2\Bigr)
\qquad\forall x\in\partial B_1.
\end{equation}
By the maximum principle, the left-hand side of \eqref{eq:disc-boundary-value}
is maximal at $y=x$ on the closed disc $y_i^2+y_j^2\le R^2$ with
$y_k=x_k$ for $k\notin\{i,j\}$. Multiplying the outward radial derivative
of the left-hand side at $y=x$ by $R$ and substituting
\eqref{eq:boundary-coordinate-derivative}, we obtain
\begin{equation}\label{eq:disc-normal-derivative}
\begin{aligned}
0&\le
\Bigl.(x_i\partial_{y_i}+x_j\partial_{y_j})
\Bigl[v^{(t)}(y)+t\bigl(y_i^4+y_j^4-\tfrac34(y_i^2+y_j^2)^2\bigr)\Bigr]
\Bigr|_{y=x}\\
&=R^2\Bigl(\partial_\nu v^{(t)}+4t\sum_kx_k^4-3tR^2\Bigr).
\end{aligned}
\end{equation}
Combining \eqref{eq:boundary-coordinate-derivative} and
\eqref{eq:disc-normal-derivative}, for $x\in\partial B_1$
with $0<|x_i|\le|x_j|$, we have 
\begin{equation}\label{eq:spherical-coordinate-sign}
\begin{aligned}
x_i\partial_iU_t
&=x_i^2\Bigl(t^{-1}\partial_\nu v^{(t)}+4\sum_kx_k^4-4x_i^2\Bigr) \ge x_i^2(3x_j^2-x_i^2)\ge2x_i^4>0.
\end{aligned}
\end{equation}
If $\liminf_{t\to\infty}\min_{\partial B_1}|DU_t|=0$, compactness of
$\mathbb S^{n-1}$ allows us to take $t_j\to\infty$ and $x_j\to x_*$ on
$\mathbb S^{n-1}$ with $|DU_{t_j}(x_j)|\to0$ as $j\to\infty$.
Then \eqref{eq:C11} implies
\[
|DU_{t_j}(x_*)|
\le|DU_{t_j}(x_j)|+C(n)|x_j-x_*|\longrightarrow0
\quad\text{as }j\to\infty.
\]
If $x_*=\pm e_i$, then by \eqref{eq:axis-normal-sign},
$\liminf_j|DU_{t_j}(x_*)|\ge-\partial_\nu w_h(e_1)/C(n,h)>0$,
a contradiction. Otherwise, choose $i\ne\ell$ with
$0<|(x_*)_i|\le|(x_*)_\ell|$;
then by \eqref{eq:spherical-coordinate-sign}, 
$|DU_{t_j}(x_*)|\ge2|(x_*)_i|^3>0$, again a contradiction.
Consequently,
\begin{equation}\label{eq:boundary-gradient-nonvanishing}
\liminf_{t\to\infty}\min_{\partial B_1}|DU_t|>0.
\end{equation}
Hence
\[
\liminf_{t\to\infty}\frac{\dist(0,\spt\partial M_t)}t
\ge\liminf_{t\to\infty}\frac{\min_{\partial B_1}|Dv^{(t)}|}t
=\liminf_{t\to\infty}\min_{\partial B_1}|DU_t|>0.
\]
Thus \textup{(i)} holds.

\textup{(ii)} For each $x\in \mathbb S^{n-1}$, choose $i$ with
$0<|x_i|=\min_{x_k\ne0}|x_k|$.
If $x_j=0$ for all $j\ne i$, then $x=\pm e_i$.
Otherwise, choose $j\ne i$ with $x_j\ne0$; then
$0<|x_i|\le|x_j|$ and $2x_i^2\le x_i^2+x_j^2\le1$.
Using \eqref{eq:axis-normal-sign} for $x=\pm e_i$ and
\eqref{eq:spherical-coordinate-sign} for $0<|x_i|\le|x_j|$ with $j\ne i$,
we have for $0\le\lambda\le1$
\[
\begin{aligned}
x_i\bigl[(1-\lambda)\partial_iU_t(x)+\lambda(3x_i-4x_i^3)\bigr]
&=(1-\lambda)x_i\partial_iU_t(x)
 +\lambda x_i^2(3-4x_i^2)\ne0.
\end{aligned}
\]
Then
\begin{equation}\label{eq:boundary-homotopy}
(1-\lambda)DU_t(x)+\lambda(3x_k-4x_k^3)_{k=1}^n\ne0
\qquad\forall x\in \mathbb S^{n-1},\quad0\le\lambda\le1.
\end{equation}
Choose $0<\delta<\liminf_{t\to\infty}\min_{\partial B_1}|DU_t|$ by
\eqref{eq:boundary-gradient-nonvanishing}.
For $t$ with $\min_{\partial B_1}|DU_t|>\delta$ and $p\in B_\delta$, we have
\[
|\lambda p|\le |p|<\delta<\min_{x\in\partial B_1}|DU_t(x)|
\qquad\forall\lambda\in[0,1].
\]
By homotopy invariance along the targets $\lambda p$ and along
\eqref{eq:boundary-homotopy}, we obtain
\begin{equation}\label{eq:degree}
\begin{aligned}
\degmap(DU_t,B_1,p)
&=\degmap(DU_t,B_1,0) =\degmap\bigl((3x_k-4x_k^3)_{k=1}^n,B_1,0\bigr)=1-2n
\qquad\forall p\in B_\delta,
\end{aligned}
\end{equation}

By \eqref{eq:dirichlet}, $M_t$ is calibrated by
$\operatorname{Re}(e^{-i(n-2)\pi/2}dz_1\wedge\cdots\wedge dz_n)$.
Since $e^{-i(n-2)\pi/2}i^n=-1$, we obtain
\begin{equation}\label{eq:critical-density}
\begin{aligned}
t^{-n}\sqrt{\det(I+(D^2v^{(t)})^2)}
&=t^{-n}\operatorname{Re}\bigl(e^{-i(n-2)\pi/2}\det(I+itD^2U_t)\bigr)\\
&=-\operatorname{Re}\det(D^2U_t-it^{-1}I)\\
&=-\det D^2U_t+\sum_{k=1}^{\lfloor n/2\rfloor}(-1)^{k+1}t^{-2k}\sigma_{n-2k}(D^2U_t)\\
&=-\det D^2U_t+O(t^{-2}),
\end{aligned}
\end{equation}
The $O(t^{-2})$ in \eqref{eq:critical-density} is uniform on $B_1$ by \eqref{eq:C11}.
If $|(x,Dv^{(t)}(x))|<\delta t/2$, then $|DU_t(x)|<\delta/2$.
Using the graph area formula, \eqref{eq:critical-density}, the signed area
formula for $DU_t$, and the degree in \eqref{eq:degree}, we obtain
\begin{equation}\label{eq:ball-mass-estimate}
\begin{aligned}
t^{-n}\|M_t\|(\mathbf B_{\delta t/2})
&\le\int_{DU_t^{-1}(B_{\delta/2})}
  \sqrt{\det((D^2U_t)^2+t^{-2}I)}\,dx\\
&=-\int_{DU_t^{-1}(B_{\delta/2})}\det D^2U_t\,dx+O(t^{-2})\\
&=-\int_{B_{\delta/2}}\degmap(DU_t,B_1,p)\,dp+O(t^{-2})\\
&\le(2n-1)\omega_n(\delta/2)^n+C(n)t^{-2},
\end{aligned}
\end{equation}
Since $\min_{\partial B_1}|Dv^{(t)}|=t\min_{\partial B_1}|DU_t|>\delta t$,
$M_t$ has no boundary in $\mathbf B_{\delta t/2}$.
Monotonicity and \eqref{eq:ball-mass-estimate} imply
\[
\limsup_{t\to\infty}\sup_{0<R<\delta t/2}
\frac{\|M_t\|(\mathbf B_R)}{\omega_nR^n}
\le \limsup_{t\to\infty}\frac{t^{-n}\|M_t\|(\mathbf B_{\delta t/2})}
{\omega_n(\delta/2)^n}\le2n-1.
\]
Thus \textup{(ii)} holds.
\end{proof}

\begin{lemma}\label{lem:interior-estimate}
Let $\mathcal U\subset\R^n$ be a domain and let $w_j\in C^\infty(\mathcal U)$
solve $\Theta(D^2w_j)=\frac{n-2}{2}\pi$ in $\mathcal U$.
Assume that
\[
\sup_{j\ge1}\sup_{x\in K}w_j(x)<\infty
\qquad\forall K\Subset\mathcal U
\]
and that $\inf_{j\ge1}w_j(x_0)>-\infty$ for some $x_0\in\mathcal U$.
Then $\{w_j\}$ is precompact in $C^\infty_{\mathrm{loc}}(\mathcal U)$.
\end{lemma}
\begin{proof}
By Lemma~\ref{lem:critical-admissibility},
$D^2w_j\in\Gamma_{n-1}$, so $\Delta w_j>0$ in $\mathcal U$.
On $B_r(x_0)\Subset\mathcal U$, the mean-value inequality implies
$\int_{B_r(x_0)}w_j\ge\omega_nr^n w_j(x_0)$. Hence
\[
\begin{aligned}
\int_{B_r(x_0)}|w_j|
 &=2\int_{B_r(x_0)}\max\{w_j,0\}-\int_{B_r(x_0)}w_j 
\le 2\omega_nr^n\sup_{B_r(x_0)}\max\{w_j,0\}-\int_{B_r(x_0)}w_j\\
&\le \omega_nr^n\Bigl(2\sup_{B_r(x_0)}\max\{w_j,0\}-w_j(x_0)\Bigr).
\end{aligned}
\]
Thus $\sup_j\int_{B_r(x_0)}|w_j|<\infty$.
Since $D^2w_j\in\Gamma_{n-1}$ in $\mathcal U$ and $n-1>n/2$, by \cite[Theorem 2.7]{TW},
\[
\sup_j\|w_j\|_{C^{0,(n-2)/(n-1)}(B_{r/2}(x_0))}
\le C(n,r)\sup_j\|w_j\|_{L^1(B_r(x_0))}<\infty.
\]
Because $\mathcal U$ is connected, we further obtain 
$\sup_j\|w_j\|_{C^{0,(n-2)/(n-1)}(K')}<\infty$, $\forall K'\Subset\mathcal U$.
Therefore, for $K\Subset K'\Subset\mathcal U$, by \cite[Theorems 1.1 and 1.2]{BMS},
\[
\sup_j\|D^2w_j\|_{C^0(K)}<\infty.
\]
Note that, by Lemma~\ref{lem:critical-admissibility} and continuity of $\Theta$,
\[
\overline{\bigcup_{j\ge1}D^2w_j(K)}\Subset\Gamma_{n-1}.
\]
Thus, by Lemma~\ref{lem:gauge} and \eqref{eq:gauge-derivative},
$w_j$ satisfy $\mathcal F_n(D^2w_j)=1$ with ellipticity constants
independent of $j$ on $K$.
Since $\mathcal F_n$ is concave by Lemma~\ref{lem:gauge},
applying  Evans--Krylov estimate~\cite[Theorem 17.14]{GT}
and the Schauder estimates~\cite[Theorem 6.2]{GT}, we obtain
\[
\sup_{j\ge1}\|w_j\|_{C^k(K)}<\infty
\qquad\forall K\Subset\mathcal U,\quad\forall k=0,1,\ldots.
\]
Hence $\{w_j\}$ is precompact in $C^\infty_{\mathrm{loc}}(\mathcal U)$
by the Arzel\`a--Ascoli theorem.
\end{proof}

For $x,p\in\R^n$, let $\pi_x(x,p)\coloneqq x$
and $\pi_p(x,p)\coloneqq p$.
\begin{lemma}\label{lem:bounded-source}
Let $S\ne0$ be an $n$-dimensional locally integral current in $\C^n$
with $\partial S=0$, calibrated by
$\operatorname{Re}(e^{-i(n-2)\pi/2}dz_1\wedge\cdots\wedge dz_n)$, and satisfying
\[
(O\times O)_\sharp S=(\det O)S\qquad\forall O\in G.
\]
Assume that
\begin{equation}\label{eq:bounded-source-mass-ratio}
\sup_{R>0}\frac{\|S\|(\mathbf B_R)}{\omega_nR^n}<2n,
\end{equation}
and that there exist $c>0$ and ${M_*}<\infty$ such that
\begin{equation}\label{eq:center-cylinder}
\spt S\cap\{|x|<c\}\subset\{|p|\le {M_*}\}.
\end{equation}
Then $\pi_x(\spt S)$ is unbounded.
\end{lemma}
\begin{proof}
Assume not. Suppose $\pi_x(\spt S)$ is bounded.
Then $\pi_p:\spt S\to\R^n$ is proper.

After reversing the orientation of $S$, let $P$ and $X$ be the
matrices of $d\pi_p$ and $d\pi_x$ in an oriented orthonormal tangent
basis at $\|S\|$-a.e. $z\in\spt S$. Then
$J_p\coloneqq\det P$, $|X|^2=|\nabla_S x|^2$, and
$\operatorname{Re}\det(P+iX)=1$ by the calibration in the coordinates
$p+ix$. Since $\det(P-isX)=\overline{\det(P+isX)}$, we have
$\left.\frac{d}{ds}\operatorname{Re}\det(P+isX)\right|_{s=0}=0$.
Using $J_p\le1$ and Taylor's formula, we obtain
\begin{equation}\label{eq:calibrated-height-jacobian}
\begin{aligned}
0\le1-J_p&=\operatorname{Re}\det(P+iX)-\det P\\
&=\int_0^1(1-s)\frac{d^2}{ds^2}\operatorname{Re}\det(P+isX)\,ds \le C(n)|X|^2=C(n)|\nabla_S x|^2,
\end{aligned}
\end{equation}
where each determinant in $\frac{d^2}{ds^2}\det(P+isX)$ contains two
columns of $iX$ and $n-2$ columns of $P+isX$, which are uniformly bounded
for $0\le s\le1$.
Since $\pi_p$ is proper on $\spt S$ and $\partial S=0$,
$\partial((\pi_p)_\sharp S)=(\pi_p)_\sharp(\partial S)=0$.
By \cite[Theorem 26.27]{Simon},
\begin{equation}\label{eq:bounded-projection-current}
(\pi_p)_\sharp S=Q_0\llbracket\R^n\rrbracket,
\qquad Q_0\in\mathbb Z.
\end{equation}
For $R>0$, let ${\chi_R}={\chi_R}(p)\in C_c^\infty(B_{2R})$ satisfy
$0\le{\chi_R}\le1$, ${\chi_R}|_{B_R}=1$ and $|D{\chi_R}|\le C(n)/R$.
Because $S$ is calibrated and $\partial S=0$,
\begin{equation}\label{eq:first-variation}
\int\operatorname{div}_S X\,d\|S\|=0
\qquad\forall X\in C_c^\infty(\R^{2n};\R^{2n}).
\end{equation}
Since $\spt S\cap\pi_p^{-1}(\spt\chi_R)$ is compact, there exists
$X\in C_c^\infty(\R^{2n};\R^{2n})$ such that
$X(x,p)=(\chi_R(p)^2x,0)$ near $\spt S$.
By \eqref{eq:first-variation} and Young's inequality,
\[
\begin{aligned}
\int{\chi_R}^2|\nabla_Sx|^2\,d\|S\|
&=-2\int{\chi_R}\sum_{i=1}^n
x_i\nabla_Sx_i\cdot\nabla_S{\chi_R}\,d\|S\|\\
&\le\frac12\int{\chi_R}^2|\nabla_Sx|^2\,d\|S\|
+2\sup_{\pi_p^{-1}(B_{2R})\cap\spt S}|x|^2
\int|\nabla_S{\chi_R}|^2\,d\|S\|.
\end{aligned}
\]
Hence, using ${\chi_R}|_{B_R}=1$ and $|D{\chi_R}|\le C(n)/R$,
\begin{equation}\label{eq:bounded-height-energy}
\begin{aligned}
\int_{\pi_p^{-1}(B_R)}|\nabla_Sx|^2\,d\|S\|
&\le\int{\chi_R}^2|\nabla_Sx|^2\,d\|S\| \le4\sup_{\pi_p^{-1}(B_{2R})\cap\spt S}|x|^2
\int|\nabla_S{\chi_R}|^2\,d\|S\|\\
&\le C(n)R^{-2}\sup_{\pi_p^{-1}(B_{2R})\cap\spt S}|x|^2
\|S\|(\pi_p^{-1}(B_{2R}))\qquad\forall R>0.
\end{aligned}
\end{equation}
Let $A_0\coloneqq\sup_{\spt S}|x|<\infty$, then
$\pi_p^{-1}(B_{2R})\cap\spt S\subset\mathbf B_{3R}$ whenever $R>0$ and $R\ge A_0$.
Thus, by \eqref{eq:bounded-projection-current},
\eqref{eq:calibrated-height-jacobian} and \eqref{eq:bounded-height-energy},
\begin{equation}\label{eq:bounded-source-excess}
\begin{aligned}
0\le\|S\|(\pi_p^{-1}(B_R))-{Q_0}\omega_nR^n
&=\int_{\pi_p^{-1}(B_R)}(1-J_p)\,d\|S\| 
 \le C(n)\int_{\pi_p^{-1}(B_R)}|\nabla_Sx|^2\,d\|S\|\\
&\le C(n)R^{-2}\sup_{\pi_p^{-1}(B_{2R})\cap\spt S}|x|^2
\|S\|(\pi_p^{-1}(B_{2R}))\\
&\le C(n)A_0^2R^{-2}\|S\|(\mathbf B_{3R})
\le C(n)A_0^2R^{n-2}.
\end{aligned}
\end{equation}
Therefore,
\[
0\le\frac{\|S\|(\pi_p^{-1}(B_R))}{\omega_nR^n}-Q_0
\le\frac{C(n)A_0^2}{\omega_nR^2}
\longrightarrow0\quad\text{as }R\to\infty.
\]
Since $S\ne0$, there exists $z_0\in\spt S$ at which $\|S\|$ has
$n$-dimensional density at least one. Since
\[
\mathbf B_{R-|z_0|}(z_0)\subset\mathbf B_R
\subset\pi_p^{-1}(B_R)\qquad\forall R>|z_0|,
\]
monotonicity of the mass ratio of $S$ at $z_0$ yields
\[
\begin{aligned}
\Bigl(1-\frac{|z_0|}{R}\Bigr)^n
&\le\frac{\|S\|(\mathbf B_{R-|z_0|}(z_0))}{\omega_nR^n}
\le\frac{\|S\|(\pi_p^{-1}(B_R))}{\omega_nR^n} \\&\le\Bigl(1+\frac{A_0^2}{R^2}\Bigr)^{n/2}
\sup_{\rho>0}\frac{\|S\|(\mathbf B_\rho)}{\omega_n\rho^n}
\qquad\forall R>|z_0|.
\end{aligned}
\]
Letting $R\to\infty$ and using \eqref{eq:bounded-source-mass-ratio}, we have
\[
1\le Q_0
\le\sup_{\rho>0}\frac{\|S\|(\mathbf B_\rho)}{\omega_n\rho^n}<2n.
\]
Hence $Q_0\in\{1,\ldots,2n-1\}$.
Since $1-J_p\ge0$, we obtain from \eqref{eq:bounded-source-excess},
for $j\ge1$, $R>0$ and $4jR\ge A_0$,
\begin{equation}\label{eq:bounded-approximation-mass}
\begin{aligned}
0&\le\frac{\|(R^{-1}\operatorname{Id})_\sharp S\|(\pi_p^{-1}(B_{4j}))}
{\omega_n(4j)^n}-Q_0\\
&=\frac{\|S\|(\pi_p^{-1}(B_{4jR}))-Q_0\omega_n(4jR)^n}
{\omega_n(4jR)^n} =\frac{\displaystyle\int_{\pi_p^{-1}(B_{4jR})}(1-J_p)\,d\|S\|}
{\omega_n(4jR)^n}
\le\frac{C(n)A_0^2}{16\omega_nj^2R^2}.
\end{aligned}
\end{equation}
For integers $N\ge1$, we define
\[
\A_N(\R^n)\coloneqq
\Bigl\{\sum_{\ell=1}^N\llbracket\mathbf a_\ell\rrbracket:\mathbf a_\ell\in\R^n\Bigr\},
\]
where $\llbracket \mathbf a\rrbracket$ denotes unit mass at $\mathbf a\in\R^n$.
For $\mathbf a_\ell,\mathbf b_\ell\in\R^n$, set
\[
\mathsf d\Bigl(\sum_{\ell=1}^N\llbracket\mathbf a_\ell\rrbracket,
\sum_{\ell=1}^N\llbracket\mathbf b_\ell\rrbracket\Bigr)
\coloneqq\min_\sigma
\Bigl(\sum_{\ell=1}^N|\mathbf a_\ell-\mathbf b_{\sigma(\ell)}|^2\Bigr)^{1/2},
\]
where $\sigma$ ranges over all permutations of $\{1,\ldots,N\}$
\cite[Definitions 0.1 and 0.2]{DSQ}.
By \eqref{eq:bounded-approximation-mass}, choose an increasing
sequence $\{\widehat R_j\}$ with $\widehat R_j>\max\{j^n,A_0/(4j)\}$
so that \cite[Theorems 1.4 and 1.6]{DS} apply to
$(\widehat R_j^{-1}\operatorname{Id})_\sharp S$ on $\pi_p^{-1}(B_{4j})$
with $\bar\eta=j^{-(n+2)}$. By (1.5), the Lipschitz approximation
$\widehat R_j^{-1}\widehat F_j$ agrees with the slices over a closed
set $\widehat K_j\subset B_j$ satisfying
\begin{equation}\label{eq:bounded-slice-complement}
|B_j\setminus\widehat K_j|
\le C(n)j^n\frac{A_0^2}{j^2\widehat R_j^2}
\le C(n)A_0^2j^{-n-2}\longrightarrow0\quad\text{as }j\to\infty.
\end{equation}
For a.e. $p\in\widehat K_j$, the slice of $S$ over $\widehat R_jp$ is
\begin{equation}\label{eq:bounded-slice-approximation}
\sum_{\mathbf a\in \widehat F_j(p)}\llbracket(\mathbf a,\widehat R_jp)\rrbracket
=\langle S,\pi_p,\widehat R_jp\rangle.
\end{equation}
By $|x|\le A_0$ on $\spt S$, (1.7), and \eqref{eq:bounded-slice-approximation},
$\max_{\mathbf a\in\widehat F_j(p)}|\mathbf a|\le C(n)A_0$ for all $p\in B_j$.
Multiplying (1.8) for $\widehat R_j^{-1}\widehat W_j$ by $\widehat R_j^2$
and using \eqref{eq:bounded-approximation-mass}, we obtain
\begin{equation}\label{eq:bounded-dirichlet-approximation}
\int_{B_j}\mathsf d(\widehat F_j,\widehat W_j)^2
\le C(n)A_0^2j^{-2}\longrightarrow0\quad\text{as }j\to\infty.
\end{equation}
Consequently, by \eqref{eq:bounded-dirichlet-approximation},
\begin{equation}\label{eq:bounded-approximant-L2}
\begin{aligned}
j^{-n}\int_{B_j}\mathsf d(\widehat W_j,Q_0\llbracket0\rrbracket)^2
&\le2j^{-n}\int_{B_j}\mathsf d(\widehat W_j,\widehat F_j)^2
 +2j^{-n}\int_{B_j}\mathsf d(\widehat F_j,Q_0\llbracket0\rrbracket)^2\\
&\le C(n)A_0^2.
\end{aligned}
\end{equation}
By \cite[Proposition 3.1]{DSQ}, the mean-value inequality and
\eqref{eq:bounded-approximant-L2},
\[
\mathsf d(\widehat W_j(0),Q_0\llbracket0\rrbracket)^2
\le C(n)j^{-n}\int_{B_j}\mathsf d(\widehat W_j,Q_0\llbracket0\rrbracket)^2
\le C(n)A_0^2.
\]
Hence, up to a subsequence, there exists
$\widehat F\in\A_{Q_0}(\R^n)$ such that $ \mathsf d(\widehat W_j(0),\widehat F)\to0$ as $j\to\infty$.   Using \eqref{eq:bounded-approximant-L2}, Caccioppoli's inequality
and \cite[Theorem 3.9]{DSQ}, for some $\alpha_0>0$,
\[
0\le\limsup_{j\to\infty}\sup_{p\in B_2}
\mathsf d(\widehat W_j(p),\widehat W_j(0))
\le C(n)A_0\lim_{j\to\infty}j^{-\alpha_0}=0.
\]
Thus $\widehat W_j\to\widehat F$ uniformly on $B_2$ as $j\to\infty$.
Together with
\eqref{eq:bounded-dirichlet-approximation}, we obtain
\begin{equation}\label{eq:bounded-fixed-ball-limit}
\begin{aligned}
\int_{B_2}\mathsf d(\widehat F_j,\widehat F)^2
&\le2\int_{B_2}\mathsf d(\widehat F_j,\widehat W_j)^2
+2|B_2|\sup_{p\in B_2}\mathsf d(\widehat W_j(p),\widehat F)^2 \longrightarrow0\quad\text{as }j\to\infty.
\end{aligned}
\end{equation}
Since $O(B_2)=B_2$ for $O\in G$, \eqref{eq:bounded-slice-complement} yields
\[
\sum_j|B_2\setminus(\widehat K_j\cap O^{-1}\widehat K_j)|
\le 2\sum_j|B_2\setminus\widehat K_j|<\infty.
\]
By \eqref{eq:bounded-fixed-ball-limit} and the Borel--Cantelli lemma,
after passing to a subsequence we have, for every $O\in G$ and a.e. $p\in B_2$,
\[
\widehat F_j(p)\to\widehat F,\qquad
\widehat F_j(Op)\to\widehat F,\qquad
\mathbf 1_{\widehat K_j\cap O^{-1}\widehat K_j}(p)\to1
\qquad\text{as }j\to\infty.
\]
Using \eqref{eq:bounded-slice-approximation} and $G$-symmetry of $S$,
we have $\widehat F_j(Op)=O\widehat F_j(p)$ for a.e.
$p\in B_2\cap\widehat K_j\cap O^{-1}\widehat K_j$.
Letting $j\to\infty$, we obtain $O\widehat F=\widehat F$, $\forall O\in G$.
For $\mathbf a=(a_1,\ldots,a_n)\in\R^n\setminus\{0\}$, let $k_{\mathbf a}\coloneqq\#\{i:a_i\ne0\}$.
Note that, 
\begin{equation}\label{eq:nonzero-orbit-size}
\#(G\mathbf a)\ge2^{k_{\mathbf a}}\binom n{k_{\mathbf a}}
=\frac{n2^{k_{\mathbf a}}}{k_{\mathbf a}}\binom{n-1}{k_{\mathbf a}-1}\ge2n.
\end{equation}
Since $\widehat F$ is $G$-invariant and has multiplicity $Q_0<2n$, we obtain ${\widehat F}=Q_0\llbracket0\rrbracket$
from \eqref{eq:nonzero-orbit-size}.
For a.e. $p\in B_2\setminus\{0\}$, we have
$\mathbf 1_{\widehat K_j}(p)\to1$ and $|\widehat R_jp|\to\infty$.
By \eqref{eq:bounded-slice-approximation} and \eqref{eq:center-cylinder},
\[
\liminf_{j\to\infty}\mathsf d(\widehat F_j(p),\widehat F)
=\liminf_{j\to\infty}\mathsf d(\widehat F_j(p),Q_0\llbracket0\rrbracket)
\ge\sqrt{Q_0}\,c>0,
\]
contradicting $\widehat F_j(p)\to\widehat F$ for a.e. $p\in B_2$.
Hence $\pi_x(\spt S)$ is unbounded.
\end{proof}

Now, we turn to the proof of Proposition \ref{prop:first-limit}. 
\begin{proof}[Proof of Proposition~\ref{prop:first-limit}]
\textup{(i)} By uniqueness in \eqref{eq:dirichlet},
$v^{(t)}\circ O=v^{(t)}$,  $\forall O\in G$.
By coordinate reflections, $Dv^{(t)}(0)=0$ and
$\partial_{ij}v^{(t)}(0)=0$ for $i\ne j$.
Hence, by coordinate permutations and \eqref{eq:dirichlet} at $x=0$,
\[
D^2v^{(t)}(0)
=\partial_{11}v^{(t)}(0)I
=\tan\Bigl(\frac{\Theta(D^2v^{(t)}(0))}{n}\Bigr)I
=\cot(\pi/n)I.
\]
Expanding $v^{(t)}$ at $0$ with $Dv^{(t)}(0)=0$, we obtain
\[
\frac{v^{(t)}(ry)-v^{(t)}(0)}{r^2}
=\frac{\cot(\pi/n)}2|y|^2+o(1)
\qquad\text{uniformly for }y\in B_1\text{ as }r\downarrow0.
\]
Then
\[
\lim_{r\downarrow0}\frac{\osc_{B_r}v^{(t)}}{r^2}=\frac{\cot(\pi/n)}2<\cot(\pi/n).
\]
Thus $r_t>0$ by \eqref{eq:first-radius}.

Next, if $r_{t_j}\ge c_0>0$ for some $t_j\to\infty$, then
\[
\|v^{(t_j)}-v^{(t_j)}(0)\|_{C^0(B_{c_0}(0))}\le\cot(\pi/n) c_0^2.
\]
By Lemma~\ref{lem:interior-estimate}, the sequence
$\{v^{(t_j)}-v^{(t_j)}(0)\}$ is precompact in
$C^\infty_{\mathrm{loc}}(B_{c_0}(0))$.
For each fixed $R>0$, $R/t_j\to0$ as $j\to\infty$. Thus
\eqref{eq:ball-area} and \eqref{eq:ball-boundary} imply
\[
\limsup_{j\to\infty}\frac{\|M_{t_j}\|(\mathbf B_R)}{\omega_nR^n}\le2n-1,
\qquad
\dist(0,\spt\partial M_{t_j})\longrightarrow\infty
\quad\text{as }j\to\infty.
\]
Using current compactness \cite[Theorem 27.3]{Simon} and
Lemma~\ref{lem:interior-estimate}, we obtain, up to a subsequence,
$M_{t_j}\rightharpoonup M$ locally as currents as $j\to\infty$,
with $M\ne0$, $\partial M=0$ and $\pi_x(\spt M)\subset\overline B_1$.
By lower semicontinuity,
\[
\|M\|(\mathbf B_R)\le(2n-1)\omega_nR^n\qquad\forall R>0.
\]
Using $M_{t_j}=\Graph(Dv^{(t_j)})$, $M_{t_j}\rightharpoonup M$ locally as
currents as $j\to\infty$, and Lemma~\ref{lem:interior-estimate}, we obtain
\[
\begin{aligned}
\sup_{(x,p)\in\spt M,\ |x|<c_0/2}|p|
&\le\limsup_{j\to\infty}\sup_{x\in B_{c_0/2}(0)}
|Dv^{(t_j)}(x)| \le C(n,c_0).
\end{aligned}
\]
Thus \eqref{eq:center-cylinder} holds with $S=M$ and $c=c_0/2$.
Finally, by $v^{(t_j)}\circ O=v^{(t_j)}$ and
$M_{t_j}\rightharpoonup M$ locally as currents as $j\to\infty$,
\[
(O\times O)_\sharp M
=\lim_{j\to\infty}(O\times O)_\sharp M_{t_j}
=\lim_{j\to\infty}(\det O)M_{t_j}
=(\det O)M\qquad\forall O\in G.
\]
By Lemma~\ref{lem:bounded-source}, $\pi_x(\spt M)$ is unbounded, which is a contradiction with $\pi_x(\spt M)\subset\overline B_1$.
Thus $r_t\to0$ as $t\to\infty$, which proves \textup{(i)}.

\textup{(ii)} Let $t_j\to\infty$ as $j\to\infty$. Since $r_{t_j}\to0$,
we may discard finitely many terms so that $r_{t_j}<1/2$.
For $t=t_j$, the infimum in \eqref{eq:first-radius} equals $r_t$ because
$r_t<1/2$. By continuity of $\osc_{B_r}v^{(t)}$ in $r$,
$\osc_{B_{r_t}}v^{(t)}=\cot(\pi/n)r_t^2$.
Moreover, by \eqref{eq:rescaled-potentials},
\begin{equation}\label{eq:gradient-rescaling}
Du_t(y)=r_t^{-1}Dv^{(t)}(r_ty)
=\frac{t}{r_t}DU_t(r_ty).
\end{equation}
By \eqref{eq:gradient-rescaling}, $\Graph(Du_t)=(r_t^{-1}\operatorname{Id})_\sharp M_t$.
Rescaling \eqref{eq:ball-area}, we obtain
\begin{equation}\label{eq:first-area}
\limsup_{t\to\infty}\sup_{0<R<\delta t/(2r_t)}
\frac{\|\Graph(Du_t)\|(\mathbf B_R)}{\omega_nR^n}\le2n-1.
\end{equation}
Moreover, by \eqref{eq:ball-boundary},
\begin{equation}\label{eq:first-boundary}
\liminf_{t\to\infty}\frac{r_t}{t}\dist(0,\spt\partial\Graph(Du_t))
=\liminf_{t\to\infty}\frac{\dist(0,\spt\partial M_t)}t>0.
\end{equation}
Since $r_{t_j}R/t_j\to0$ as $j\to\infty$ for fixed $R>0$, using
\eqref{eq:first-area}, \eqref{eq:first-boundary} and current compactness
\cite[Theorem 27.3]{Simon}, we obtain a subsequence,
still denoted by $t_j$, with $\Graph(Du_{t_j})\rightharpoonup S$ as $j\to\infty$,
where $S$ is a current without boundary and
$\|S\|(\mathbf B_R)\le(2n-1)\omega_nR^n$.
Since $\Graph(Du_{t_j})$ has phase $(n-2)\pi/2$, the current limit $S$
also has phase $(n-2)\pi/2$.
For $\chi\in C_c^\infty(\C^n)$,
\begin{equation}\label{eq:first-mass-convergence}
\begin{aligned}
\|\Graph(Du_{t_j})\|(\chi)-\|S\|(\chi)
&=(\Graph(Du_{t_j})-S)\Bigl(\chi\operatorname{Re}
(e^{-i(n-2)\pi/2}dz_1\wedge\cdots\wedge dz_n)\Bigr)\\
&\longrightarrow0\quad\text{as }j\to\infty.
\end{aligned}
\end{equation}
If $z_j\in\spt \Graph(Du_{t_j})$ and $z_j\to z$ as $j\to\infty$, then, combining
\eqref{eq:first-boundary}, monotonicity and \eqref{eq:first-mass-convergence},
we obtain
\[
0<\omega_n\rho^n
\le\limsup_{j\to\infty}\|\Graph(Du_{t_j})\|(\mathbf B_\rho(z_j))
\le\|S\|(\overline{\mathbf B_{2\rho}(z)})
\qquad\forall\rho>0.
\]
Thus $z\in\spt S$.
Conversely, by local current convergence, we have
$\dist(z,\spt \Graph(Du_{t_j}))\to0$ as $j\to\infty$, $\forall z\in\spt S$.
Consequently,
$\spt \Graph(Du_{t_j})\to\spt S$ locally in Hausdorff distance as $j\to\infty$.

On the other hand, using $u_t(0)=0$ and $\osc_{B_1}u_t=\cot(\pi/n)$ in Lemma~\ref{lem:interior-estimate},
up to a subsequence we have
\[
u_{t_j}\longrightarrow u\quad\text{in }C^k(\overline B_{1/2}(0))
\quad\text{as }j\to\infty,\qquad\forall k\ge0.
\]
Since $\limsup_{j\to\infty}\|Du_{t_j}\|_{C^0(\overline B_{1/2}(0))}<\infty$,
we obtain from $\Graph(Du_{t_j})\rightharpoonup S$ locally as currents as $j\to\infty$
\begin{equation}\label{eq:graph-near-origin}
S\res\bigl(B_{1/2}(0)\times\R^n\bigr)
=\Graph\bigl(Du\big|_{B_{1/2}(0)}\bigr).
\end{equation}
By \eqref{eq:graph-near-origin}, $(0,0)$ is a regular point of $S$.
Let $\Lambda$ be the connected component of $\spt S\setminus\Sing S$
containing $(0,0)$.
As $\Lambda$ is closed in $\spt S\setminus\Sing S$,
$\overline\Lambda\setminus\Lambda\subset\Sing S$.
Hence, by \cite[Theorem 1.1]{DMS}, $\mathcal H^{n-1}(\overline\Lambda\setminus\Lambda)=0$.  
For any compactly supported smooth $(n-1)$-form $\varphi$ and any $\epsilon>0$,
we cover $(\overline\Lambda\setminus\Lambda)\cap\spt\varphi$ by finitely
many balls $\mathbf B_{r_i}(z_i)$ with $0<r_i<1$ and
$\sum_i r_i^{n-1}<\epsilon$, and  choose
$\psi_i\in C_c^\infty(\mathbf B_{2r_i}(z_i))$ in $\C^n$ with
\[
0\le\psi_i\le1,\qquad
\psi_i|_{\mathbf B_{r_i}(z_i)}=1,\qquad |D\psi_i|\le2/r_i,
\]
such that $\prod_i(1-\psi_i)\varphi$ has compact support on $\Lambda$.
Using $0\le1-\prod_i(1-\psi_i)\le\sum_i\psi_i$,
$\|S\res\Lambda\|\le\|S\|$ and monotonicity for $S$, we obtain
\begin{equation}\label{eq:component-boundary-estimate}
\begin{aligned}
|\partial(S\res\Lambda)(\varphi)|
&=\Bigl|(S\res\Lambda)\bigl(d[(1-\prod_i(1-\psi_i))\varphi]\bigr)\Bigr|\\
&=\Bigl|(S\res\Lambda)\bigl((1-\prod_i(1-\psi_i))d\varphi
 +\sum_i\prod_{k\ne i}(1-\psi_k)\,d\psi_i\wedge\varphi\bigr)\Bigr|\\
&\le\|d\varphi\|_\infty\sum_i\int\psi_i\,d\|S\res\Lambda\|
 +\|\varphi\|_\infty\sum_i\int|d\psi_i|\,d\|S\res\Lambda\|\\
&\le\sum_i\bigl(\|d\varphi\|_\infty
               +2r_i^{-1}\|\varphi\|_\infty\bigr)
                  \|S\res\Lambda\|(\mathbf B_{2r_i}(z_i))\\
&\le C(n,\varphi)\sum_i r_i^{-1}\|S\res\Lambda\|(\mathbf B_{2r_i}(z_i))
 \le C(n,\varphi)\sum_i r_i^{n-1}<C(n,\varphi)\epsilon,
\end{aligned}
\end{equation}
where $C(n,\varphi)$ is independent of $\epsilon$.
Letting $\epsilon\downarrow0$, we obtain $\partial(S\res\Lambda)(\varphi)=0$.
Thus $\partial(S\res\Lambda)=0$.
Since $(O\times O)_\sharp S=(\det O)S$ and $O(0)=0$, $\forall O\in G$,
$\Lambda$ is $G$-invariant. Hence $S\res\Lambda$ has the calibration and
$G$-symmetry of $S$.
Moreover, \eqref{eq:center-cylinder} holds for $S\res\Lambda$ by
\eqref{eq:graph-near-origin}.
By $\|S\res\Lambda\|\le\|S\|$, \eqref{eq:first-area}
and Lemma~\ref{lem:bounded-source},
$\pi_x\bigl(\spt(S\res\Lambda)\bigr)$ is unbounded, and hence
$\pi_x(\Lambda)$ is unbounded.

For $n\ge1$ and $\theta\in\R$, we denote by $\Gr_n(\theta)$ the closure
in the oriented Grassmannian of Lagrangian graph planes with principal phase $\theta$, oriented by $x$.
Let $L\in\Gr_n(\theta)$. The principal angles
$\theta_1,\ldots,\theta_n$ can be chosen so that
\begin{equation}\label{eq:closed-principal-angles}
-\frac\pi2\le\theta_i\le\frac\pi2\quad(1\le i\le n),
\qquad\sum_{i=1}^n\theta_i=\theta.
\end{equation}
Next, for each $z_0\in\Lambda$, we choose a ball
$\mathbf B(z_0)\subset\C^n$ with
$\spt S\cap\mathbf B(z_0)\subset\Lambda$.
Let $\xi$ be the oriented unit tangent $n$-vector of $S$ on
$\Lambda\cap\mathbf B(z_0)$, smoothly extended to $\mathbf B(z_0)$ with $|\xi|=1$,
and let $\xi_{t_j}$ be the oriented unit tangent $n$-vector of $\Graph(Du_{t_j})$.
Taking $\chi\in C_c^\infty(\mathbf B(z_0))$ with $0\le\chi\le1$, and using
$\Graph(Du_{t_j})\rightharpoonup S$ locally as currents as $j\to\infty$
and \eqref{eq:first-mass-convergence}, we obtain
\begin{equation}\label{eq:gauss-plane-convergence}
\begin{aligned}
\int\chi|\xi_{t_j}-\xi|^2\,d\|\Graph(Du_{t_j})\|
&=2\int\chi\bigl(1-\langle\xi_{t_j},\xi\rangle\bigr)\,d\|\Graph(Du_{t_j})\|\\
&\longrightarrow2\int\chi\bigl(1-|\xi|^2\bigr)\,d\|S\|=0
\quad\text{as }j\to\infty.
\end{aligned}
\end{equation}
Since $\xi_{t_j}\in\Gr_n((n-2)\pi/2)$ and $\Gr_n((n-2)\pi/2)$ is closed,
\eqref{eq:gauss-plane-convergence} implies
\[
T_{z_0}\Lambda\in\Gr_n((n-2)\pi/2)\qquad\forall z_0\in\Lambda.
\]
By Proposition~\ref{prop:finite-projection}, local invertibility of
$\pi_x$ at $(0,0)\in\Lambda$ extends to all of $\Lambda$.
The graph $\Graph(Du_t)$ has multiplicity one, so for
$\varrho\in C_c^\infty(\R^n)$ and $\chi\in C_c^\infty(\C^n)$ we have
\begin{equation}\label{eq:one-sheet}
0\le \Graph(Du_t)\bigl(\varrho(x)\chi(x,p)\,dx_1\wedge\cdots\wedge dx_n\bigr)
\le\int_{\R^n}\varrho(x)\,dx,
\qquad \varrho\ge0,\quad0\le\chi\le1,
\end{equation}
By local current convergence, \eqref{eq:one-sheet} also holds for
$S$. If $\pi_x$ is not injective on $\Lambda$, then, since $\pi_x$ is locally
invertible, $\Lambda$ contains two disjoint graphs over the same ball $B\subset\R^n$.
For $0\le \varrho\in C_c^\infty(B)$ with $\int_B\varrho>0$, we choose
$\chi\in C_c^\infty(\C^n)$ with $0\le\chi\le1$ and $\chi=1$ on
both graphs over $\spt \varrho$. By \eqref{eq:one-sheet},
\[
2\int_B\varrho\le S\bigl(\varrho(x)\chi(x,p)\,dx_1\wedge\cdots\wedge dx_n\bigr)
\le\int_B\varrho,
\]
a contradiction. Hence $\pi_x|_\Lambda$ is a diffeomorphism onto
the connected unbounded open set $\pi_x(\Lambda)\subset\R^n$ containing $0$.
By \eqref{eq:one-sheet}, $S$ has multiplicity one on $\Lambda$.
Since $\Theta(D^2u_{t_j})=(n-2)\pi/2$, the graphs $\Graph(Du_{t_j})$ are minimal.
Moreover, by \eqref{eq:first-boundary},
$\dist(0,\spt\partial\Graph(Du_{t_j}))\to\infty$ as $j\to\infty$.
By \eqref{eq:first-mass-convergence}, \cite[Theorem 24.2]{Simon}
and local invertibility of $\pi_x$, we obtain
\[
Du_{t_j}\longrightarrow p\quad\text{locally uniformly on }\pi_x(\Lambda)
\quad\text{as }j\to\infty,
\]
where $(x,p(x))\in\Lambda$ for $x\in\pi_x(\Lambda)$.
Integrating along a piecewise $C^1$ path $\gamma_x$ in $\pi_x(\Lambda)$
from $0$ to $x$ and using $u_{t_j}(0)=0$, we have
\[
u(x)\coloneqq\lim_{j\to\infty}u_{t_j}(x)
=\lim_{j\to\infty}\int_{\gamma_x}Du_{t_j}(y)\cdot dy
=\int_{\gamma_x}p(y)\cdot dy.
\]
Thus $u$ extends the $C^\infty$ limit of $\{u_{t_j}\}$ on $B_{1/2}(0)$
and satisfies $Du=p$, so $\Lambda=\Graph(Du)$.
Moreover, $u_{t_j}\to u$ locally uniformly by the local uniform convergence
$Du_{t_j}\to p$. By Lemma~\ref{lem:interior-estimate}, we obtain
\begin{equation}\label{eq:local-potential-limit}
u_{t_j}\longrightarrow u\quad\text{in }C^\infty_{\mathrm{loc}}(\pi_x(\Lambda))
\quad\text{as }j\to\infty.
\end{equation}
Thus \textup{(ii)} holds by \eqref{eq:local-potential-limit}
and $\Graph(Du_{t_j})\rightharpoonup S$.
\end{proof}

\section{The tangent cone at infinity}\label{sec:cones}
\label{sec:moving}
\label{subsec:splitting}
\label{subsec:models}

In this section, we prove the following.
\begin{proposition}\label{prop:noescape}
Let $S$ be any local current limit of $\Graph(Du_t)$ as $t\to\infty$. Then
\begin{enumerate}[label=\textup{(\roman*)}]
\item $\pi_p:\spt S\to\R^n$ is proper, and
\[
(\pi_p)_\sharp S=-(2n-1)\llbracket\R^n\rrbracket.
\]
\item $S$ has the unique tangent cone $-(2n-1)\llbracket\{0\}\times\R^n\rrbracket$ at infinity.
\end{enumerate}
\end{proposition}
We first decompose the possible tangent cones into calibrated planes.
\begin{proposition}\label{prop:planar}
Let $n\ge3$, and let $T$ be a locally integral cone in $\C^n$,
without boundary, calibrated by
$\operatorname{Re}(e^{-i(n-2)\pi/2}dz_1\wedge\cdots\wedge dz_n)$.
Assume that $\|T\|(\mathbf B_1)<\infty$ and that
\[
T_z(\spt T)\in\Gr_n((n-2)\pi/2)
\qquad\forall z\in\spt T\setminus\Sing T.
\]
Then $T$ is a finite sum of calibrated linear planes with positive
integer multiplicities.
\end{proposition}
\begin{proof}
We first consider a locally integral $n$-cone $T$ without boundary,
$n\ge3$, with $\|T\|(\mathbf B_1)<\infty$, calibrated by
$\operatorname{Re}(e^{-i\theta}dz_1\wedge\cdots\wedge dz_n)$,
where $(n-2)\pi/2\le\theta<n\pi/2$.
Let $\Lambda$ be a connected component of $\spt T\setminus\Sing T$,
locally represented by $\Graph(Du)$ with $\Theta(D^2u)=\theta$.
We claim that $\lambda_{\max}(D^2u)$ is constant on $\Lambda$
and that $T\res\Lambda$ splits off a fixed line.

Let $\Delta_\Lambda$ and $\nabla_\Lambda$ denote the Laplacian and
gradient with respect to the metric induced on $\Lambda$. 
Since $\Lambda$ is open in $\spt T$ and closed in
$\spt T\setminus\Sing T$, $\overline\Lambda\setminus\Lambda$ is a closed
subset of $\Sing T$. By \cite[Theorem 1.1]{DMS}, $\overline\Lambda\setminus\Lambda$ has $\sigma$-finite
$\mathcal H^{n-2}$ measure. Hence $\mathcal H^{n-1}(\overline\Lambda\setminus\Lambda)=0$.

Fix $0<r<R$, and let
$K\subset (\overline\Lambda\setminus\Lambda)\cap\{r\le|z|\le R\}$ be compact with
$\mathcal H^{n-2}(K)<\infty$. For $0<h<r/4$, we cover $K$ by finitely many balls
$\mathbf B_{r_i}(z_i)$ with $z_i\in K$, $r_i<h$ and
$\sum_i r_i^{n-2}\le C(n,K)$, and choose functions
$\psi_i\in C_c^\infty(\mathbf B_{2r_i}(z_i))$ in $\C^n$ satisfying
\[
0\le\psi_i\le1,\qquad
\psi_i|_{\mathbf B_{r_i}(z_i)}=1,\qquad |D\psi_i|\le2/r_i.
\]
We then set $\zeta_h\coloneqq\max_i\psi_i$. Thus $0\le\zeta_h\le1$,
$\zeta_h=1$ near $K$ and
$\spt\zeta_h\subset\{r/2<|z|<2R\}$. Therefore
\begin{equation}\label{eq:cutoff-energy-control}
\begin{aligned}
\int\bigl(\zeta_h^2+h^2|\nabla_\Lambda\zeta_h|^2\bigr)\,d\|T\res\Lambda\|
&\le\sum_i\int\bigl(\psi_i^2+h^2|\nabla_\Lambda\psi_i|^2\bigr)\,d\|T\res\Lambda\|\\
&\le C(n)h^2\sum_i r_i^{-2}\|T\res\Lambda\|\bigl(\mathbf B_{2r_i}(z_i)\bigr)\\
&\le C(n,r,R,T)h^2\sum_i r_i^{n-2}
\le C(n,r,R,T)C(n,K)h^2.
\end{aligned}
\end{equation}
Since the covering radii can be made arbitrarily small, choose
$h_\ell<2^{-\ell}$ with $\spt\zeta_{h_{\ell+1}}$
contained in the open set where $\zeta_{h_\ell}=1$.
By \eqref{eq:cutoff-energy-control} and $\sum_\ell h_\ell<\infty$,
the $L^2(\|T\res\Lambda\|)$ norm of $N^{-1}\sum_{\ell=1}^N\zeta_{h_\ell}$
is at most $C(n,r,R,K,T)/N$.
The gradients $\nabla_\Lambda\zeta_{h_\ell}$ have disjoint supports and
bounded $L^2(\|T\res\Lambda\|)$ norms, so
\[
\Bigl\|\frac1N\sum_{\ell=1}^N\zeta_{h_\ell}\Bigr\|_{L^2(\|T\res\Lambda\|)}^2
+\frac1{N^2}\sum_{\ell=1}^N
 \|\nabla_\Lambda\zeta_{h_\ell}\|_{L^2(\|T\res\Lambda\|)}^2
\le\frac{C(n,r,R,K,T)}N.
\]
$N^{-1}\sum_{\ell=1}^N\zeta_{h_\ell}=1$ near $K$, so
$\operatorname{Cap}_{W^{1,2}(\|T\res\Lambda\|)}(K)=0$.

Next, by the $\sigma$-finiteness of $\mathcal H^{n-2}$ on $\overline\Lambda\setminus\Lambda$,
there exist compact sets  $K_\ell\subset (\overline\Lambda\setminus\Lambda)\cap \{r\le|z|\le R\}$ with
$\mathcal H^{n-2}(K_\ell)<\infty$ and a set ${\mathcal N}\subset (\overline\Lambda\setminus\Lambda)\cap \{r\le|z|\le R\}$ such that
\[
(\overline\Lambda\setminus\Lambda)\cap \{r\le|z|\le R\}={\mathcal N}\cup\bigcup_{\ell=1}^\infty K_\ell,\qquad
\mathcal H^{n-2}({\mathcal N})=0.
\]
On the other hand, for any $\epsilon>0$, there exists  
$\mathbf B_{r_i}(z_i)$,
$z_i\in {\mathcal N}$, such that
\[
{\mathcal N}\subset\bigcup_i\mathbf B_{r_i}(z_i),\qquad
0<r_i<\min\{1,r/4\},\qquad
\sum_i r_i^{n-2}<\epsilon.
\]
Taking $\psi_i\in C_c^\infty(\mathbf B_{2r_i}(z_i))$ with 
$0\le\psi_i\le1$, $\psi_i|_{\mathbf B_{r_i}(z_i)}=1$ and
$|D\psi_i|\le2/r_i$, we obtain
\[
\begin{aligned}
0&\le\operatorname{Cap}_{W^{1,2}(\|T\res\Lambda\|)}({\mathcal N})
\le\sum_i\operatorname{Cap}_{W^{1,2}(\|T\res\Lambda\|)}
\bigl(\mathbf B_{r_i}(z_i)\bigr)\\
&\le\sum_i\int\bigl(\psi_i^2+|\nabla_\Lambda\psi_i|^2\bigr)\,d\|T\res\Lambda\| \le\sum_i(1+4r_i^{-2})\|T\res\Lambda\|\bigl(\mathbf B_{2r_i}(z_i)\bigr)\\
&\le C(n,r,R,T)\sum_i(r_i^2+4)r_i^{n-2}
\le5C(n,r,R,T)\sum_i r_i^{n-2}<5C(n,r,R,T)\epsilon
\longrightarrow0\quad\text{as }\epsilon\downarrow0.
\end{aligned}
\]
Since $\operatorname{Cap}_{W^{1,2}(\|T\res\Lambda\|)}(K_\ell)=0$, $\forall\ell\ge1$,
\[
\begin{aligned}
\operatorname{Cap}_{W^{1,2}(\|T\res\Lambda\|)}\bigl((\overline\Lambda\setminus\Lambda)\cap \{r\le|z|\le R\}\bigr)
&\le\operatorname{Cap}_{W^{1,2}(\|T\res\Lambda\|)}({\mathcal N}) +\sum_{\ell=1}^\infty
\operatorname{Cap}_{W^{1,2}(\|T\res\Lambda\|)}(K_\ell)=0.
\end{aligned}
\]
We therefore can choose compactly supported Lipschitz functions
$\zeta_j:\C^n\to[0,1]$ with
$\spt\zeta_j\subset\{r/2<|z|<2R\}$ and $\zeta_j=1$ near
$(\overline\Lambda\setminus\Lambda)\cap \{r\le|z|\le R\}$, such that
\begin{equation}\label{eq:cutoff-strong-convergence}
\|\zeta_j\|_{W^{1,2}(\|T\res\Lambda\|)}\longrightarrow0
\quad\text{as }j\to\infty.
\end{equation}
Then, up to a subsequence, $\zeta_j(z)\to0$ as $j\to\infty$
for $\|T\res\Lambda\|$-a.e. $z\in\Lambda$.

Now, on $\Lambda$, we define
\[
{b_\Lambda}\coloneqq (1+\lambda_{\max}(D^2u)^2)^{-\frac14+\frac1{\sqrt{4n+1}+1}}.
\]
By \cite[Proposition 2.1]{WY}, 
\begin{equation}\label{eq:Jacobi}
\Delta_\Lambda {b_\Lambda}\le0,\qquad0<{b_\Lambda}\le1.
\end{equation}
Note that, since $T$ is a cone, $\lambda\Lambda=\Lambda$, $\forall\lambda>0$, and
${b_\Lambda}(\lambda z)={b_\Lambda}(z)$ for $z\in\Lambda$. Hence
$z\cdot\nabla_\Lambda {b_\Lambda}=0$.
For a radial $0\le\chi_0\in C_c^\infty(\{z\in\C^n:r<|z|<R\})$,
we test \eqref{eq:Jacobi} with $\chi_0(1-\zeta_j)^2(1-b_\Lambda)$,
which has compact support on $\Lambda$. Since $\chi_0$ is radial and
$z\cdot\nabla_\Lambda b_\Lambda=0$,
$\nabla_\Lambda b_\Lambda\cdot\nabla_\Lambda\chi_0=0$. Integrating
against $d\|T\res\Lambda\|$ and applying Young's inequality with
$0\le1-b_\Lambda\le1$, we obtain
\[
\begin{aligned}
\int\chi_0(1-\zeta_j)^2|\nabla_\Lambda b_\Lambda|^2
&\le-2\int\chi_0(1-\zeta_j)(1-b_\Lambda)
   \nabla_\Lambda b_\Lambda\cdot\nabla_\Lambda\zeta_j\\
&\le\frac12\int\chi_0(1-\zeta_j)^2|\nabla_\Lambda b_\Lambda|^2
+2\int\chi_0|\nabla_\Lambda\zeta_j|^2.
\end{aligned}
\]
By Fatou's lemma and \eqref{eq:cutoff-strong-convergence},
\begin{equation}\label{eq:jacobi-local-energy}
\int\chi_0|\nabla_\Lambda b_\Lambda|^2
\le4\lim_{j\to\infty}\int\chi_0|\nabla_\Lambda\zeta_j|^2=0.
\end{equation}
Since $\chi_0$ can be chosen positive near each $z\in\Lambda\setminus\{0\}$,
$\nabla_\Lambda b_\Lambda=0$ a.e. on $\Lambda$. By connectedness,
$b_\Lambda$ is constant on $\Lambda$.
The phase $\Theta(D^2u)=\theta>0$ implies $\lambda_{\max}(D^2u)>0$.
Hence $\lambda_{\max}(D^2u)$ is constant on $\Lambda$ by the definition of $b_\Lambda$.

Next, near any $(x_0,Du(x_0))\in\Lambda$, we represent $\Lambda$ by
$\{(x,Du(x)):x\in B(x_0)\subset\R^n\}$. Let ${e_*}$ 
be a unit eigenvector of $D^2u(x_0)$
corresponding to $\lambda_{\max}(D^2u(x_0))$, then
\begin{equation}\label{eq:phase-tangent-variation}
0=\partial_{e_*}\Theta(D^2u)
 =\sum_{i,j=1}^n\bigl(I+(D^2u)^2\bigr)^{-1}_{ij}
                    \partial_{e_*}u_{ij}.
\end{equation}
Differentiating \eqref{eq:phase-tangent-variation} in the direction $e_*$, we obtain
\[
\begin{aligned}
\sum_{i,j=1}^n\bigl(I+(D^2u)^2\bigr)^{-1}_{ij}
                  \partial_{ij}u_{{e_*}{e_*}}
&=\partial_{e_*e_*}\Theta(D^2u) -D^2\Theta(D^2u)
       [\partial_{e_*}D^2u,\partial_{e_*}D^2u]\\
&=-D^2\Theta(D^2u)
       [\partial_{e_*}D^2u,\partial_{e_*}D^2u]\ge0.
\end{aligned}
\]
Since $\lambda_{\max}(D^2u)$ is constant on $B(x_0)$, we have
$u_{{e_*}{e_*}}\le\lambda_{\max}(D^2u(x_0))$ on $B(x_0)$ and
$u_{{e_*}{e_*}}(x_0)=\lambda_{\max}(D^2u(x_0))$.
By the strong maximum principle,
$u_{{e_*}{e_*}}=\lambda_{\max}(D^2u)$ on $B(x_0)$.
Using $\lambda_{\max}(D^2u)I-D^2u\ge0$, we have
\[
\bigl(\lambda_{\max}(D^2u)I-D^2u\bigr){e_*}=0.
\]
Consequently,
\[
\begin{aligned}
\partial_{e_*}(x,Du(x))
&=(e_*,D^2u(x)e_*) =(e_*,\lambda_{\max}(D^2u(x_0))e_*)
\in T_{(x,Du(x))}\Lambda,
\qquad\forall x\in B(x_0).
\end{aligned}
\]
By analytic continuation from the graph over $B(x_0)$ to the connected
real analytic component $\Lambda$,
\begin{equation}\label{eq:fixed-eigenvector-tangency}
\bigl({e_*},\lambda_{\max}(D^2u(x_0)){e_*}\bigr)\in T_z\Lambda
\qquad\forall z\in\Lambda.
\end{equation}
Set
\[
{v_\Lambda}\coloneqq
\frac{\bigl({e_*},\lambda_{\max}(D^2u(x_0)){e_*}\bigr)}
{\sqrt{1+\lambda_{\max}(D^2u(x_0))^2}}\,\in T_z\Lambda.
\]
Applying the proof of \eqref{eq:component-boundary-estimate} to $T$
in place of $S$, we have $\partial(T\res\Lambda)=0$.
For a compactly supported smooth $n$-form $\varphi$ on $\C^n$, we obtain
\[
\begin{aligned}
\frac{d}{d\lambda}\bigl[(\operatorname{Id}+\lambda {v_\Lambda})_\sharp(T\res\Lambda)\bigr](\varphi)
&=(T\res\Lambda)\bigl(
 d\iota_{v_\Lambda}((\operatorname{Id}+\lambda {v_\Lambda})^*\varphi)
 +\iota_{v_\Lambda}d((\operatorname{Id}+\lambda {v_\Lambda})^*\varphi)\bigr)\\
&=\partial(T\res\Lambda)
\bigl(\iota_{v_\Lambda}((\operatorname{Id}+\lambda {v_\Lambda})^*\varphi)\bigr)=0.
\end{aligned}
\]
The term containing $\iota_{v_\Lambda}d((\operatorname{Id}+\lambda v_\Lambda)^*\varphi)$
is zero on $T\res\Lambda$ because $v_\Lambda\in T_z\Lambda$ by
\eqref{eq:fixed-eigenvector-tangency}. Hence
\begin{equation}\label{eq:component-translations}
(\operatorname{Id}+\lambda {v_\Lambda})_\sharp(T\res\Lambda)=T\res\Lambda
\qquad\forall \lambda\in\R.
\end{equation}
Since $z,v_\Lambda\in T_z\Lambda$ and $\Lambda$ is Lagrangian,
$\langle z,iv_\Lambda\rangle=0$, $\forall z\in\Lambda$.
Combining \eqref{eq:component-translations} with $\langle z,iv_\Lambda\rangle=0$,
we obtain
\begin{equation}\label{eq:component-line-product}
(T\res\Lambda)=\llbracket\R {v_\Lambda}\rrbracket\times T',\qquad
\spt T'\subset(\C {v_\Lambda})^\perp,
\end{equation}
where $T'$ is a locally integral $(n-1)$-cone without boundary calibrated by
$\iota_{v_\Lambda}\operatorname{Re}(e^{-i\theta}dz_1\wedge\cdots\wedge dz_n)$.

We now take $\theta=(n-2)\pi/2$ and prove the proposition by
induction on $n$.
If $T_2$ is a two-dimensional calibrated cone in $\C^2$ and
$z\in\Sing T_2\setminus\{0\}$, then
$\R_+z\subset\Sing T_2$ since $T_2$ is a cone. Thus, by
\cite[Theorem 1.1]{DMS},
\[
1\le\dim_{\mathcal H}\Sing T_2\le0,
\]
a contradiction.
Thus every connected component of $\spt T_2\cap \mathbb S^3$ is a geodesic
of $\mathbb S^3$, and hence lies in a great circle of $\mathbb S^3$.
Since $\partial T_2=0$, by \cite[Theorem 26.27]{Simon}, each connected
component of $\spt T_2\setminus\{0\}$ extends to a plane in $\C^2$.

Now let $\Lambda$ be a connected component of $\spt T\setminus\Sing T$.
If $\operatorname{rank}d(\pi_x|_\Lambda)_z=n$ at some $z\in\Lambda$,
then $\pi_x$ is locally invertible at every point of $\Lambda$ by
Proposition~\ref{prop:finite-projection}.
By \eqref{eq:fixed-eigenvector-tangency}, with respect to the fixed decomposition
$\R e_*\oplus e_*^\perp$ we have
\[
D^2u=\diag\bigl(\lambda_{\max}(D^2u),\mathsf H'\bigr).
\]
Since $\lambda_n(D^2u)>-\lambda_{n-1}(D^2u)\ge-\lambda_{\max}(D^2u)$
by \eqref{eq:closed-principal-angles}, we obtain
\begin{equation}\label{eq:quotient-slopes-phase}
|\lambda_i(\mathsf H')|\le\lambda_{\max}(D^2u),\qquad
\Theta(\mathsf H')=\frac{n-2}{2}\pi
-\arctan\lambda_{\max}(D^2u)>\frac{n-3}{2}\pi.
\end{equation}
By the product decomposition \eqref{eq:component-line-product}, the tangent
planes at regular points of the $(n-1)$-cone $T'$ are limits of
$T_z\Lambda\cap v_\Lambda^\perp$ with $z\in\Lambda$.
By \eqref{eq:quotient-slopes-phase}, every regular tangent plane of $T'$
is a graph with phase greater than $(n-3)\pi/2$.
Repeating \eqref{eq:component-line-product}, we split off $n-2$ lines.
As for $T_2$, the link of the remaining two-dimensional calibrated cone in
$\mathbb S^3$ consists of geodesics, so each regular component of the
remaining cone lies in a plane. Hence $\Lambda$ is contained in a calibrated plane.

Suppose the maximal rank $n_0$ of $d(\pi_x|_\Lambda)$ satisfies $0<n_0<n$.
By Lemma~\ref{lem:partial-propagation}, $n_0\ge2$. On a small neighborhood
$U\subset\Lambda$ of a point of rank $n_0$ with $x\ne0$, the image $\pi_x(U)$ is open
in an affine $n_0$-plane if $n_0\ge3$ and is a minimal surface if $n_0=2$.
Since $(x,p)\in T_{(x,p)}\Lambda$, we have $x\in T_x\pi_x(U)$ for
$x\in\pi_x(U)$. If $n_0\ge3$, the affine plane therefore passes through
the origin. If $n_0=2$, the projected surface is locally a minimal cone.
The link in $\mathbb S^{n-1}$ lies in a great circle, so the projected
surface lies in a linear two-plane.

If $\nu$ is normal to the plane containing $\pi_x(U)$, then
$(\nu,0)\perp T_z\Lambda$ for $z\in U$. Since $\Lambda$ is Lagrangian,
$(0,\nu)\in T_z\Lambda$ on $U$. By analyticity and
\eqref{eq:component-line-product}, we split off the $n-n_0$ lines
from $T\res\Lambda$. At rank $n_0$, the remaining angles sum to
$(n_0-2)\pi/2$ by \eqref{eq:closed-principal-angles}. By density of
the rank-$n_0$ points, the regular tangent planes of the remaining
cone belong to $\Gr_{n_0}((n_0-2)\pi/2)$. We conclude by induction
for $n_0\ge3$ and by the case of $T_2$ for $n_0=2$.
Thus $\Lambda$ is locally planar and, by connectedness, lies in a
single linear $n$-plane. If $d(\pi_x|_\Lambda)=0$, then $x=0$ since
$(x,p)\in T_{(x,p)}\Lambda$.

Thus each connected component $\Lambda_a$ of $\spt T\setminus\Sing T$
lies in a calibrated plane $\Pi_a$. Since $\partial(T\res\Lambda_a)=0$,
$T\res\Lambda_a=m_a\llbracket\Pi_a\rrbracket$ with $m_a\in\mathbb Z_{>0}$ by
\cite[Theorem 26.27]{Simon}. As $\mathcal H^n(\Sing T)=0$,
$T=\sum_a m_a\llbracket\Pi_a\rrbracket$ and
\[
\|T\|(\mathbf B_1)=\omega_n\sum_a m_a<\infty.
\]
Hence the sum over $a$ is finite.
\end{proof}

We now determine the calibrated planes $\Pi_a$ and positive integers $m_a$ in
the finite decomposition from Proposition~\ref{prop:planar},
\[
T=\lim_{k\to\infty}(R_k^{-1}\operatorname{Id})_\sharp\Graph(Du_{t_k})
=\sum_a m_a\llbracket\Pi_a\rrbracket,
\]
where the limit of $(R_k^{-1}\operatorname{Id})_\sharp\Graph(Du_{t_k})$
is a cone and convergence is local in the sense of currents.
\begin{lemma}\label{lem:cone-models}
Let $n\ge4$, $t_k\to\infty$ and $R_k>0$ satisfy $r_{t_k}R_k\to0$.
Suppose that $T$ is a cone and
\[
(R_k^{-1}\operatorname{Id})_\sharp\Graph(Du_{t_k})\rightharpoonup T
\quad\text{locally as currents as }k\to\infty.
\]
Then
\begin{equation}\label{eq:two-models}
T\in\Bigl\{\Iplane,\ -(2n-1)\llbracket\{0\}\times\R^n\rrbracket,\ -(n-1)\llbracket\{0\}\times\R^n\rrbracket+\sum_{i=1}^n\Pplane_i\Bigr\},
\end{equation}
where
\[
\begin{aligned}
\Iplane&\coloneqq\llbracket\{p=\cot(\pi/n)x\}\rrbracket,\\
\Pplane_i&\coloneqq\llbracket\{x_i=0,\ p_j=x_j\cot(\pi/(n-1))\text{ for }j\ne i\}\rrbracket,
\qquad 1\le i\le n,
\end{aligned}
\]
with all planes oriented by
$\operatorname{Re}(e^{-i(n-2)\pi/2}dz_1\wedge\cdots\wedge dz_n)$.
\end{lemma}
\begin{proof}
Let $i,j\in\{1,\ldots,n\}$, $i\ne j$, and define $R_s:\R^n\to\R^n$ by
\[
(R_sx)_k\coloneqq
\begin{cases}
x_i\cos s-x_j\sin s,&k=i,\\
x_i\sin s+x_j\cos s,&k=j,\\
x_k,&k\notin\{i,j\}.
\end{cases}
\]
Since $R_s^TR_s=I$, we obtain from \eqref{eq:dirichlet}
\[
\Theta\bigl(D^2(v^{(t)}\circ R_s)(x)\bigr)
=\Theta\bigl(D^2v^{(t)}(R_sx)\bigr)=\frac{n-2}{2}\pi
\qquad\forall x\in B_1,\ s\in\R.
\]
Moreover,
\[
\frac{d}{ds}v^{(t)}(R_sx)\bigg|_{s=0}
=-x_j\partial_iv^{(t)}+x_i\partial_jv^{(t)}.
\]
Therefore
\begin{equation}\label{eq:linearized-operator}
\begin{aligned}
0&=\frac{d}{ds}\Theta\bigl(D^2(v^{(t)}\circ R_s)\bigr)\bigg|_{s=0} =\tr\bigl((I+(D^2v^{(t)})^2)^{-1}
D^2(x_i\partial_jv^{(t)}-x_j\partial_iv^{(t)})\bigr)
\quad\text{in }B_1.
\end{aligned}
\end{equation}
Since $v^{(t)}$ is invariant under reflection in $x_i=0$ and
interchange of $x_i,x_j$,
\[
x_i\partial_jv^{(t)}-x_j\partial_iv^{(t)}=0
\quad\text{on }B_1\cap(\{x_i=0\}\cup\{x_i=x_j\}).
\]
Moreover, since $x_i\partial_j-x_j\partial_i$ is tangent to $\mathbb S^{n-1}$,
we obtain from the boundary condition in \eqref{eq:dirichlet}
\[
\begin{aligned}
x_i\partial_jv^{(t)}-x_j\partial_iv^{(t)}
&=-t(x_i\partial_j-x_j\partial_i)H_n\\
&=4tx_ix_j(x_i^2-x_j^2)<0
\quad\text{on }\mathbb S^{n-1}\cap\{0<x_i<x_j\}.
\end{aligned}
\]
Applying the maximum principle to \eqref{eq:linearized-operator} on
$B_1\cap\{0<x_i<x_j\}$ with zero data on $x_i=0$ and $x_i=x_j$
and negative data on the spherical face, we obtain
\begin{equation}\label{eq:angular-sign}
x_i\partial_jv^{(t)}-x_j\partial_iv^{(t)}<0\qquad\text{in }B_1\cap\{0<x_i<x_j\}.
\end{equation}

Now let $x\in B_1$ satisfy $x_i=x_j>0$. On the circle
$y_i^2+y_j^2=x_i^2+x_j^2$ with $y_k=x_k$ for $k\notin\{i,j\}$,
$v^{(t)}$ attains a maximum at $y=x$. Indeed, \eqref{eq:angular-sign}
shows that $v^{(t)}$ decreases as $\arctan(y_j/y_i)$ increases from
$\pi/4$ toward $\pi/2$. Reflections and interchanges of $y_i,y_j$
send the arc $0\le y_i\le y_j$ to all eight arcs of the circle.
By \eqref{eq:least-pair}, $\Delta_{ij}v^{(t)}>0$. The Hopf lemma for
$v^{(t)}$ on $y_i^2+y_j^2<x_i^2+x_j^2$, $y_k=x_k$ for $k\notin\{i,j\}$,
yields at $y=x$
\begin{equation}\label{eq:tied-coordinate-sign}
0<\frac{\partial_iv^{(t)}(x)+\partial_jv^{(t)}(x)}{\sqrt2}
=\sqrt2\,\partial_iv^{(t)}(x)=\sqrt2\,\partial_jv^{(t)}(x),
\end{equation}
where we used $\partial_iv^{(t)}(x)=\partial_jv^{(t)}(x)$ by symmetry.

Finally, $\partial_i v^{(t)}$ satisfies
\[
\begin{cases}
\tr\bigl((I+(D^2v^{(t)})^2)^{-1}D^2(\partial_iv^{(t)})\bigr)=0
&\text{in }B_1\cap\{0<x_i<\max_{k\ne i}|x_k|\},\\
\partial_iv^{(t)}=0&\text{on }B_1\cap\{x_i=0\},\\
\partial_iv^{(t)}>0
&\text{on }\partial B_1\cap\{0<x_i\le\max_{k\ne i}|x_k|\},\\
\partial_iv^{(t)}>0&\text{on }B_1\cap\{x_i=|x_j|>0\},\quad\forall j\ne i,
\end{cases}
\]
where we differentiated \eqref{eq:dirichlet} with respect to $x_i$ in the
first line, used reflection across $x_i=0$ in the second,
used $v^{(t)}=tU_t$ and \eqref{eq:spherical-coordinate-sign} in the third,
and combined \eqref{eq:tied-coordinate-sign} with reflection across
$x_j=0$ in the fourth.
On the boundary of each connected component of
$B_1\cap\{0<x_i<\max_{j\ne i}|x_j|\}$, we have
$\partial_i v^{(t)}\ge0$, with $\partial_i v^{(t)}>0$ at some boundary point
on $\mathbb S^{n-1}$. By the strong maximum principle and reflection in the
coordinate hyperplanes, we obtain
\begin{equation}\label{eq:signs}
x_i\partial_iv^{(t)}(x)>0
\quad\text{if }x_i\ne0,\quad
|x_i|\le\max_{j\ne i}|x_j|.
\end{equation}
By \eqref{eq:signs}, we have $|x_i|>\max_{j\ne i}|x_j|$ whenever
$x_ip_i<0$ at any point
\[
(x,p)=(\lambda y,\lambda Du_t(y)),\qquad
\lambda>0,\quad y\in B_{1/r_t}.
\]
In particular,
\begin{equation}\label{eq:closed-coordinate-signs}
|x_i|=\max_k|x_k|\quad\text{if }\,\,
x_ip_i<0,
\quad \#\{k:x_kp_k<0\}\le1.
\end{equation}

Let
\[
\widetilde v_k(x)\coloneqq R_k^{-2}u_{t_k}(R_kx)
\qquad\forall x\in B_{1/(r_{t_k}R_k)}.
\]
Then $\widetilde S_k\coloneqq(R_k^{-1}\operatorname{Id})_\sharp\Graph(Du_{t_k})
=\Graph(D\widetilde v_k)$.
Since $r_{t_k}R_kR/t_k\to0$ as $k\to\infty$ for fixed $R>0$,
we obtain from \eqref{eq:first-area} at radius $R_kR$
\[
\begin{aligned}
\|T\|(\mathbf B_R)
&\le\liminf_{k\to\infty}\|\widetilde S_k\|(\mathbf B_R) =\liminf_{k\to\infty}R_k^{-n}
\|\Graph(Du_{t_k})\|(\mathbf B_{R_kR})
\le(2n-1)\omega_nR^n.
\end{aligned}
\]
Since $r_{t_k}R_k/t_k\to0$ as $k\to\infty$, \eqref{eq:first-boundary} implies
\[
\dist(0,\spt\partial\widetilde S_k)
=R_k^{-1}\dist(0,\spt\partial\Graph(Du_{t_k}))
\longrightarrow\infty\quad\text{as }k\to\infty.
\]
Hence $T$ is locally integral and $\partial T=0$ by
\cite[Theorem 27.3]{Simon}.
Local current convergence preserves the calibration of $\widetilde S_k$,
so $T$ has phase $(n-2)\pi/2$.
Applying \eqref{eq:gauss-plane-convergence} to the oriented tangent
$n$-vectors of $\widetilde S_k$ and $T$ near each regular point of $\spt T$,
we obtain
\[
T_z(\spt T)\in\Gr_n((n-2)\pi/2)
\qquad\forall z\in\spt T\setminus\Sing T.
\]
Thus $T$ is a sum of calibrated planes by
Proposition~\ref{prop:planar}.
Since $\widetilde S_k\rightharpoonup T$ locally as currents as $k\to\infty$,
\[
\dist(z,\spt\widetilde S_k)\longrightarrow0
\quad\text{as }k\to\infty,\qquad\forall z\in\spt T.
\]
For $(x,p)\in\spt T$ with $x_ip_i<0$, choose
$(x_k,p_k)\in\spt\widetilde S_k$ tending to $(x,p)$. Then
$(x_k)_i(p_k)_i<0$ for large $k$, so \eqref{eq:closed-coordinate-signs}
implies $|(x_k)_i|=\max_j|(x_k)_j|$. Letting $k\to\infty$, we obtain
$|x_i|=\max_j|x_j|$.
Now fix a plane $\Pi\ne\{x=0\}$ in the decomposition of $T$ and set
\[
W\coloneqq\{w\in\R^n:(0,w)\in\Pi\}.
\]
Since $\Pi$ is Lagrangian, $W=\bigl(\pi_x(\Pi)\bigr)^\perp$.
By \eqref{eq:closed-principal-angles}, all vertical principal angles of
$\Pi$ equal $\pi/2$. If $\dim \pi_x(\Pi)=1$, the single finite angle equals
$-\pi/2$, a contradiction. Hence $\dim \pi_x(\Pi)\ge2$.
Suppose that some $w\in W$ has $w_i\ne0$ and that $x_i$ is not
identically zero on $\pi_x(\Pi)$.
For $(x,p)\in\Pi$ with $x_i\ne0$, we have $(x,p+\tau w)\in\Pi$ and
$x_i(p_i+\tau w_i)<0$ for some $\tau\in\R$.
By \eqref{eq:closed-coordinate-signs} and continuity,
\[
\max_j|x_j|=|x_i|\qquad\forall x\in \pi_x(\Pi).
\]
Hence $\pi_x(\Pi)\cap\{x_i=0\}=\{0\}$, contradicting $\dim \pi_x(\Pi)\ge2$.
Thus $w_i\ne0$ forces $\pi_x(\Pi)\subset\{x_i=0\}$ and hence $e_i\in W$.
Therefore $W$ is a coordinate subspace.

By \eqref{eq:one-sheet}, at most one plane with $\dim W=0$ occurs in $T$,
with multiplicity one.
When $\dim W\ge1$, $\#(G\cdot W)
\ge n$. If two distinct planes in the decomposition of $T$ have vertical
subspace $W$, or the coefficient of $\Pi$ is at least $2$, then, by $G$-symmetry,
\[
2n\le2\#(G\cdot W)
\le\frac{\|T\|(\mathbf B_1)}{\omega_n}\le2n-1,
\]
a contradiction.
Hence $\Pi$ is unique over $W$ and has coefficient $1$, and so $\Pi=\{(x,Ax+w):x\in \pi_x(\Pi),\ w\in W\}$ where $A:\pi_x(\Pi)\to \pi_x(\Pi)$ is symmetric.
By uniqueness and $G$-symmetry, $OAO^T=A$ whenever $O\in G$ and $O(W)=W$.
In the coordinate basis of $\pi_x(\Pi)$, we obtain $A_{ij}=-A_{ij}=0$
for $i\ne j$ by coordinate reflections, and $A_{ii}=A_{jj}$ by
coordinate permutations, so $A=A_{11}I_{n-\dim W}$.
By \eqref{eq:closed-principal-angles},
\[
(n-\dim W)\arctan A_{11}+\frac\pi2\dim W=\frac{n-2}{2}\pi.
\]
Hence
\[
A_{11}=\tan\Bigl(\frac{n-\dim W-2}{2(n-\dim W)}\pi\Bigr)
=\cot\frac\pi{n-\dim W}.
\]
If $n\ge5$ and $2\le\dim W\le n-2$, then, by $G$-symmetry,
\[
2n>\frac{\|T\|(\mathbf B_1)}{\omega_n}
\ge\#(G\cdot W)=\binom n{\dim W}\ge\binom n2\ge2n,
\]
a contradiction. Thus, for $\Pi\ne\{x=0\}$,
\[
A=
\begin{cases}
\cot(\pi/n)I_n,&\dim W=0,\\
\cot(\pi/(n-1))I_{n-1},&\dim W=1,\\
0,&n=4,\ \dim W=2.
\end{cases}
\]
If $\dim W\le1$,
\[
\llbracket\Pi\rrbracket=
\begin{cases}
\Iplane,&W=\{0\},\\
\Pplane_i,&W=\R e_i,\quad1\le i\le n.
\end{cases}
\]
If $n=4$ and $\dim W=2$,
$\#(G\cdot W)
=6$.
We orient each plane $\Pi\ne\{x=0\}$ consistently with $T$, and
$\{0\}\times\R^n$ by $dp_1\wedge\cdots\wedge dp_n$.
Since $\operatorname{Re}(e^{-i(n-2)\pi/2}i^n)=-1$,  
\begin{equation}\label{eq:plane-list}
\begin{aligned}
T&=m_0\Iplane+m_1\sum_{i=1}^n\Pplane_i-m_n\llbracket\{0\}\times\R^n\rrbracket +m_2\sum_{i<j}
\llbracket\{x_i=x_j=0,\ p_k=0\text{ for }k\notin\{i,j\}\}\rrbracket,
\end{aligned}
\end{equation}
where $m_0\in\{0,1\}$ and $m_1,m_2,m_n\in\mathbb N_0$, with $m_2=0$ for
$n\ge5$ and $m_0+nm_1+m_n+6m_2\le2n-1$.
Moreover, we claim that 
\begin{equation}\label{eq:oblique-congruence}
m_0+nm_1+6m_2-m_n\equiv1\pmod{2n}.
\end{equation}
Indeed, let $\pi_+(x,p)\coloneqq x+p$.
On every plane in \eqref{eq:plane-list}, we have $x\cdot p\ge0$, so
\begin{equation}\label{eq:cone-oblique-separation}
|\pi_+(z)|^2=|x|^2+|p|^2+2x\cdot p\ge|z|^2
\qquad\forall z=(x,p)\in\spt T.
\end{equation}
Since $\spt\widetilde S_k\to\spt T$ locally in Hausdorff distance and
$r_{t_k}R_k\to0$ as $k\to\infty$, we may assume that
$\overline B_1\subset B_{1/(r_{t_k}R_k)}$ and
\begin{equation}\label{eq:oblique-separation}
|\pi_+(z)|\ge\tfrac12
\qquad\forall z\in\spt\widetilde S_k\text{ with }1\le|z|\le2.
\end{equation}
If $(x,p)\in\spt\widetilde S_k$ satisfies $|x|\le1$ and $|x+p|<1/2$, then
\[
|(x,p)|^2\le |x|^2+(|x+p|+|x|)^2\le 1+(3/2)^2<4.
\]
Hence $|(x,p)|<1$ by \eqref{eq:oblique-separation}, so
$(\operatorname{Id}+D\widetilde v_k)(\partial B_1)\cap B_{1/2}=\varnothing$.

If $x\in B_1$ and $x+D\widetilde v_k(x)=0$, then
$x_i\partial_i\widetilde v_k(x)=-x_i^2$. By \eqref{eq:closed-coordinate-signs},
\[
\#\{i:x_i\ne0\}=\#\{i:x_i\partial_i\widetilde v_k(x)<0\}\le1.
\]
Since $D\widetilde v_k(0)=0$ and
$I+D^2\widetilde v_k(0)=(1+\cot(\pi/n))I>0$, choose $0<\rho<1$
so that $\operatorname{Id}+D\widetilde v_k$ has only the zero $0$ in
$\overline B_\rho$, of local degree $1$. Each nonzero root $x\in B_1$ satisfies
$x\in D_i^\pm$ for exactly one of the disjoint open sets
\[
D_i^\pm\coloneqq\{x\in B_1:\ \pm x_i>\rho,\ |x_j|<\rho
\text{ for all }j\ne i\},\qquad i=1,\ldots,n.
\]
Since $Ox+D\widetilde v_k(Ox)=O(x+D\widetilde v_k(x))$ for $O\in G$,
the degrees $\degmap(\operatorname{Id}+D\widetilde v_k,D_i^\pm,0)$ are equal by the transitive action of $G$
and $(\det O)^2=1$.
By additivity,
\begin{equation}\label{eq:finite-oblique-degree}
\degmap(\operatorname{Id}+D\widetilde v_k,B_1,0)
=1+2n\degmap(\operatorname{Id}+D\widetilde v_k,D_1^+,0)
\in1+2n\mathbb Z.
\end{equation}

To transfer \eqref{eq:finite-oblique-degree} to $T$, take a smooth
$n$-form $\omega$ supported in $B_{1/2}$ with
$\int\omega=1$, and let $0\le\chi\le1$ with
$\chi\in C_c^\infty(\mathbf B_2)$ and $\chi=1$ on $\mathbf B_1$.
By \eqref{eq:oblique-separation}, on $\widetilde S_k$ the integrand
$\chi\pi_+^*\omega$ agrees with $\pi_+^*\omega$ for $x\in B_1$ and
$\chi\pi_+^*\omega=0$ for $|x|\ge1$. On $\spt T$, $\chi=1$ wherever $\pi_+(z)\in B_{1/2}$
by \eqref{eq:cone-oblique-separation}. Since
$(\operatorname{Id}+D\widetilde v_k)(\partial B_1)\cap B_{1/2}=\varnothing$,
$\degmap(\operatorname{Id}+D\widetilde v_k,B_1,p)$ is constant for
$p\in B_{1/2}$. Thus the degree formula with $\int\omega=1$ and local
current convergence imply
\begin{equation}\label{eq:oblique-degree-limit}
\begin{aligned}
\degmap(\operatorname{Id}+D\widetilde v_k,B_1,0)
&=\widetilde S_k(\chi\pi_+^*\omega)\\
&\longrightarrow T(\chi\pi_+^*\omega)
=m_0+nm_1+6m_2-m_n\qquad\text{as }k\to\infty.
\end{aligned}
\end{equation}
The value of $T(\chi\pi_+^*\omega)$ in \eqref{eq:oblique-degree-limit}
follows from the degrees $1$, $1$, $-1$ and $1$ of $\pi_+$ on
$\Iplane$, $\Pplane_i$, $-\llbracket\{0\}\times\R^n\rrbracket$ and
the coordinate two-vertical planes, respectively.
By \eqref{eq:finite-oblique-degree} and \eqref{eq:oblique-degree-limit},
\eqref{eq:oblique-congruence} holds.
If $n=4$ and $m_2\ge1$, then, by $m_0+4m_1+6m_2+m_n\le7$,
$m_2=1$, $m_1=0$ and $m_0+m_n\le1$. Hence
$m_0+4m_1+6m_2-m_n\in\{5,6,7\}$, contradicting
\eqref{eq:oblique-congruence}. Thus $m_2=0$, $\forall n\ge4$.
Now $m_1\le1$. If $m_1=0$, combining $m_0+nm_1+m_n\le2n-1$ with
\eqref{eq:oblique-congruence}, we obtain
$(m_0,m_n)=(1,0)$ or $(0,2n-1)$. If  $m_1=1$, we obtain
$m_0=0$ and $m_n=n-1$.
Substituting $(m_0,m_1,m_n)$ into \eqref{eq:plane-list}, we obtain
\eqref{eq:two-models}.
\end{proof}

Now we prove Proposition~\ref{prop:noescape}.
\begin{proof}[Proof of Proposition~\ref{prop:noescape}]
\textup{(i)} Let $S$ be a local current limit of $\Graph(Du_t)$.
For a cone $T$ at infinity of $S$,
there exists a sequence $\{R_k\}$ with $R_k\to\infty$ such that
$(R_k^{-1}\operatorname{Id})_\sharp S\rightharpoonup T$ as $k\to\infty$.
Since $r_t\to0$, we can choose a sequence $\{t_k\}$ with $t_k\to\infty$
and $\Graph(Du_{t_k})\rightharpoonup S$ such that
\[
(R_k^{-1}\operatorname{Id})_\sharp \Graph(Du_{t_k})\rightharpoonup T,\qquad
r_{t_k}R_k\longrightarrow0\quad\text{as }k\to\infty.
\]
If $T=\Iplane$, the density of $S$ at infinity is $1$. The density at
$(0,0)$ is also $1$ by \eqref{eq:graph-near-origin}, so monotonicity implies
\[
1\le\frac{\|S\|(\mathbf B_R)}{\omega_nR^n}\le1
\qquad\forall R>0.
\]
By the monotonicity formula,
\begin{equation}\label{eq:annular-monotonicity}
\int_{\mathbf B_R\setminus\mathbf B_r}
\frac{\dist(z,T_z(\spt S))^2}{|z|^{n+2}}\,d\|S\|(z)
=\frac{\|S\|(\mathbf B_R)}{R^n}
-\frac{\|S\|(\mathbf B_r)}{r^n}=0
\qquad\forall\,0<r<R.
\end{equation}
Hence $S=T=\Iplane$ and $\pi_x(\spt S)=\R^n$.
By Proposition~\ref{prop:first-limit}\textup{(ii)}, up to
a subsequence,
\[
u_{t_k}\longrightarrow\frac{\cot(\pi/n)}2|x|^2
\quad\text{in }C^\infty_{\mathrm{loc}}(\R^n)\quad\text{as }k\to\infty,
\]
since $u_{t_k}(0)=0$ and $\Iplane$ is the graph of the gradient of
$\cot(\pi/n)|x|^2/2$.
Consequently,
\[
\cot(\pi/n)=\lim_{k\to\infty}\osc_{B_1}u_{t_k}
=\osc_{B_1}\Bigl(\frac{\cot(\pi/n)}2|x|^2\Bigr)=\frac{\cot(\pi/n)}2,
\]
which contradicts $\cot(\pi/n)>0$. Therefore, by Lemma~\ref{lem:cone-models},
\begin{equation}\label{eq:limit-density-at-infinity}
\lim_{R\to\infty}\frac{\|S\|(\mathbf B_R)}{\omega_nR^n}=2n-1.
\end{equation}

Next, we prove that, for each $P_*>0$, there exist $R_{P_*},t_{P_*}<\infty$ such that
\begin{equation}\label{eq:inverse-bound}
|y|\le R_{P_*}\qquad\forall t\ge t_{P_*},\quad
y\in B_{1/r_t}\text{ with }|Du_t(y)|\le {P_*}.
\end{equation}
If \eqref{eq:inverse-bound} fails, we can choose $t_j\to\infty$
and $y_j\in B_{1/r_{t_j}}$ with $|Du_{t_j}(y_j)|\le P_*$ and
$|y_j|\to\infty$. By Proposition~\ref{prop:first-limit},
up to a subsequence, $\Graph(Du_{t_j})\rightharpoonup S^\prime$ locally
in the sense of currents as $j\to\infty$.
Let $\widehat S_j\coloneqq (|y_j|^{-1}\operatorname{Id})_\sharp \Graph(Du_{t_j})$.
Since $y_j\in B_{1/r_{t_j}}$, we have $r_{t_j}|y_j|<1$. For each fixed
$r>0$, $r_{t_j}|y_j|r/t_j\le r/t_j\to0$ as $j\to\infty$.
Taking $R=|y_j|r$ in \eqref{eq:first-area}, we obtain
\[
\begin{aligned}
\limsup_{j\to\infty}\frac{\|\widehat S_j\|(\mathbf B_r)}{\omega_nr^n}
&=\limsup_{j\to\infty}\frac{\|\Graph(Du_{t_j})\|(\mathbf B_{|y_j|r})}
        {\omega_n(|y_j|r)^n}\le2n-1.
\end{aligned}
\]
On the other hand, by \eqref{eq:first-boundary} and
$r_{t_j}|y_j|<1$,
\[
\liminf_{j\to\infty}\frac{\dist(0,\spt\partial\widehat S_j)}{t_j}
\ge\liminf_{j\to\infty}\frac{r_{t_j}}{t_j}
\dist(0,\spt\partial\Graph(Du_{t_j}))>0.
\]
Thus $\dist(0,\spt\partial\widehat S_j)\to\infty$ as $j\to\infty$.
By current compactness \cite[Theorem 27.3]{Simon}, up to a subsequence,
$\widehat S_j\rightharpoonup T$ locally in the sense of currents as
$j\to\infty$, with $\partial T=0$.
All $\widehat S_j$ have phase $(n-2)\pi/2$, and so does $T$.

On the other hand, let ${R^\prime}>0$ satisfy $\|S^\prime\|(\partial\mathbf B_{R^\prime})=0$.
Since ${R^\prime}/(|y_j|r)\to0$ as $j\to\infty$, we obtain by monotonicity
\begin{equation}\label{eq:moving-comparison}
\begin{aligned}
\liminf_{j\to\infty}\frac{\|\widehat S_j\|(\mathbf B_r)}{\omega_nr^n}
&\ge\lim_{j\to\infty}\frac{\|\Graph(Du_{t_j})\|(\mathbf B_{R^\prime})}{\omega_n{R^\prime}^n} =\frac{\|S^\prime\|(\mathbf B_{R^\prime})}{\omega_n{R^\prime}^n}.
\end{aligned}
\end{equation}
Letting ${R^\prime}\to\infty$ in
\eqref{eq:moving-comparison}, and using \eqref{eq:first-area} and
\eqref{eq:limit-density-at-infinity} for $S^\prime$, we have
$2n-1\le\|T\|(\mathbf B_r)/(\omega_nr^n)\le2n-1$
whenever $r>0$ and $\|T\|(\partial\mathbf B_r)=0$.
Moreover, up to a subsequence,
$|y_j|^{-1}\bigl(y_j,Du_{t_j}(y_j)\bigr)\to(a,0)$ with $|a|=1$.
By monotonicity at $|y_j|^{-1}\bigl(y_j,Du_{t_j}(y_j)\bigr)$,
$(a,0)\in\spt T$.
Since $\mathbf B_s\uparrow\mathbf B_r$ as $s\uparrow r$ and
$\|T\|(\partial\mathbf B_s)=0$ for a.e. $s>0$, we obtain
\begin{equation}\label{eq:saturation}
\|T\|(\mathbf B_r)
=\lim_{s\uparrow r,\ \|T\|(\partial\mathbf B_s)=0}\|T\|(\mathbf B_s)
=(2n-1)\omega_n\lim_{s\uparrow r}s^n
=(2n-1)\omega_nr^n\quad\forall r>0.
\end{equation}
For $0<r<R$, applying the monotonicity formula to $T$ and substituting
\eqref{eq:saturation}, we obtain
\[
\int_{\mathbf B_R\setminus\mathbf B_r}
\frac{\dist(z,T_z(\spt T))^2}{|z|^{n+2}}\,d\|T\|(z)
=\frac{\|T\|(\mathbf B_R)}{R^n}-\frac{\|T\|(\mathbf B_r)}{r^n}=0.
\]
Thus $z\in T_z(\spt T)$ for $\|T\|$-a.e. $z$, and $T$ is a cone
centered at $(0,0)$.

After passing to a subsequence, either $\inf_j(r_{t_j}|y_j|)>0$ or
$r_{t_j}|y_j|\to0$. In the first case,
\[
\pi_x(\spt T)\subset\overline B_{1/\inf_j(r_{t_j}|y_j|)}.
\]
Since $T$ is invariant under positive dilations, boundedness of
$\pi_x(\spt T)$ implies $\pi_x(\spt T)=\{0\}$, contradicting
$(a,0)\in\spt T$ with $|a|=1$.
In the second case, $r_{t_j}|y_j|\to0$, so the convergence
$\widehat S_j\rightharpoonup T$ satisfies the hypotheses of
Lemma~\ref{lem:cone-models}. By \eqref{eq:saturation}, $T\ne\Iplane$.
For either remaining cone in \eqref{eq:two-models},
$\spt T\cap\{p=0\}=\{0\}$, again contradicting $(a,0)\in\spt T$ with $|a|=1$.
Hence \eqref{eq:inverse-bound} holds.

Let $(x,p)\in\spt S$ with $|p|\le P_*$. By local Hausdorff convergence
of the supports, choose $(x_k,Du_{t_k}(x_k))\to(x,p)$ as $k\to\infty$.
Since $|Du_{t_k}(x_k)|<P_*+1$ for large $k$, \eqref{eq:inverse-bound}
at radius $P_*+1$ implies
\[
|x|=\lim_{k\to\infty}|x_k|\le R_{P_*+1}.
\]
Since $\spt S$ is closed, $\pi_p:\spt S\to\R^n$ is proper.

Consequently, for any sequence $\{\Graph(Du_{t_j})\}$ converging locally as
currents to $S$ with $t_j\to\infty$, and any smooth
compactly supported $n$-form $\omega$ in the $p$-variables, we choose,
by \eqref{eq:inverse-bound} and the bounded $x$-projections of the finitely
many initial graphs, $\chi\in C_c^\infty(\C^n)$ with
$0\le\chi\le1$ and $\chi=1$ on
$\spt \Graph(Du_{t_j})\cap\pi_p^{-1}(\spt\omega)$ for all $j$ and on
$\spt S\cap\pi_p^{-1}(\spt\omega)$.
By \eqref{eq:gradient-rescaling}, $Du_{t_j}(y)=p$ is equivalent to
$DU_{t_j}(r_{t_j}y)=\frac{r_{t_j}}{t_j}p$.
Since $r_{t_j}/t_j\to0$ as $j\to\infty$, we have
\[
\sup_{p\in\spt\omega}\frac{r_{t_j}}{t_j}|p|
\longrightarrow0\quad\text{as }j\to\infty.
\]
Hence, by current convergence and \eqref{eq:degree},
\[
S(\chi\pi_p^*\omega)
=\lim_{j\to\infty}\Graph(Du_{t_j})(\chi\pi_p^*\omega)
=-(2n-1)\int_{\R^n}\omega.
\]
Therefore
\begin{equation}\label{eq:complete-degree}
(\pi_p)_\sharp S=-(2n-1)\llbracket\R^n\rrbracket.
\end{equation}

\textup{(ii)} Let $T$ be a cone at infinity of $S$, and let $R_j\to\infty$
satisfy $(R_j^{-1}\operatorname{Id})_\sharp S\rightharpoonup T$ locally as currents
as $j\to\infty$.
By Lemma~\ref{lem:cone-models} and \eqref{eq:limit-density-at-infinity},
\begin{equation}\label{eq:cone-projection-estimate}
|(x,p)|^2=|x|^2+|p|^2
\le\Bigl(1+\tan^2\frac\pi{n-1}\Bigr)|p|^2
\le4|p|^2,\quad\forall(x,p)\in\spt T.
\end{equation}
We claim that
\begin{equation}\label{eq:momentum-proper}
|(x,p)|\le3|p|\qquad
\text{on }\spt S\text{ outside a compact set}.
\end{equation}
If \eqref{eq:momentum-proper} fails, there exist $z_k\in\spt S$ with
$|z_k|\to\infty$ as $k\to\infty$ and $|z_k|>3|\pi_p(z_k)|$.
Up to a subsequence, $(|z_k|^{-1}\operatorname{Id})_\sharp S\rightharpoonup T'$
locally as currents as $k\to\infty$, where $T'$ is a cone at infinity of $S$,
and $z_k/|z_k|\to z\in\spt T'$ as $k\to\infty$.
Applying \eqref{eq:cone-projection-estimate} to $T'$, we obtain
\[
1=|z|\le2|\pi_p(z)|
=2\lim_{k\to\infty}\frac{|\pi_p(z_k)|}{|z_k|}
\le\frac23<1,
\]
a contradiction.

To compute the $p$-projection degree of $T$, fix a smooth compactly
supported $n$-form $\omega$ on $\R_p^n$ and choose $C(\omega)>0$ with
$\spt\omega\subset B_{C(\omega)}$.
By \eqref{eq:momentum-proper}, 
$|z|\le C(n,S)+3|\pi_p(z)|$, $\forall z\in\spt S$, for some $C(n,S)>0$.
For $R>1$ and $z\in\spt((R^{-1}\operatorname{Id})_\sharp S)$ with
$\pi_p(z)\in\spt\omega$, we have $Rz\in\spt S$, and hence
\[
|z|\le\frac{C(n,S)}R+3|\pi_p(z)|
\le C(n,S)+3C(\omega).
\]
On the other hand, by \eqref{eq:cone-projection-estimate},
\[
|z|\le2|\pi_p(z)|\le2C(\omega)
\qquad\forall z\in\spt T\text{ with }\pi_p(z)\in\spt\omega.
\]
Choose $\chi\in C_c^\infty(\C^n)$ with
$0\le\chi\le1$ and $\chi=1$ on $\mathbf B_{C(n,S)+3C(\omega)+1}$.
For $z\in\spt T\cup\spt((R^{-1}\operatorname{Id})_\sharp S)$ with $R>1$,
$\pi_p(z)\in\spt\omega$ implies $|z|\le C(n,S)+3C(\omega)$.
Hence $\chi\pi_p^*\omega=\pi_p^*\omega$ on both supports. Since
$\pi_p\circ(R_j^{-1}\operatorname{Id})=R_j^{-1}\pi_p$,
local current convergence and \eqref{eq:complete-degree} imply
\[
\begin{aligned}
\bigl((\pi_p)_\sharp T\bigr)(\omega)
&=T(\chi\pi_p^*\omega)
=\lim_{j\to\infty}\bigl((R_j^{-1}\operatorname{Id})_\sharp S\bigr)(\chi\pi_p^*\omega)\\
&=\lim_{j\to\infty}
\bigl((R_j^{-1}\operatorname{Id}_{\R^n})_\sharp
      ((\pi_p)_\sharp S)\bigr)(\omega)
=-(2n-1)\int_{\R^n}\omega,
\end{aligned}
\]
Hence $(\pi_p)_\sharp T=-(2n-1)\llbracket\R^n\rrbracket$.
On the other hand, using $(\pi_p)_\sharp\llbracket\{0\}\times\R^n\rrbracket=\llbracket\R^n\rrbracket$ and
$(\pi_p)_\sharp\Pplane_i=\llbracket\R^n\rrbracket$, we obtain
\[
\begin{aligned}
(\pi_p)_\sharp\Bigl(-(n-1)\llbracket\{0\}\times\R^n\rrbracket+\sum_{i=1}^n\Pplane_i\Bigr)
&=\bigl(-(n-1)+n\bigr)\llbracket\R^n\rrbracket =\llbracket\R^n\rrbracket\ne-(2n-1)\llbracket\R^n\rrbracket.
\end{aligned}
\]
Since $\|T\|(\mathbf B_1)=(2n-1)\omega_n>\omega_n
=\|\Iplane\|(\mathbf B_1)$ by \eqref{eq:limit-density-at-infinity},
we have $T\ne\Iplane$.
Hence $T=-(2n-1)\llbracket\{0\}\times\R^n\rrbracket$ by \eqref{eq:two-models}.
\end{proof}

\section{Proof of Theorem~\ref{thm:main}}\label{sec:profile}
\label{sec:completion}
\label{sec:nozero}
\label{subsec:approximation}
\label{subsec:jets}
\label{subsec:support-transfer}

Now we prove Theorem~\ref{thm:main}.
\begin{proof}[Proof of Theorem~\ref{thm:main}]
Let $t_j\to\infty$. By Proposition~\ref{prop:first-limit}, up to a
subsequence, $\Graph(Du_{t_j})\rightharpoonup S$ and
$u_{t_j}\to u$ in $C^\infty_{\mathrm{loc}}(\pi_x(\Lambda))$ as $j\to\infty$.

\emph{Step 1.}\textit{ We construct a Dirichlet minimizer $F$ on $\R^n$ such that
each $F(p)$ consists of $2n-1$ vectors counted with multiplicity,
and $0\notin\spt F(p)$ for $p\ne0$.}

\emph{Step 1.1.} \textit{We construct $F$ with $F\not\equiv0$ and
$F(0)=(2n-1)\llbracket0\rrbracket$.}

Let
\begin{equation}\label{eq:height-function}
	\rho_x(R)\coloneqq \max\{|x|:(x,p)\in\spt S,\ |p|\le R\}
	\qquad\forall R>0.
\end{equation}
Since $\pi_p$ is proper by \eqref{eq:momentum-proper},
$\spt S\cap\{|p|\le R\}$ is compact. Hence the maximum in
\eqref{eq:height-function} is attained and $\rho_x$ is upper semicontinuous.
Moreover, $\pi_x(\Lambda)$ is unbounded and $\Graph(Du)\subset\spt S$ by
Proposition~\ref{prop:first-limit}, so $\rho_x(R)\to\infty$ as $R\to\infty$.
By \eqref{eq:momentum-proper},
\[
\limsup_{R\to\infty}\sup_{(x,p)\in\spt S,\ |p|\le R}
\frac{|(x,p)|}R\le3<4.
\]
By local convergence of the supports of $(R^{-1}\operatorname{Id})_\sharp S$
to $\{x=0\}$ in Proposition~\ref{prop:noescape}\textup{(ii)},
\begin{equation}\label{eq:sublinear-x-height}
0\le\limsup_{R\to\infty}\frac{\rho_x(R)}R
\le\limsup_{R\to\infty}\sup\bigl\{|x|:(x,p)\in\spt((R^{-1}\operatorname{Id})_\sharp S)
\cap\overline{\mathbf B}_4\bigr\}=0.
\end{equation}
For the rescaling, choose $\underline R_j\to\infty$ with $\rho_x(\underline R_j)>0$ and
\begin{equation}\label{eq:profile-rescaling-scale}
\sup_{R\ge\underline R_j}\frac{\rho_x(R)}R
\le\min\Bigl\{\frac18,\frac{2^{-j/2}}{j^{(n+2)/2}}\Bigr\}.
\end{equation}
As $\rho_x(R)/R\to0$ when $R\to\infty$, we can choose $L_j>\underline R_j$ with
\[
\sup_{R\ge L_j}\frac{\rho_x(R)}R
<\frac{\rho_x(\underline R_j)}{\underline R_j}.
\]
By upper semicontinuity of $\rho_x(R)/R$ on $[\underline R_j,L_j]$,
there exists $R_j\in[\underline R_j,L_j]$ such that
\[
\frac{\rho_x(R_j)}{R_j}
=\max_{\underline R_j\le R\le L_j}\frac{\rho_x(R)}R
=\max_{R\ge\underline R_j}\frac{\rho_x(R)}R.
\]
Since $R_j\ge\underline R_j\to\infty$, $\rho_x(R)\to\infty$ and
$\rho_x(R)/R\to0$ as $R\to\infty$, we have
\[
\rho_x(R_j)\to\infty,\qquad \frac{\rho_x(R_j)}{R_j}\to0\quad\text{as }j\to\infty.
\]
By maximality of $\rho_x(R_j)/R_j$ on $[\underline R_j,\infty)$,
\begin{equation}\label{eq:height-selection}
\rho_x(LR_j)\le LR_j\frac{\rho_x(R_j)}{R_j}
=L\rho_x(R_j),\qquad\forall L\ge1.
\end{equation}

Next, let $T_j$ be the current in $\R_p^n\times\R_x^n$ defined by $T_j\coloneqq-\bigl(R_j^{-1}\pi_p,R_j^{-1}\pi_x\bigr)_\sharp S$. 
Since $\rho_x(R_j)/R_j\le1/8$, \eqref{eq:height-selection} yields
\[
\sup_{\spt T_j\cap\pi_p^{-1}(B_{2L})}|x|
\le\frac{\rho_x(2LR_j)}{R_j}\le2L\frac{\rho_x(R_j)}{R_j}\le L
\qquad\forall L\ge1.
\]
Since $|p|<2L$ on $\pi_p^{-1}(B_{2L})$,
$\spt T_j\cap\pi_p^{-1}(B_{2L})\subset\mathbf B_{3L}$. Therefore
\[
\|T_j\|(\pi_p^{-1}(B_{2L}))\le\|T_j\|(\mathbf B_{3L})
\le(2n-1)\omega_n(3L)^n.
\]
By \eqref{eq:complete-degree},
$(\pi_p)_\sharp T_j=(2n-1)\llbracket\R^n\rrbracket$.
As $S$ has phase $(n-2)\pi/2$, $T_j$ is calibrated by
$\operatorname{Re}\bigl((dp_1+i\,dx_1)\wedge\cdots\wedge(dp_n+i\,dx_n)\bigr)$,
so by \eqref{eq:calibrated-height-jacobian},
$0\le1-J_p\le C(n)|\nabla_{T_j}x|^2$ at $\|T_j\|$-a.e. point.
Since $\partial T_j=0$ and $\pi_p$ is proper on $\spt T_j$,
combining \eqref{eq:bounded-height-energy} and \eqref{eq:height-selection}, we obtain
\begin{equation}\label{eq:profile-excess}
\begin{aligned}
\frac{\|T_j\|(\pi_p^{-1}(B_L))}{\omega_nL^n}-(2n-1)
&=\frac1{\omega_nL^n}\int_{\pi_p^{-1}(B_L)}(1-J_p)\,d\|T_j\| \le\frac{C(n)}{L^n}
\int_{\pi_p^{-1}(B_L)}|\nabla_{T_j}x|^2\,d\|T_j\|\\
&\le\frac{C(n)}{L^{n+2}}
\sup_{\pi_p^{-1}(B_{2L})\cap\spt T_j}|x|^2
\|T_j\|(\pi_p^{-1}(B_{2L}))
\le C(n)\frac{\rho_x(R_j)^2}{R_j^2}\quad\forall L\ge1.
\end{aligned}
\end{equation}

Choose $r\in(2j,3j)$ so that $T_j\res\pi_p^{-1}(B_{4r})$ is integral.
By properness of $\pi_p$, $\partial T_j=0$, and
$(\pi_p)_\sharp T_j=(2n-1)\llbracket\R^n\rrbracket$, the current
$(r^{-1}\operatorname{Id})_\sharp\bigl(T_j\res\pi_p^{-1}(B_{4r})\bigr)$ satisfies
\cite[Assumption 1.1]{DS}. By \eqref{eq:profile-excess} at $L=4r$, we
take $\underline R_j$ large enough in \eqref{eq:profile-rescaling-scale}
for \cite[Theorems 1.4 and 1.6]{DS} with $\bar\eta=j^{-(n+4)}$.
By (1.5), the graph current of the Lipschitz map $\frac{\rho_x(R_j)}{R_j}F_j$
agrees with $T_j$ over a closed set $K_j\subset B_r$. In particular,
\begin{equation}\label{eq:profile-slice-identity}
\sum_{\mathbf a\in F_j(p)}\llbracket\Bigl(p,\frac{\rho_x(R_j)}{R_j}\mathbf a\Bigr)\rrbracket
=\langle T_j,\pi_p,p\rangle
\quad\text{for a.e. }p\in K_j,
\end{equation}
By (1.5), \cite[Theorem 6.1]{DS}, \cite[p.~1, correction to (5.14)]{DSerr},
and \eqref{eq:profile-rescaling-scale}, for some $\gamma_0>0$,
\begin{equation}\label{eq:current-bad-set-error}
|B_r\setminus K_j|
+\|T_j\|\bigl((B_r\setminus K_j)\times\R^n\bigr)
\le C(n)j^n\Bigl(\frac{\rho_x(R_j)}{R_j}\Bigr)^{2+2\gamma_0}
\le C(n)j^{-2}2^{-j}.
\end{equation}
By \eqref{eq:height-selection}, the $x$-coordinates of $\spt T_j$ over
$B_{4r}$ have diameter at most $8r\rho_x(R_j)/R_j\le r$. Thus (1.7) and
\eqref{eq:profile-excess} imply $\mathsf d(F_j(p),F_j(p'))\le C(n)r$
for $p,p'\in B_r$. By \eqref{eq:current-bad-set-error} and $r>2j$,
$|K_j|/|B_r|\to1$ as $j\to\infty$. By \eqref{eq:profile-slice-identity}
and \eqref{eq:height-selection}, every value of $F_j$ has length at most
$r$ at a.e. point of $K_j$. Hence, for $|K_j|>0$,
\[
\sup_{p\in B_r}\max_{\mathbf a\in F_j(p)}|\mathbf a|\le C(n)j.
\]
Since $R_j|x|/\rho_x(R_j)\le r<3j$ on $\spt T_j$ over $B_r$ by
\eqref{eq:height-selection}, \eqref{eq:current-bad-set-error} implies
\begin{equation}\label{eq:profile-bad-set-height}
\begin{aligned}
&\int_{B_r\setminus K_j}\mathsf d(F_j,(2n-1)\llbracket0\rrbracket)^2\,dp
+\frac{R_j^2}{\rho_x(R_j)^2}
   \int_{(B_r\setminus K_j)\times\R^n}|x|^2\,d\|T_j\|\\
\le &  C(n)j^2\bigl(|B_r\setminus K_j|
      +\|T_j\|((B_r\setminus K_j)\times\R^n)\bigr) \le C(n)2^{-j}.
\end{aligned}
\end{equation}
Next, let $W_j$ be $R_j/\rho_x(R_j)$ times the Dirichlet minimizer from Theorem 1.6.
By (1.8) and \eqref{eq:profile-excess},
\begin{equation}\label{eq:harmonic-approximation-error}
r^{-2}\int_{B_{2j}}\mathsf d(F_j,W_j)^2
+\int_{B_{2j}}\bigl(|DF_j|-|DW_j|\bigr)^2
\le C(n)j^{-n-4}r^n
\le C(n)j^{-4}.
\end{equation}
Since $r>2j$, $B_{2R}\Subset B_r$ for fixed $R>0$ and $j>R$.
Applying Caccioppoli's inequality to $W_j$ on $B_{2R}$ and using
\eqref{eq:height-selection}, \eqref{eq:profile-bad-set-height},
\eqref{eq:harmonic-approximation-error} and $r<3j$, we obtain
\begin{equation}\label{eq:normalized-value-control}
\begin{aligned}
&\int_{B_R}\Bigl(\mathsf d(W_j,(2n-1)\llbracket0\rrbracket)^2+R^2|DW_j|^2\Bigr)\\
\le & C(n)\int_{B_{2R}}\mathsf d(W_j,(2n-1)\llbracket0\rrbracket)^2\\
\le & C(n)\int_{B_{2R}}\mathsf d(F_j,(2n-1)\llbracket0\rrbracket)^2
+C(n)\int_{B_{2R}}\mathsf d(F_j,W_j)^2\\
\le &C(n)\int_{B_{2R}\cap K_j}\mathsf d(F_j,(2n-1)\llbracket0\rrbracket)^2
+C(n)2^{-j}+C(n)j^{-2}\le C(n,R).
\end{aligned}
\end{equation}
By \eqref{eq:normalized-value-control} and \cite[Proposition 3.20]{DSQ},
up to a subsequence, $W_j$ converges
in $L^2_{\mathrm{loc}}(\R^n)$ to an entire Dirichlet minimizer
$F:\R^n\to\A_{2n-1}(\R^n)$, with
$\int_{B_R}|DW_j|^2\to\int_{B_R}|DF|^2$ as $j\to\infty$, $\forall R>0$.
For $R>0$ and $j>R$, the triangle inequality,
\eqref{eq:harmonic-approximation-error} on $B_R\subset B_{2j}$, and $r<3j$
imply
\begin{equation}\label{eq:normalized-profile-convergence}
\|\mathsf d(F_j,F)\|_{L^2(B_R)}
\le C(n)j^{-1}+\|\mathsf d(W_j,F)\|_{L^2(B_R)}\longrightarrow0
\qquad\text{as }j\to\infty.
\end{equation}
For the derivatives, the reverse triangle inequality in $L^2(B_R)$
and \eqref{eq:harmonic-approximation-error} imply
\begin{equation}\label{eq:profile-energy-convergence}
\|DF_j\|_{L^2(B_R)}
=\|DW_j\|_{L^2(B_R)}+O(j^{-2})\longrightarrow\|DF\|_{L^2(B_R)}
\qquad\text{as }j\to\infty.
\end{equation}

To compare $F$ with $\spt T_j$, fix $p_0\in\R^n$ and $0<\rho\le1$.
By \cite[Theorem 3.9]{DSQ},
$F\in C^{0,\alpha}(\overline B_2(p_0))$ for some $\alpha>0$.
For $j>|p_0|+2\rho$, we have $B_{2\rho}(p_0)\Subset B_r$, since $r>2j$.
The squared distance from $R_jx/\rho_x(R_j)$ to $\spt F(p_0)$ is at most
$2R_j^2|x|^2/\rho_x(R_j)^2+2\max_{\mathbf a\in F(p_0)}|\mathbf a|^2$.
Integrating over $\pi_p^{-1}(B_{2\rho}(p_0)\setminus K_j)$ and using
\eqref{eq:profile-bad-set-height} and \eqref{eq:current-bad-set-error}, we obtain
\[
\int_{\pi_p^{-1}(B_{2\rho}(p_0)\setminus K_j)}
\dist\Bigl(\frac{R_jx}{\rho_x(R_j)},\spt F(p_0)\Bigr)^2\,d\|T_j\|
\le C(n,F,p_0)2^{-j}\longrightarrow0\qquad\text{as }j\to\infty.
\]
Applying the area formula to the graph of $\rho_x(R_j)F_j/R_j$ over
$B_{2\rho}(p_0)\cap K_j$ and using
$\frac{\rho_x(R_j)}{R_j}\operatorname{Lip}(F_j)\to0$ as $j\to\infty$
by \cite[(1.4)]{DS} and \eqref{eq:profile-excess}, we obtain
\begin{equation}\label{eq:profile-distance-integral}
\begin{aligned}
&\limsup_j\int_{\pi_p^{-1}(B_{2\rho}(p_0))}
\dist\Bigl(\frac{R_jx}{\rho_x(R_j)},\spt F(p_0)\Bigr)^2\,d\|T_j\|\\
=&\limsup_j\int_{\pi_p^{-1}(B_{2\rho}(p_0)\cap K_j)}
\dist\Bigl(\frac{R_jx}{\rho_x(R_j)},\spt F(p_0)\Bigr)^2\,d\|T_j\|\\
 \le& C(n)\limsup_j\int_{B_{2\rho}(p_0)\cap K_j}
\sum_{\mathbf b\in F_j(p)}\dist(\mathbf b,\spt F(p_0))^2\,dp\\
 \le& C(n)\limsup_j\int_{B_{2\rho}(p_0)\cap K_j}
\mathsf d(F_j(p),F(p_0))^2\,dp\\
 \le& C(n)\limsup_j\int_{B_{2\rho}(p_0)}\mathsf d(F_j(p),F(p_0))^2\,dp 
 = C(n) \int_{B_{2\rho}(p_0)}\mathsf d(F(p),F(p_0))^2\,dp\\
\le& C(n,p_0)\int_{B_{2\rho}(p_0)}|p-p_0|^{2\alpha}\,dp
\le C(n,p_0)\rho^{n+2\alpha},
\end{aligned}
\end{equation}
where \eqref{eq:normalized-profile-convergence} implies convergence of
$\mathsf d(F_j,F(p_0))$ to $\mathsf d(F,F(p_0))$ in $L^2(B_{2\rho}(p_0))$.
Since $\rho_x(R_j)^2/R_j^2\to0$ as $j\to\infty$,
\eqref{eq:profile-distance-integral} implies
\begin{equation}\label{eq:profile-local-height}
\int_{\pi_p^{-1}(B_{2\rho}(p_0))}
\min_{\mathbf a\in\spt F(p_0)}\Bigl|x-\frac{\rho_x(R_j)}{R_j}\mathbf a\Bigr|^2\,d\|T_j\|(p,x)
\longrightarrow0\qquad\text{as }j\to\infty.
\end{equation}
Let $L>|p_0|+2\rho+1$, so that
$\pi_p^{-1}(B_{2\rho}(p_0))\subset\pi_p^{-1}(B_L)$.
Since $(\pi_p)_\sharp T_j=(2n-1)\llbracket\R^n\rrbracket$, 
\begin{equation}\label{eq:profile-local-projection}
(\pi_p)_\sharp(T_j\res\pi_p^{-1}(B_{2\rho}(p_0)))
=(2n-1)\llbracket B_{2\rho}(p_0)\rrbracket.
\end{equation}
By \eqref{eq:profile-excess} and $1-J_p\ge0$,
\begin{equation}\label{eq:profile-local-mass}
\begin{aligned}
0&\le\|T_j\|(\pi_p^{-1}(B_{2\rho}(p_0)))-(2n-1)\omega_n(2\rho)^n\\
&=\int_{\pi_p^{-1}(B_{2\rho}(p_0))}(1-J_p)\,d\|T_j\|\\
&\le\int_{\pi_p^{-1}(B_L)}(1-J_p)\,d\|T_j\|
\le C(n)\frac{\rho_x(R_j)^2}{R_j^2}L^n
\longrightarrow0\quad\text{as }j\to\infty.
\end{aligned}
\end{equation}
Moreover, by \eqref{eq:height-selection},
\[
\limsup_{j\to\infty}
\sup_{(p,x)\in\spt T_j,\ |p-p_0|<2\rho}
|(p,x)-(p_0,0)|^2
\le\lim_{j\to\infty}\Bigl(4\rho^2+\frac{\rho_x(R_j)^2}{R_j^2}L^2\Bigr)
=4\rho^2<9\rho^2.
\]
Together with \eqref{eq:profile-local-height} and
\eqref{eq:profile-local-mass}, the strict inequality allows us to choose
$j_0$ such that for $j\ge j_0$ the integral in
\eqref{eq:profile-local-height} is at most $\rho^{n+2}$, the mass over
$B_{2\rho}(p_0)$ is at most $(2n-\tfrac12)\omega_n(2\rho)^n$, and
\begin{equation}\label{eq:profile-local-support}
\spt T_j\cap\pi_p^{-1}(B_{2\rho}(p_0))
\subset\mathbf B_{3\rho}\bigl((p_0,0)\bigr).
\end{equation}
Note also that $T_j$ is calibrated in $\R^n_p\times\R^n_x$ and
$\partial T_j=0$.
We choose $\rho_j\in(8n\rho,(8n+1)\rho)$ for each $j$ so that
$T_j\res\mathbf B_{\rho_j}\bigl((p_0,0)\bigr)$ is integral.
Since $3\rho<\rho_j$, the restriction to $\mathbf B_{\rho_j}((p_0,0))$
retains all of $T_j$ over $B_{2\rho}(p_0)$ by \eqref{eq:profile-local-support}. Also,
\[
\partial\Bigl(T_j\res\mathbf B_{\rho_j}\bigl((p_0,0)\bigr)\Bigr)
\res\mathbf B_{8n\rho}\bigl((p_0,0)\bigr)=0.
\]
The height and mass inequalities for $j\ge j_0$ and the
projection identity \eqref{eq:profile-local-projection} imply that
\[
\Bigl(\rho^{-1}\bigl(\operatorname{Id}-(p_0,0)\bigr)\Bigr)_\sharp
\Bigl(T_j\res\mathbf B_{\rho_j}\bigl((p_0,0)\bigr)\Bigr)
\]
satisfies \cite[Assumptions 2.1 and 3.1]{DMS}.
In the coordinate $R_jx/\rho_x(R_j)$, the estimate of
\cite[Theorem 3.2]{DMS} and \eqref{eq:profile-distance-integral} imply
\begin{equation}\label{eq:all-support}
\begin{aligned}
&\limsup_j\sup_{(p,x)\in\spt T_j,\ |p-p_0|<\rho}
\dist\Bigl(\frac{R_jx}{\rho_x(R_j)},\spt F(p_0)\Bigr)\\
\le&  C(n)\rho^{-n/2}\Bigl(\limsup_j\int_{\pi_p^{-1}(B_{2\rho}(p_0))}
\min_{\mathbf a\in\spt F(p_0)}\Bigl|\frac{R_jx}{\rho_x(R_j)}-\mathbf a\Bigr|^2\,d\|T_j\|(p,x)\Bigr)^{1/2}
\le C(n,p_0)\rho^\alpha,
\end{aligned}
\end{equation}
To estimate the distance from $R_jx/\rho_x(R_j)$ to $\spt F(p)$ uniformly
for $|p|\le1$, fix $\eta>0$.
By continuity of $F$, we cover $\overline B_1$ by finitely
many balls $B_\rho(p_0)$ with $0<\rho\le1$, $C(n,p_0)\rho^\alpha<\eta/2$ and
$\mathsf d(F(p),F(p_0))<\eta/2$ for $p\in B_\rho(p_0)$.
For $(p,x)\in\spt T_j$ with $p\in B_\rho(p_0)$,
\[
\dist\Bigl(\frac{R_jx}{\rho_x(R_j)},\spt F(p)\Bigr)
\le\dist\Bigl(\frac{R_jx}{\rho_x(R_j)},\spt F(p_0)\Bigr)
+\mathsf d(F(p),F(p_0)).
\]
By \eqref{eq:all-support} on the finitely many balls,
\[
\limsup_{j\to\infty}\sup_{(p,x)\in\spt T_j,\ |p|\le1}
\dist\Bigl(\frac{R_jx}{\rho_x(R_j)},\spt F(p)\Bigr)\le\eta.
\]
Letting $\eta\downarrow0$, we obtain
\begin{equation}\label{eq:uniform-support-distance}
\sup_{(p,x)\in\spt T_j,\ |p|\le1}
\dist\Bigl(\frac{R_jx}{\rho_x(R_j)},\spt F(p)\Bigr)
\longrightarrow0\qquad\text{as }j\to\infty.
\end{equation}
By \eqref{eq:height-function},
$\sup_{(p,x)\in\spt T_j,\ |p|\le1}R_j|x|/\rho_x(R_j)=1$.
Hence $F\not\equiv0$ by \eqref{eq:uniform-support-distance}.
By \eqref{eq:normalized-profile-convergence}, up to a subsequence,
$\mathsf d(F_j(p),F(p))\to0$ as $j\to\infty$ for a.e. $p\in\R^n$.
For fixed $R>0$ and $j>R$, $B_R\subset B_r$ since $r>2j$. Hence
\eqref{eq:current-bad-set-error} yields
\[
|B_R\setminus K_j|\le |B_r\setminus K_j|\le C(n)j^{-2}2^{-j}.
\]
Since $\sum_j|B_R\setminus K_j|<\infty$ for every $R>0$ and $O$
preserves measure,
\[
\mathbf 1_{K_j\cap O^{-1}K_j}(p)\longrightarrow1
\qquad\text{as }j\to\infty,\quad \forall O\in G,\quad\text{a.e. }p\in\R^n.
\]
By $G$-symmetry of $T_j$ and \eqref{eq:profile-slice-identity},
$F_j(Op)=OF_j(p)$ for a.e. $p\in K_j\cap O^{-1}K_j$. The a.e. convergence
of $F_j$ at $p$ and $Op$ and continuity of $F$ imply
\begin{equation}\label{eq:profile-covariance}
F(Op)=\sum_{\mathbf a\in F(p)}\llbracket O\mathbf a\rrbracket
\qquad\forall O\in G\text{ and }p\in\R^n.
\end{equation}

Furthermore, if $F(0)$ contains a nonzero
value, then, by \eqref{eq:profile-covariance} and
\eqref{eq:nonzero-orbit-size}, $F(0)$ contains at least $2n$ distinct values.
Since $F(0)$ consists of $2n-1$ vectors counted with multiplicity, we obtain
\begin{equation}\label{eq:profile-origin}
F\not\equiv0,\qquad F(0)=(2n-1)\llbracket0\rrbracket.
\end{equation}

\emph{Step 1.2.} \textit{We establish the differential and sign conditions
for $F$ and prove that the values of $F$ sum to zero.}

For each regular point $p$ of $F$ and each pair of distinct values
$\mathbf a,\mathbf a'\in F(p)$, the orthonormal bases with a vector
perpendicular to $\mathbf a-\mathbf a'$ form a set of Haar measure zero
in $O(n)$, since each basis vector is uniformly distributed on
$\mathbb S^{n-1}$. Thus, by Fubini's theorem, there exists an orthonormal basis
$b_1,\ldots,b_n$ such that, at a.e. regular point $p$ of $F$,
$\langle b_k,\mathbf a\rangle\ne\langle b_k,\mathbf a'\rangle$
for all distinct values $\mathbf a,\mathbf a'\in F(p)$ and all $k$.

For $\sum_{\ell=1}^{2n-1}\llbracket\mathbf a_\ell\rrbracket
\in\A_{2n-1}(\R^n)$, let $\Xi$ list the numbers
$\langle b_k,\mathbf a_\ell\rangle$ in increasing order for each $k$.
Matching in increasing order minimizes the squared distance in each
coordinate \cite[Section 2.1.2, inequality (Re)]{DSQ}, so $\Xi$ is
$1$-Lipschitz with respect to $\mathsf d$.
Let $DF_{j,\ell}$ denote the derivatives of the values of $F_j$, counted
with multiplicity. Sorting permutes the derivatives of the scalar
projections a.e. Hence
\[
|D(\Xi\circ F_j)|^2
=\sum_{k=1}^n\sum_{\ell=1}^{2n-1}|b_k^T DF_{j,\ell}|^2
=\sum_{\ell=1}^{2n-1}|DF_{j,\ell}|^2=|DF_j|^2
\]
by orthonormality.

Near a regular point $p$ of $F$, $F=\sum_J q_J\llbracket f_J\rrbracket$
with distinct $f_J(p)$. At a.e. such $p$, the projections
$\langle b_k,f_J(p)\rangle$ are distinct as $J$ varies, for each $k$.
If $\langle b_k,f_J(p)\rangle<\langle b_k,f_{J'}(p)\rangle$, continuity implies
\[
\langle b_k,f_J(p')\rangle<\langle b_k,f_{J'}(p')\rangle
\qquad\text{for }p'\text{ near }p.
\]
Thus for $p'$ near $p$, the projections $\langle b_k,f_J(p')\rangle$ appear
in the $k$th list of $\Xi(F(p'))$ in the order determined by
$\langle b_k,f_J(p)\rangle$. The function
$\langle b_k,f_J\rangle$ occurs $q_J$ times, each with derivative
$b_k^T Df_J(p)$ at $p$.
Hence, by orthonormality,
\[
|D(\Xi\circ F)(p)|^2
=\sum_J q_J\sum_{k=1}^n |b_k^T Df_J(p)|^2
=\sum_J q_J|Df_J(p)|^2=|DF(p)|^2.
\]
Since $\Xi$ is Lipschitz, $\Xi\circ F_j\to\Xi\circ F$ in $L^2(B_R)$
as $j\to\infty$, $\forall R>0$, by \eqref{eq:normalized-profile-convergence}.
By $|D(\Xi\circ F_j)|^2=|DF_j|^2$ and \eqref{eq:profile-energy-convergence},
$\limsup_{j\to\infty}\|D(\Xi\circ F_j)\|_{L^2(B_R)}<\infty$.
The $L^2(B_R)$ convergence then implies
$\Xi\circ F_j\rightharpoonup\Xi\circ F$ in $W^{1,2}(B_R)$ as $j\to\infty$.
Expanding the square and using the weak $W^{1,2}(B_R)$ convergence, we obtain
\[
\begin{aligned}
\int_{B_R}|D(\Xi\circ F_j)-D(\Xi\circ F)|^2
&=\int_{B_R}\bigl(|D(\Xi\circ F_j)|^2-|D(\Xi\circ F)|^2\bigr)+o(1)\\
&=\int_{B_R}(|DF_j|^2-|DF|^2)+o(1)
\longrightarrow0\quad\text{as }j\to\infty.
\end{aligned}
\]
The last limit follows from \eqref{eq:profile-energy-convergence}.
Thus, up to a subsequence, $D(\Xi\circ F_j)(p)\to D(\Xi\circ F)(p)$
as $j\to\infty$ for a.e. $p\in\R^n$.

At a.e. regular point $p$ of $F$, the convergence $F_j(p)\to F(p)$
allows us to index the values of $F_j(p)$ so that
$F_{j,\ell}(p)\to f_J(p)$ as $j\to\infty$ for $\ell\in J$, $\#J=q_J$.
For each $b_k$, the $\langle b_k,F_{j,\ell}(p)\rangle$ with $\ell\in J$
converge to $\langle b_k,f_J(p)\rangle$. The limits for distinct $J$
are different, so the $q_J$ rows at the corresponding sorted positions
of $D(\Xi\circ F_j)(p)$ converge to $b_k^T Df_J(p)$. Consequently,
\begin{equation}\label{eq:cluster-variance}
\sum_{J,\ell\in J}|DF_{j,\ell}(p)-Df_J(p)|^2
\longrightarrow0\qquad\text{as }j\to\infty
\end{equation}
by orthonormality of $b_1,\ldots,b_n$.

By \eqref{eq:profile-slice-identity}, for a.e. $p\in K_j$ and
$1\le\ell\le2n-1$, the tangent plane of $T_j$ at
$\bigl(p,\frac{\rho_x(R_j)}{R_j}F_{j,\ell}(p)\bigr)$ is the graph
$x=Ap$, where $A=\frac{\rho_x(R_j)}{R_j}DF_{j,\ell}(p)=A^T$
and $J_p>0$.
The corresponding oriented tangent plane of $S$ satisfies
\[
T_{(\rho_x(R_j)F_{j,\ell}(p),R_jp)}(\spt S)\in\Gr_n((n-2)\pi/2).
\]
Thus there exist $\mathsf H_k$ with
$\Theta(\mathsf H_k)=(n-2)\pi/2$, $\det\mathsf H_k<0$ and
$\mathsf H_k^{-1}\to A$ as $k\to\infty$.
By \eqref{eq:two-smallest-angles},
\[
n_-(DF_{j,\ell}(p))=n_-(A)
\le\liminf_{k\to\infty}n_-(\mathsf H_k^{-1})=1,
\]
where $n_-$ counts negative eigenvalues with multiplicity.
Since $T_j$ has phase zero, $\Theta(A)\in2\pi\mathbb Z$.
On the other hand,
\[
|\Theta(A)|\le n\arctan\Bigl(\frac{\rho_x(R_j)}{R_j}\operatorname{Lip}(F_j)\Bigr)
\longrightarrow0\quad\text{as }j\to\infty.
\]
Since $0$ is isolated in $2\pi\mathbb Z$,
$R_j\Theta(A)/\rho_x(R_j)\to0$ as $j\to\infty$. Using
$\arctan t=t+O(t^3)$ as $t\to0$ and \eqref{eq:cluster-variance},
we obtain
\begin{equation}\label{eq:normalized-slice-trace}
\begin{aligned}
0&=\lim_{j\to\infty}\frac{R_j}{\rho_x(R_j)}
\Theta\Bigl(\frac{\rho_x(R_j)}{R_j}DF_{j,\ell}(p)\Bigr)\\
&=\lim_{j\to\infty}\left[\tr DF_{j,\ell}(p)
  +O\Bigl(\frac{\rho_x(R_j)^2}{R_j^2}|DF_{j,\ell}(p)|^3\Bigr)\right]
=\tr Df_J(p),\qquad\forall\ell\in J.
\end{aligned}
\end{equation}
Consequently, by \eqref{eq:cluster-variance} and
\eqref{eq:normalized-slice-trace}, each regular analytic branch $f$ of $F$
satisfies
\begin{equation}\label{eq:branch-jets}
Df=(Df)^T,\qquad \tr Df=0,\qquad n_-(Df)\le1
\end{equation}
at a.e. point in the domain of $f$, and hence at every point by continuity
of $Df$.

Moreover, by \eqref{eq:height-selection}, \eqref{eq:profile-slice-identity}
and continuity of $F$,
\[
\max_{\mathbf a\in F(p)}|\mathbf a|
=\lim_{j\to\infty}\max_{\mathbf a\in F_j(p)}|\mathbf a|
\le\limsup_{j\to\infty}
\frac{\rho_x(R_j\max\{1,|p|\})}{\rho_x(R_j)}
\le\max\{1,|p|\}
\]
for a.e. $p$. By continuity of $F$,
$\max_{\mathbf a\in F(p)}|\mathbf a|\le\max\{1,|p|\}$ for all $p\in\R^n$.
Thus $\sum_{\mathbf a\in F(p)}\mathbf a$ grows at most linearly in $|p|$
and is affine because it is harmonic \cite[Lemma 3.23(a)]{DSQ}.
By \eqref{eq:profile-covariance}, the affine sum is a scalar multiple of $p$.
By $\tr Df=0$ in \eqref{eq:branch-jets},
$\operatorname{div}_p\bigl(\sum_{\mathbf a\in F(p)}\mathbf a\bigr)=0$.
Therefore
\begin{equation}\label{eq:zero-barycenter}
\sum_{\mathbf a\in F(p)}\mathbf a=0\qquad\forall p\in\R^n.
\end{equation}
Likewise, transferring \eqref{eq:closed-coordinate-signs} through
\eqref{eq:profile-slice-identity} and the a.e. convergence of $F_j$ to $F$,
then using continuity, we obtain
\begin{equation}\label{eq:value-cones}
\begin{gathered}
\max_k|a_k|\le-a_i\quad\text{if }p_i>0,\ a_i<0,\ \mathbf a\in F(p),\\
\#\{i:a_i<0\}\le1\quad\forall p\in(0,\infty)^n,\ \mathbf a\in F(p).
\end{gathered}
\end{equation}

\emph{Step 1.3.} \textit{We show that}
\begin{equation}\label{eq:nozero}
0\notin\spt F(p)\qquad\forall p\in\R^n\setminus\{0\}.
\end{equation}

Let $f$ be a regular branch of $F$ on a connected open subset of
$\{p_i>0\}\setminus\Sing F$.
If $f_i<0$ at one point and $f_i=0$ at another, choose a ball
$B\subset\{f_i<0\}$ with closure compactly contained in the domain of $f$
and with $f_i(q)=0$ for some $q\in\partial B$.
By the Hopf lemma for the harmonic function $f_i$, $Df_i(q)\ne0$.
The inequalities $|f_k|\le-f_i$ on $\{f_i<0\}$ from \eqref{eq:value-cones}
imply $f=0$ on the smooth hypersurface $\{f_i=0\}$ near $q$.
Hence $\operatorname{rank}Df(q)\le1$, and \eqref{eq:branch-jets} implies
\[
0<|Df_i(q)|^2\le|Df(q)|^2=(\tr Df(q))^2=0,
\]
a contradiction. Thus, counting with multiplicity, the number of vectors
in $F(p)$ with negative $i$th coordinate is locally constant as $p$ varies
in $\{p_i>0\}\setminus\Sing F$.
By \cite[Theorems 1.1 and 1.3]{DMSV}, $\Sing F$ is closed and countably
$\mathcal H^{n-2}$-finite, so $\Sing F$ has zero $W^{1,2}$ capacity locally
and $\{p_i>0\}\setminus\Sing F$ is connected.
Thus, for each $i$, the number of vectors in $F(p)$ with negative $i$th
coordinate, counted with multiplicity, is constant as $p$ varies in
$\{p_i>0\}\setminus\Sing F$.
By \eqref{eq:profile-origin}, \eqref{eq:profile-covariance} and continuity,
$F(p)$ contains a nonzero vector for some $p\in(0,\infty)^n\setminus\Sing F$.
Since $\sum_{\mathbf a\in F(p)}\mathbf a=0$ by \eqref{eq:zero-barycenter},
there exist $\mathbf a\in F(p)$ and $i\in\{1,\ldots,n\}$ such that $a_i<0$.
To compare the counts across coordinates, take $O\in G$
interchanging coordinates $i$ and $k$. For
$p\in(0,\infty)^n\setminus\Sing F$, we have
\[
\#\{\mathbf a\in F(p):a_i<0\}
=\#\{\mathbf a\in F(Op):a_k<0\}
=\#\{\mathbf a\in F(p):a_k<0\},
\]
where we used \eqref{eq:profile-covariance} in the first equality.
For the second equality, $p,Op\in\{q:q_k>0\}\setminus\Sing F$, and
$\#\{\mathbf a\in F(q):a_k<0\}$ is constant on $\{q:q_k>0\}\setminus\Sing F$.
Consequently, for any $p\in(0,\infty)^n\setminus\Sing F$, counting with multiplicity,
\[
\begin{aligned}
1\le\#\{\mathbf a\in F(p):a_i<0\}
&=\frac1n\sum_{k=1}^n\#\{\mathbf a\in F(p):a_k<0\}\\
&=\frac1n\sum_{\mathbf a\in F(p)}\#\{k:a_k<0\}
\le\frac{2n-1}{n}<2.
\end{aligned}
\]
By constancy of the count on $\{p_i>0\}\setminus\Sing F$, for each $p$
in the half-space outside $\Sing F$ there is a unique $\mathbf a\in F(p)$
with $a_i<0$, denoted by $f^{(i)}(p)$.
Locally, $f^{(i)}$ is a regular branch of $F$ and hence harmonic.
Local boundedness of $F$ and zero $W^{1,2}$ capacity of $\Sing F$
locally allow us to extend $f^{(i)}$ with $\Delta f^{(i)}=0$ in $\{p_i>0\}$.
By the strong minimum principle and continuity of $Df^{(i)}$,
we have $f_i^{(i)}<0$ and \eqref{eq:branch-jets} on $\{p_i>0\}$.
Moreover, by uniqueness and \eqref{eq:profile-covariance},
\begin{equation}\label{eq:negative-branch-covariance}
f^{(j)}(Op)=Of^{(i)}(p)
\qquad\forall O\in G\text{ with }Oe_i=e_j,\quad p_i>0.
\end{equation}

Since $f_i^{(i)}<0$ and every vector in $F(p)$ has at most one negative
coordinate by \eqref{eq:value-cones}, we have $f^{(i)}(p)\ne f^{(k)}(p)$
for $i\ne k$ and $p\in(0,\infty)^n$. Hence
\[
F_+(p)\coloneqq F(p)-\sum_{i=1}^n\llbracket f^{(i)}(p)\rrbracket
\in\A_{n-1}([0,\infty)^n)
\qquad\forall p\in(0,\infty)^n.
\]
By \eqref{eq:zero-barycenter},
\[
f_+\coloneqq\sum_{\mathbf a\in F_+}\mathbf a=-\sum_{k=1}^nf^{(k)}.
\]
If $f_+\equiv0$, then $f^{(i)}=-\sum_{j\ne i}f^{(j)}$ extends smoothly
across $\{p_i=0,\ p_k>0\ (k\ne i)\}$.
Using reflection in \eqref{eq:negative-branch-covariance}, we have
$f_i^{(j)}=0$ for $j\ne i$, and hence $f_i^{(i)}=0$ on
$\{p_i=0,\ p_k>0\ (k\ne i)\}$.
By \eqref{eq:value-cones}, $f^{(i)}=0$ on
$\{p_i=0,\ p_k>0\ (k\ne i)\}$.
Using \eqref{eq:branch-jets} and  Hopf lemma,
we obtain 
\[
0=\tr Df^{(i)}=\partial_if_i^{(i)}<0
\quad\text{on }\{p_i=0,\ p_k>0\ \forall k\ne i\},
\]
a contradiction. Hence $f_+\not\equiv0$.
By $f_+=-\sum_k f^{(k)}$, each $f_{+,i}$ is harmonic. Since
$f_+\not\equiv0$ and permutations exchange the nonnegative coordinates
$f_{+,i}$, each coordinate is positive somewhere. The strong minimum
principle then implies $f_{+,i}>0$ on $(0,\infty)^n$, $\forall i$.
On the positive diagonal $p=re_+$, where $e_+\coloneqq(1,\ldots,1)/\sqrt n$,
$F_+(re_+)$ is invariant under coordinate permutations by
\eqref{eq:profile-covariance} and \eqref{eq:negative-branch-covariance}.
If a vector has two unequal coordinates, then we obtain at least $n$
distinct vectors by permuting the coordinates.
As $F_+$ has only $n-1$ values,
\begin{equation}\label{eq:positive-diagonal-values}
\spt F_+(re_+)\subset[0,\infty)e_+\qquad\forall r>0.
\end{equation}
Since $\langle f_+(re_+),e_+\rangle>0$ for $r>0$ and
$f_+(re_+)\to0$ as $r\downarrow0$, there exists $r>0$ with
$\langle Df_+(re_+)e_+,e_+\rangle>0$.
Differentiating $f_+(Op)=Of_+(p)$ at $p=re_+$ for coordinate permutation
matrices $O$, and using
$\tr Df_+=-\sum_{k=1}^n\tr Df^{(k)}=0$ from
\eqref{eq:branch-jets}, we obtain
\begin{equation}\label{eq:diagonal-derivative}
Df_+(re_+)=\langle Df_+(re_+)e_+,e_+\rangle
\Bigl[e_+\otimes e_+-\frac{I-e_+\otimes e_+}{n-1}\Bigr].
\end{equation}
Thus $n_-(Df_+(re_+))=n-1$.

On the other hand, suppose $q\in\{p_i>0\}$ and $a_i=0$ for some
$\mathbf a\in F(q)$. By continuity of $F$, we can choose a ball
$B(q)\Subset\{p_i>0\}$ on which the vectors $\mathbf a\in F(q)$ with $a_i=0$
determine a separated local component $F^{(0)}$ of $F$.
Since $f^{(i)}(p)$ is outside $F^{(0)}(p)$, we have $a_i\geq 0$ for
$\mathbf a\in F^{(0)}(p)$ and $p\in B(q)$.
The sum $\sum_{\mathbf b\in F^{(0)}(p)}\mathbf b$ is locally bounded and
harmonic on $B(q)\setminus\Sing F$, where the summands are locally regular
branches of $F$. The vector sum $\sum_{\mathbf b\in F^{(0)}(p)}\mathbf b$
extends harmonically across $\Sing F\cap B(q)$,
which has zero $W^{1,2}$ capacity locally.
At $q$, $\sum_{\mathbf b\in F^{(0)}(q)}b_i=0$. By the strong minimum principle,
\begin{equation}\label{eq:component-zero-coordinate}
0\le a_i\le\sum_{\mathbf b\in F^{(0)}(p)}b_i=0
\qquad\forall p\in B(q),\quad\mathbf a\in F^{(0)}(p).
\end{equation}
The set
\[
\bigl\{p\in\R^n:p_i>0,\ \min_{\mathbf a\in F(p)}|a_i|=0\bigr\}
\]
is open in $\{p_i>0\}$ by \eqref{eq:component-zero-coordinate} and closed
in $\{p_i>0\}$ by continuity of $F$.
Suppose $(p_0)_i>0$ and $a_i=0$ for some $\mathbf a\in F(p_0)$.
Since $\{p_i>0\}$ is connected, 
\begin{equation}\label{eq:zero-coordinate-existence}
\min_{\mathbf a\in F(p)}|a_i|=0
\qquad\forall p\text{ with }p_i>0.
\end{equation}
If $f_i^{(j)}(p)=0$ at some $p\in(0,\infty)^n$ for $j\ne i$,
we obtain $f_k^{(j)}\equiv0$ for all $k\ne j$ by the minimum principle
and \eqref{eq:negative-branch-covariance} for permutations fixing $e_j$.
Using symmetry and zero trace in \eqref{eq:branch-jets}, we obtain
\[
|Df^{(j)}|^2=(\partial_jf_j^{(j)})^2
=\Bigl(\sum_{k\ne j}\partial_kf_k^{(j)}\Bigr)^2=0.
\]
Hence $f^{(j)}$ is constant, so $f^{(j)}\equiv0$ by \eqref{eq:profile-origin},
contradicting $f_j^{(j)}<0$.
Therefore $f_i^{(j)}>0$ on $(0,\infty)^n$ for $i\ne j$.
By \eqref{eq:positive-diagonal-values} and the signs of the coordinates
of $f^{(j)}(re_+)$, if $\mathbf a\in F(re_+)$ and $\mathbf a\ne0$, then
$a_i\ne0$ for all $i$. Thus \eqref{eq:zero-coordinate-existence} implies
$0\in\spt F_+(re_+)$.
By continuity of $F_+$, choose $B(re_+)\Subset(0,\infty)^n$ so that the
vectors of $F_+(p)$ that tend to $0$ as $p\to re_+$ form a separated
component of $F_+$ for $p\in B(re_+)$. The coordinate sums of the separated
component extend harmonically across $\Sing F\cap B(re_+)$. Since
$F_+(p)\in\A_{n-1}([0,\infty)^n)$, each coordinate sum is $\geq0$
on $B(re_+)$ and equals $0$ at $re_+$. By the strong minimum principle the coordinate
sums equal $0$ on $B(re_+)$, so each vector of the separated component
also equals $0$ on $B(re_+)$. Near each
regular point of $F$ in $B(re_+)$, we can therefore label the branches of $F_+$, counted
with multiplicity, so that
\[
F_+=\llbracket0\rrbracket+\sum_{\ell=1}^{n-2}\llbracket f_\ell\rrbracket,
\qquad f_+=\sum_{\ell=1}^{n-2}f_\ell.
\]
Using \eqref{eq:branch-jets},
\begin{equation}\label{eq:remaining-branch-index}
n_-(Df_+)\le\sum_{\ell=1}^{n-2}n_-(Df_\ell)\le n-2.
\end{equation}
By density of the regular points of $F$ and continuity of the eigenvalues
of $Df_+$, we obtain from
\eqref{eq:diagonal-derivative} and \eqref{eq:remaining-branch-index}
\[
n-1=n_-(Df_+(re_+))
\le\liminf_{p\to re_+,\ p\notin\Sing F}n_-(Df_+(p))\le n-2,
\]
a contradiction. Thus $a_i\ne0$ for $\mathbf a\in F(p)$ when $p_i>0$.
Reflecting the $i$th coordinate in \eqref{eq:profile-covariance}, we obtain
\[
a_i\ne0\qquad\forall p\in\R^n\text{ with }p_i\ne0,\quad\mathbf a\in F(p).
\]
For $p\ne0$, choose $i$ with $p_i\ne0$. Then $a_i\ne0$, and hence
$\mathbf a\ne0$ for every $\mathbf a\in F(p)$, proving \eqref{eq:nozero}.

\emph{Step 2.} \textit{We extend $u$ smoothly from $\pi_x(\Lambda)$ to $\R^n$.}

If $(x,p)\in\spt S$ and $|p|=R_j$, then $(p/R_j,x/R_j)\in\spt T_j$ and
\[
\frac{|x|}{\rho_x(R_j)}\ge
\min_{\mathbf a\in\spt F(p/R_j)}|\mathbf a|
-\dist\Bigl(\frac{x}{\rho_x(R_j)},\spt F(p/R_j)\Bigr).
\]
By \eqref{eq:uniform-support-distance},
$\dist(x/\rho_x(R_j),\spt F(p/R_j))\to0$ uniformly on
$\spt S\cap\{|p|=R_j\}$ as $j\to\infty$.
Since $p/R_j\in\mathbb S^{n-1}$, the minimum over $\mathbf a\in\spt F(p/R_j)$
is bounded below by $\min_{q\in\mathbb S^{n-1},\ \mathbf a\in\spt F(q)}|\mathbf a|>0$
by \eqref{eq:nozero} and continuity of $F$.
Since $\rho_x(R_j)\to\infty$,
\begin{equation}\label{eq:sphere-barrier}
\inf_{(x,p)\in\spt S,\ |p|=R_j}|x|\longrightarrow\infty
\qquad\text{as }j\to\infty.
\end{equation}

On the other hand, let $y\in\pi_x(\Lambda)$ and choose $k$ with
$|y_k|=\max_i|y_i|$.
By reflection symmetry, we may assume $y_i\ge0$ for all $1\le i\le n$.
For $z\in\R^n$, define
\[
\Psi_y(z)\coloneqq
\Bigl(\min\bigl\{\max_{1\le j\le n}|z_j|,\max\{y_i,|z_i|\}\bigr\}\Bigr)_{i=1}^n.
\]
Since ${\Psi_y}$ is continuous and ${\Psi_y}(x)=y$ whenever $x_k=y_k$ and
$|x_i|\le y_i$, $\forall i\ne k$, there exists a compact neighborhood
$K\Subset\pi_x(\Lambda)$ of $y$ and a bounded connected open set $U\subset\R^n$
such that
\[
\bigl\{x\in\R^n:x_k=y_k,\ |x_i|\le y_i,\,\,  \forall  i\ne k\bigr\}
\subset U\subset {\Psi_y}^{-1}(K).
\]
Since $U$ is bounded and $r_{t_j}\to0$ as $j\to\infty$, using
\eqref{eq:signs}, reflection symmetry and \eqref{eq:local-potential-limit},
we obtain
\[
\begin{aligned}
\limsup_{j\to\infty}\sup_{z\in U}u_{t_j}(z)
&\le\limsup_{j\to\infty}\sup_{z\in U}u_{t_j}(\Psi_y(z)) \le\lim_{j\to\infty}\sup_Ku_{t_j}=\sup_Ku<\infty.
\end{aligned}
\]
By Lemma~\ref{lem:interior-estimate}, the uniform upper estimate on $U$
and the convergence $u_{t_j}(y)\to u(y)$ from \eqref{eq:local-potential-limit}
imply, up to a further subsequence,
\[
u_{t_j}\longrightarrow u_U\quad\text{in }C^\infty_{\mathrm{loc}}(U)
\quad\text{as }j\to\infty.
\]
Since $\Graph(Du_{t_j})\rightharpoonup S$ locally as currents as $j\to\infty$,
$S\res(U\times\R^n)=\Graph(Du_U)$.
By \eqref{eq:local-potential-limit}, we have $u_U=u$ on $U\cap\pi_x(\Lambda)$.
As $U$ is connected and $y\in U\cap\pi_x(\Lambda)$,
$\Graph(Du_U)\subset\Lambda=\Graph(Du)$.
Hence $U\subset\pi_x(\Lambda)$, and therefore
\begin{equation}\label{eq:box}
\bigl\{x\in\R^n:x_k=y_k,\ |x_i|\le|y_i|\ \forall i\ne k\bigr\}
\subset\pi_x(\Lambda).
\end{equation}
Since $\pi_x(\Lambda)$ is connected and contains $0$ and $y$, the
continuous function $z\mapsto\max_i|z_i|$ takes every value in $[0,|y_k|]$
on $\pi_x(\Lambda)$. By coordinate symmetry, for $0\le\lambda\le1$,
we can choose $z\in\pi_x(\Lambda)$ with $z_k=\lambda y_k$ and
$\max_i|z_i|=|z_k|$. Applying \eqref{eq:box} to $z$ yields
$\lambda y_ke_k\in\pi_x(\Lambda)$, and applying \eqref{eq:box} to $y$ yields
$(1-\lambda)y_ke_k+\lambda y\in\pi_x(\Lambda)$, $\forall\lambda\in[0,1]$.
The ball $\overline B_{|y|}$ is convex and contains $0$, $y_ke_k$ and $y$.
Joining the two segments, we obtain a piecewise linear path $\gamma$
from $0$ to $y$ through $y_ke_k$ with
\begin{equation}\label{eq:controlled-path}
\gamma(0)=0,\qquad\gamma(1)=y,\qquad
\gamma([0,1])\subset\pi_x(\Lambda)\cap\overline B_{|y|}.
\end{equation}

Fix $R>0$ and suppose $\sup_{y\in\pi_x(\Lambda)\cap B_R}|Du(y)|=\infty$.
Starting from $Du(0)=0$ and applying the intermediate value theorem
along the paths in \eqref{eq:controlled-path}, we obtain
\[
\liminf_{j\to\infty}\inf_{(x,p)\in\spt S,\ |p|=R_j}|x|\le R,
\]
contradicting \eqref{eq:sphere-barrier}. Hence
\begin{equation}\label{eq:gradient-on-balls}
\sup_{x\in\pi_x(\Lambda)\cap B_R}|Du(x)|<\infty
\qquad\forall R>0.
\end{equation}
For $x\in\partial(\pi_x(\Lambda))$, take $x_\ell\in\pi_x(\Lambda)$ with
$x_\ell\to x$. By \eqref{eq:gradient-on-balls}, up to a subsequence,
$(x_\ell,Du(x_\ell))\to(x,p)$ for some $p\in\R^n$.
Then $(x,p)\in\overline\Lambda\setminus\Lambda\subset\Sing S$
because $x\notin\pi_x(\Lambda)$. Thus, by \cite[Theorem 1.1]{DMS},
\[
\dim_{\mathcal H}\partial(\pi_x(\Lambda))
\le\dim_{\mathcal H}\pi_x(\Sing S)\le n-2.
\]
If $\pi_x(\Lambda)$ is not dense, take a ball outside
$\overline{\pi_x(\Lambda)}$ and a ball inside $\pi_x(\Lambda)$.
Join the centers of the two balls by a segment. Every small parallel translate of the
segment meets $\partial(\pi_x(\Lambda))$, so the orthogonal projection
of $\partial(\pi_x(\Lambda))$ onto the hyperplane perpendicular to the
segment contains an $(n-1)$-dimensional ball. This contradicts
$\dim_{\mathcal H}\partial(\pi_x(\Lambda))\le n-2$.
Consequently,
\begin{equation}\label{eq:graph-domain-complement}
\R^n\setminus\pi_x(\Lambda)=\partial(\pi_x(\Lambda)).
\end{equation}

Fix $x\in\R^n$. The radial projection from $x$ of
$\partial(\pi_x(\Lambda))\setminus\{x\}$ onto $\mathbb S^{n-1}$
has Hausdorff dimension at most $n-2$. Hence there exists a unit vector
$e$ such that $x+\R e$ meets $\partial(\pi_x(\Lambda))$ only possibly at $x$.
For $(r,e')\in(1,2)\times(\mathbb S^{n-1}\cap e^\perp)$,
consider the circle $x+r(\cos\vartheta\,e+\sin\vartheta\,e')$.
If $y=x+r(\cos\vartheta\,e+\sin\vartheta\,e')$ belongs to
$\partial(\pi_x(\Lambda))$, then
$y\notin x+\R e$ and
\[
(r,e')=\Bigl(|y-x|,\ \pm\frac{P_{e^\perp}(y-x)}{|P_{e^\perp}(y-x)|}\Bigr),
\]
where $P_{e^\perp}$ is the orthogonal projection onto $e^\perp$.
For either sign, the displayed map from $y$ to $(r,e')$ is locally
Lipschitz on $\R^n\setminus(x+\R e)$. The images of
$\partial(\pi_x(\Lambda))\setminus(x+\R e)$ under the two parameter maps have
Hausdorff dimension at most $n-2$. Since
$(1,2)\times(\mathbb S^{n-1}\cap e^\perp)$ has dimension $n-1$, we can
choose $(r,e')$ such that $x+r(\cos\vartheta\,e+\sin\vartheta\,e')\in\pi_x(\Lambda)$
for $0\le\vartheta\le2\pi$, by \eqref{eq:graph-domain-complement}.
Hence there exists $\delta>0$ such that
\[
\bigl\{x'+r(\cos\vartheta\,e+\sin\vartheta\,e'):
|x'-x|\le\delta,\ 0\le\vartheta\le2\pi\bigr\}\Subset\pi_x(\Lambda).
\]
By $D^2u_{t_j}(y)=D^2v^{(t_j)}(r_{t_j}y)$ and
\eqref{eq:two-smallest-angles}, for $y\in B_{1/r_{t_j}}$,
\[
\begin{aligned}
\langle D^2u_{t_j}(y)e,e\rangle+\langle D^2u_{t_j}(y)e',e'\rangle
&\ge\lambda_{n-1}(D^2u_{t_j}(y))+\lambda_n(D^2u_{t_j}(y))>0 .
\end{aligned}
\]
Thus $u_{t_j}$ is subharmonic on the affine plane
$x'+\operatorname{span}\{e,e'\}$ wherever $u_{t_j}$ is defined.
Since $r_{t_j}(|x|+\delta+r)\to0$ as $j\to\infty$, the mean-value
inequality on $B_r(x')\cap(x'+\operatorname{span}\{e,e'\})$ for
$|x'-x|\le\delta$ and \eqref{eq:local-potential-limit} imply
\[
\begin{aligned}
\limsup_{j\to\infty}\sup_{|x'-x|<\delta}u_{t_j}(x')
&\le\limsup_{j\to\infty}\sup_{|x'-x|\le\delta}
\frac1{2\pi}\int_0^{2\pi}u_{t_j}\bigl(x'+r(\cos\vartheta\,e+\sin\vartheta\,e')\bigr)\,d\vartheta\\
&=\sup_{|x'-x|\le\delta}
\frac1{2\pi}\int_0^{2\pi}u\bigl(x'+r(\cos\vartheta\,e+\sin\vartheta\,e')\bigr)\,d\vartheta<\infty.
\end{aligned}
\] 
Covering $\overline B_R(0)$ by finitely many neighborhoods $B_\delta(x)$
on which $\limsup_{j\to\infty}\sup_{B_\delta(x)}u_{t_j}<\infty$, we obtain
\[
\limsup_{j\to\infty}\sup_{x\in B_R(0)}u_{t_j}(x)<\infty
\qquad\forall R>0.
\]
Since $\{u_{t_j}\}$ is locally bounded above and $u_{t_j}(0)=0$,
Lemma~\ref{lem:interior-estimate} implies that, up to a subsequence,
$u_{t_j}\to u$ in $C^\infty_{\mathrm{loc}}(\R^n)$ as $j\to\infty$,
where $u$ extends the limit on $\pi_x(\Lambda)$ in \eqref{eq:local-potential-limit}.
Consequently, $u$ satisfies
$\Theta(D^2u)=\frac{n-2}{2}\pi$,
$u(0)=0$, $Du(0)=0$, $D^2u(0)=\cot(\pi/n) I$ and $\osc_{B_1}u=\cot(\pi/n)$.

\emph{Step 3.} \textit{We prove that $u$ is nonpolynomial.}

Suppose $u$ is a polynomial. The symmetry $u_{t_j}\circ O=u_{t_j}$
and local uniform convergence imply $u\circ O=u$, $\forall O\in G$.
Differentiating $u(Ox)=u(x)$ twice with respect to $x$ and evaluating
at $\tau e_1$, we obtain
\[
D^2u(\tau e_1)=O^TD^2u(\tau Oe_1)O=O^TD^2u(\tau e_1)O
\qquad\forall\tau\in\R,\quad O\in G\text{ with }Oe_1=e_1.
\]
Using coordinate reflections and permutations fixing $e_1$, we have
$\partial_{ik}u(\tau e_1)=0$ for $i\ne k$ and
$\partial_{22}u(\tau e_1)=\cdots=\partial_{nn}u(\tau e_1)$. Therefore
\[
D^2u(\tau e_1)=\diag(a(\tau),b(\tau),\ldots,b(\tau))
\qquad\forall\tau\in\R,
\]
where $a(\tau),b(\tau)$ are real polynomials.
Hence
\begin{equation}\label{eq:axis-principal-phase}
\arctan a(\tau)+(n-1)\arctan b(\tau)
=\Theta(D^2u(\tau e_1))=\frac{n-2}{2}\pi\qquad\forall\tau\in\R.
\end{equation}
Since $\arctan a(\tau)<\pi/2$, \eqref{eq:axis-principal-phase} implies
$b(\tau)>0$. If $b$ is nonconstant, then $b(\tau)\to+\infty$
as $\tau\to+\infty$. By \eqref{eq:axis-principal-phase},
$\arctan a(\tau)+\pi/2=(n-1)\arctan(1/b(\tau))$, so
\[
\begin{aligned}
a(\tau)+\frac{b(\tau)}{n-1}
&=-\cot\bigl((n-1)\arctan(1/b(\tau))\bigr)+\frac{b(\tau)}{n-1} =O\bigl(b(\tau)^{-1}\bigr)\longrightarrow0
\qquad\text{as }\tau\to+\infty.
\end{aligned}
\]
Thus  $a+b/(n-1)\equiv 0$, contradicting 
$a(0)+b(0)/(n-1)>0$ by $D^2u(0)=\cot(\pi/n)I$.
Therefore $b$ is constant, and \eqref{eq:axis-principal-phase} implies
that $a$ is constant. Hence $D^2u(se_1)=\cot(\pi/n)I$ for all $s\in\R$.
For $R>0$, take $\tau\coloneqq R\tan(\pi/n)$. Since $Du(0)=0$,
\[
Du(\tau e_1)=\int_0^\tau D^2u(se_1)e_1\,ds
=\tau\cot(\pi/n)e_1=Re_1.
\]
By smooth convergence on compact subsets of $\R^n$,
$(\tau e_1,Re_1)\in\spt S$. Hence, by \eqref{eq:height-function},
\[
\frac{\rho_x(R)}R\ge\frac\tau R=\tan(\pi/n)>0\qquad\forall R>0,
\]
which contradicts \eqref{eq:sublinear-x-height}.
\end{proof}

\appendix
\section{Local invertibility of the projection onto \texorpdfstring{$x$}{x}}\label{app:partial}
\label{app:vertical}

\begin{lemma}\label{lem:partial-propagation}
Let $n\ge3$ and let $\Lambda\subset\C^n$ be a real analytic special
Lagrangian submanifold with $T_z\Lambda\in\Gr_n((n-2)\pi/2)$,
$\forall z\in\Lambda$. Suppose $\operatorname{rank}d(\pi_x|_\Lambda)_{z_0}=r$
at $z_0\in\Lambda$, with $0<r<n$. Then $r\ge2$ and
$\operatorname{rank}d(\pi_x|_\Lambda)=r$ near $z_0$. Moreover, $\pi_x$ maps
a sufficiently small neighborhood of $z_0$ in $\Lambda$ onto an open
subset of an affine $r$-plane if $r\ge3$ and onto a minimal surface if $r=2$.
\end{lemma}
\begin{proof}
Let $k\coloneqq n-r$. At $z_0$, the projection rank is
$r=\#\{i:|\theta_i|<\pi/2\}>0$. Hence, by \eqref{eq:closed-principal-angles},
\[
\theta_i=\frac{n-2}{2}\pi-\sum_{j\ne i}\theta_j>-\frac{\pi}{2}
\qquad\forall i\text{ with }|\theta_i|=\frac{\pi}{2}.
\]
Thus each of the $k$ angles with $|\theta_i|=\pi/2$ equals $\pi/2$,
and
\[
-\frac{r\pi}{2}<\sum_{|\theta_i|<\pi/2}\theta_i
=\frac{n-2}{2}\pi-\sum_{|\theta_i|=\pi/2}\theta_i
=\frac{n-2-k}{2}\pi=\frac{r-2}{2}\pi,
\]
so $r\ge2$.

Translate $z_0$ to $0$ and rotate $x,p$ together such that the first $k$
$p$-axes span $\ker d\pi_x(0)$. Since $T_0\Lambda$ is Lagrangian,
$dx_a|_{T_0\Lambda}=0$ for $1\le a\le k$, so
$(\xi,y)\coloneqq(p_1,\ldots,p_k,x_{k+1},\ldots,x_n)$ are local
coordinates. As $\Lambda$ is Lagrangian, the restriction of
$-\sum_{a=1}^k x_a\,dp_a+\sum_{i=k+1}^n p_i\,dx_i$ to $\Lambda$
is closed and hence locally equals $d\psi$ for some analytic $\psi$. Then
\[
x=(-\psi_\xi,y),\qquad p=(\xi,\psi_y),\qquad
D^2\psi=\begin{pmatrix}W&B\\ B^T&C\end{pmatrix}.
\]
By $\partial_{\xi_a}(x,p)(0)=(0,e_a)$ for $1\le a\le k$,
we have $W(0)=B(0)=0$.
From $D_{(\xi,y)}(x+ip)=\operatorname{diag}(iI_k,I_r)(I+iD^2\psi)$
and the constant phase of $\Lambda$, we obtain
\[
e^{i\Theta(D^2\psi)}
=i^{-k}\frac{\det D_{(\xi,y)}(x+ip)}{|\det D_{(\xi,y)}(x+ip)|}
=e^{i\Theta(D^2\psi(0))}.
\]
By continuity and $W(0)=B(0)=0$,
\[
\Theta(D^2\psi)=\Theta(D^2\psi(0))=\Theta(C(0))
=\sum_{|\theta_i|<\pi/2}\theta_i=\frac{r-2}{2}\pi.
\]

Fix $z\in\Lambda$ near $0$. Since
$i^k e^{i\Theta(D^2\psi)}=e^{i(n-2)\pi/2}$, $(\xi,y)$ are positively
oriented coordinates on $T_z\Lambda$. By continuity of the oriented
projection determinant, we approximate $T_z\Lambda$
by graphs over $x$ with $\Theta(D_xp)=(n-2)\pi/2$ on which
$(\xi,y)$ and $x$ induce the same orientation. The $j$th plane has
$x=(-W_j\xi-B_jy,y)$ and $p=(\xi,B_j^T\xi+C_jy)$, with
\[
\begin{pmatrix}W_j&B_j\\B_j^T&C_j\end{pmatrix}\longrightarrow D^2\psi(\xi,y)
\qquad\text{as }j\to\infty.
\]
Since the $j$th plane is a graph over $x$, $W_j$ is invertible and the
first $k$ rows and columns of $D_xp$ form $-W_j^{-1}$. By
\eqref{eq:two-smallest-angles},
\[
n_+(W_j)=n_-(-W_j^{-1})\le n_-(D_xp)\le1.
\]
Moreover, $\det(-W_j)>0$ because $(\xi,y)$ and $x$ induce the same
orientation, so $n_+(W_j)$ is even. Since $n_+(W_j)\le1$ and $W_j$ is invertible,
$W_j<0$, and hence $W\le0$.

For $H\in\operatorname{Sym}(n)$ near $D^2\psi(0)$, define
$s\coloneqq\mathscr F(H)$ by
\[
\Theta(H-s\operatorname{diag}(0_{k\times k},I_r))=\frac{r-2}{2}\pi.
\]
By $D\Theta>0$ and the implicit function theorem, $\mathscr F$ is
smooth with $D\mathscr F>0$. At $H=D^2\psi$, the phase is $(r-2)\pi/2$,
so $\mathscr F(D^2\psi)=0$.
For $t\in[0,1]$,
\[
\begin{aligned}
\det\left(I+i\begin{pmatrix}0&tB\\tB^T&C\end{pmatrix}\right)
&=\det\left[
  \begin{pmatrix}I_k&0\\-itB^T&I_r\end{pmatrix}
  \begin{pmatrix}I_k&itB\\itB^T&I_r+iC\end{pmatrix}
  \right] =\det(I_r+t^2B^TB+iC)\\
&=\det(I_r+t^2B^TB)
  \det\left(I_r+i(I_r+t^2B^TB)^{-1/2}C(I_r+t^2B^TB)^{-1/2}\right).
\end{aligned}
\]
Since $I_r+t^2B^TB>0$,
\[
\Theta\begin{pmatrix}0&tB\\tB^T&C\end{pmatrix}
-\Theta\left((I_r+t^2B^TB)^{-1/2}C(I_r+t^2B^TB)^{-1/2}\right)
\in 2\pi\mathbb Z.
\]
At $t=0$, $\Theta(\operatorname{diag}(0_{k\times k},C))-\Theta(C)=0$.
By continuity, at $t=1$ we have
\begin{equation}\label{eq:reduced-mixed-phase}
\Theta\begin{pmatrix}0&B\\B^T&C\end{pmatrix}
=\Theta((I_r+B^TB)^{-1/2}C(I_r+B^TB)^{-1/2}).
\end{equation}

For fixed $B$, $(I_r+B^TB)^{-1/2}C(I_r+B^TB)^{-1/2}$ is linear in $C$,
so \eqref{eq:reduced-mixed-phase} and \cite[Lemma 2.1]{YuanSupercritical}
imply that
\[
\left\{C\in\operatorname{Sym}(r):
\Theta\begin{pmatrix}0&B\\B^T&C\end{pmatrix}\ge\frac{r-2}{2}\pi\right\}
\]
is convex. For small $E\in\operatorname{Sym}(r)$, set
\[
s_\pm\coloneqq\mathscr F\begin{pmatrix}0&B\\B^T&C\pm E\end{pmatrix}.
\]
By the definition of $\mathscr F$,
$\Theta\bigl(\begin{psmallmatrix}0&B\\B^T&C\pm E-s_\pm I_r\end{psmallmatrix}\bigr)
=(r-2)\pi/2$. By convexity of the superlevel set,
\[
\Theta\begin{pmatrix}0&B\\ B^T&C-\frac{s_++s_-}{2}I_r\end{pmatrix}
\ge\frac{r-2}{2}\pi.
\]
Since $\Theta(H-s\operatorname{diag}(0_{k\times k},I_r))$ decreases
strictly in $s$,
\[
\mathscr F\begin{pmatrix}0&B\\B^T&C\end{pmatrix}
\ge \frac{s_++s_-}{2}.
\]
Thus $D_C^2\mathscr F\le0$ at $W=0$.
Since $\partial_{\xi_a}B=D_yW_{\cdot a}$,
\[
\left|\partial_{\xi_a}D^2\psi-
   \begin{pmatrix}0&0\\0&\partial_{\xi_a}C\end{pmatrix}\right|
=\left|\begin{pmatrix}
\partial_{\xi_a}W&D_yW_{\cdot a}\\
(D_yW_{\cdot a})^T&0
\end{pmatrix}\right|
\le C(n)|DW|.
\]
Differentiating $\mathscr F(D^2\psi)=0$ twice in $\xi_a$, we obtain
\begin{equation}\label{eq:partial-second-derivative}
\begin{aligned}
D\mathscr F(D^2\psi)[D^2W_{aa}]
&=-D^2\mathscr F(D^2\psi)[\partial_{\xi_a}D^2\psi,\partial_{\xi_a}D^2\psi]\\
&\ge -D_C^2\mathscr F(D^2\psi)[\partial_{\xi_a}C,\partial_{\xi_a}C]
  -C(n,\psi)|DW|\\
&\ge -D_C^2\mathscr F\bigl(\begin{smallmatrix}0&B\\B^T&C\end{smallmatrix}\bigr)
  [\partial_{\xi_a}C,\partial_{\xi_a}C]-C(n,\psi)(|W|+|DW|)\\
&\ge -C(n,\psi)(|W|+|DW|).
\end{aligned}
\end{equation}

If $W\not\equiv0$, let $W_m$ be the first nonzero term in the Taylor
expansion of $W=D_\xi^2\psi$ at $(\xi,y)=(0,0)$, homogeneous of degree $m$.
Then
\[
W(\epsilon\xi,\epsilon y)=\epsilon^m W_m(\xi,y)+O(\epsilon^{m+1})
\qquad\text{as }\epsilon\downarrow0.
\]
Since $W\le0$ and $W(0)=0$, we have $DW(0)=0$ and $W_m\le0$, so $m\ge2$.
Summing \eqref{eq:partial-second-derivative} over $a$ at
$(\epsilon\xi,\epsilon y)$ and multiplying by $\epsilon^{2-m}$,
\[
\begin{aligned}
 D\mathscr F(D^2\psi(0))[(D^2\tr W_m)(\xi,y)] 
=& \lim_{\epsilon\downarrow0}\epsilon^{2-m}
 D\mathscr F(D^2\psi(\epsilon\xi,\epsilon y))[(D^2\tr W)(\epsilon\xi,\epsilon y)]\\
\ge& -C(n,\psi)\lim_{\epsilon\downarrow0}\epsilon^{2-m}
 (|W|+|DW|)(\epsilon\xi,\epsilon y)\\
\ge& -\lim_{\epsilon\downarrow0} C(n,\psi,\xi,y)\epsilon=0.
\end{aligned}
\]
Since $\tr W_m\le0=\tr W_m(0,0)$, $\tr W_m$ has a maximum at $(\xi,y)=(0,0)$.
The strong maximum principle for $D\mathscr F(D^2\psi(0))[D^2(\cdot)]$
implies $\tr W_m=0$. As $W_m\le0$, we obtain $W_m=0$,
a contradiction. Thus $D_\xi^2\psi=0$ and
$\psi(\xi,y)=\xi\cdot\ell(y)+g(y)$, with analytic
$\ell(y)\coloneqq D_\xi\psi(0,y)$ and $g(y)\coloneqq\psi(0,y)$. Hence
\[
x=(-D_\xi\psi,y)=(-\ell(y),y),\qquad
p=(\xi,D_y\psi)=(\xi,D\ell(y)^T\xi+Dg(y)).
\]
Since $D_\xi x=0$ and the last $r$ rows of $D_yx$ form $I_r$,
$\operatorname{rank}d(\pi_x|_\Lambda)=r$ near $z_0$.

From $\psi(\xi,y)=\xi\cdot\ell(y)+g(y)$,
$B=D\ell(y)$ and $C=D^2g(y)+\sum_a\xi_aD^2\ell_a(y)$.
Suppose $r\ge3$, and let $e_a$ be the $a$th coordinate vector of $\R^k$.
Evaluating $B$ and $C$ at $(\xi+te_a,y)$, we obtain
\[
B(\xi+te_a,y)=B(\xi,y),\qquad
C(\xi+te_a,y)=C(\xi,y)+tD^2\ell_a(y)
\qquad\text{for }t\text{ near }0.
\]
Consequently
$(I_r+B^TB)^{-1/2}C(\xi+te_a,y)(I_r+B^TB)^{-1/2}=A+tE$, where
\[
A\coloneqq(I_r+B^TB)^{-1/2}C(\xi,y)(I_r+B^TB)^{-1/2},\qquad
E\coloneqq(I_r+B^TB)^{-1/2}D^2\ell_a(y)(I_r+B^TB)^{-1/2}.
\]
By \eqref{eq:reduced-mixed-phase},
$\Theta(A+tE)=\Theta(D^2\psi(\xi+te_a,y))=(r-2)\pi/2$ for $t$ near $0$.
$\Theta(A+tE)$ is real analytic in $t$ because
$\det(I+i(A+tE))\ne0$ for all $t\in\R$. Hence
$\Theta(A+tE)=(r-2)\pi/2$, $\forall t\in\R$.
By \eqref{eq:two-smallest-angles},
$\lambda_{r-1}(A+tE)+\lambda_r(A+tE)>0$. Since $r\ge3$,
$\sum_{i=2}^r\lambda_i(A+tE)>0$. Thus
\[
0<\max_i|\lambda_i(A+tE)|=\lambda_1(A+tE)
\le\tr(A+tE)=\tr A+t\tr E\qquad\forall t\in\R.
\]
In particular, $\tr(A+tE)>0$ for every $t\in\R$, so $\tr E=0$. Hence
$A+tE$ is bounded on $\R$, which implies $E=0$ and $D^2\ell_a(y)=0$
for each $a$ and $y$.
Therefore $\ell$ is affine, so $\{(-\ell(y),y)\}$ is open in an affine
$r$-plane.

When $r=2$, $\Theta(A)=0$ implies $\tr A=0$ for $A\in\operatorname{Sym}(2)$.
Thus \eqref{eq:reduced-mixed-phase} implies $\tr((I_2+B^TB)^{-1}C)=0$.
Since $B=D\ell(y)$, the metric of $x(y)=(-\ell(y),y)$ is $I_2+B^TB$.
Using $\partial_{\xi_a}C=D^2\ell_a(y)$, we compute the mean curvature vector
\[
\begin{aligned}
\vec H(x(y))
&=\left(-\bigl(\tr((I_2+B^TB)^{-1}D^2\ell_a(y))\bigr)_{a=1}^k,0\right)^\perp =\left(-D_\xi\tr\bigl((I_2+B^TB)^{-1}C\bigr),0\right)^\perp=0,
\end{aligned}
\]
where $(\cdot)^\perp$ is projection onto the normal space of the
surface $x=(-\ell(y),y)$ at $x(y)$.
Thus $\{(-\ell(y),y)\}$ is a minimal surface.
\end{proof}

\begin{proposition}\label{prop:finite-projection}
Let $n\ge3$ and let $\Lambda\subset\C^n$ be a connected real analytic
special Lagrangian submanifold with
$T_z\Lambda\in\Gr_{n}((n-2)\pi/2)$, $\forall z\in \Lambda$.
Suppose there exists $z_0\in\Lambda$ with
$\operatorname{rank}d(\pi_x|_\Lambda)_{z_0}=n$. Then
\[
\operatorname{rank}d(\pi_x|_\Lambda)_z=n
\qquad\forall z\in\Lambda.
\]
\end{proposition}
\begin{proof}
By analyticity and connectedness of $\Lambda$, if
$\det d(\pi_x|_\Lambda)=0$ on a nonempty open subset of $\Lambda$, then
$\det d(\pi_x|_\Lambda)=0$ on all of $\Lambda$, contrary to
$\operatorname{rank}d(\pi_x|_\Lambda)_{z_0}=n$.
If $\operatorname{rank}d(\pi_x|_\Lambda)_z=r$ with $0<r<n$ at some
$z\in\Lambda$, the rank is $r<n$ on a neighborhood of $z$ by
Lemma~\ref{lem:partial-propagation}. Thus $\det d(\pi_x|_\Lambda)$ is zero
on an open set, contrary to the rank at $z_0$.

Up to translation of $\Lambda$, suppose $\operatorname{rank}d(\pi_x|_\Lambda)_0=0$.
The tangent plane is $T_0\Lambda=\{0\}\times\R^n$, so $\Lambda$ is a
graph over $p$ near $0$.
Since $\Lambda$ is Lagrangian, $-x\cdot dp$ is closed on the graph,
so $x=-D\widehat u(p)$ for some analytic $\widehat u$ with
$\widehat u(0)=D\widehat u(0)=D^2\widehat u(0)=0$.
Since $D_p(p-ix)=-iD_p(x+ip)=I+iD^2\widehat u$, the constant phase of
$\Lambda$ implies
\[
e^{i\Theta(D^2\widehat u(p))}
=\frac{\det D_p(p-ix)}{|\det D_p(p-ix)|}
=e^{i\Theta(D^2\widehat u(0))}=1
\qquad\text{near }p=0.
\]

Since $\pi_p$ is invertible on $T_{(x(p),p)}\Lambda$, we approximate
$T_{(x(p),p)}\Lambda$ by planes of phase $(n-2)\pi/2$ that are graphs over
both $x$ and $p$. On the approximating planes, $n_+(-(D_xp)^{-1})\le1$
by \eqref{eq:two-smallest-angles}, and $-(D_xp)^{-1}\to D^2\widehat u(p)$
as the approximating planes converge to $T_{(x(p),p)}\Lambda$.
By continuity,
\begin{equation}\label{eq:vertical-graph-phase}
\Theta(D^2\widehat u)=0,\qquad n_+(D^2\widehat u)\le1
\qquad\text{near }p=0.
\end{equation}
If $\widehat u\equiv0$ near $p=0$, then $\operatorname{rank}d(\pi_x|_\Lambda)=0$
on some open subset of $\Lambda$, contradicting
$\operatorname{rank}d(\pi_x|_\Lambda)_{z_0}=n$ by analyticity. Let $P_d$,
$d\ge3$, be the first nonzero homogeneous Taylor term of $\widehat u$.
Since $\epsilon^{2-d}D^2\widehat u(\epsilon p)=D^2P_d(p)+O(\epsilon)$, by 
\eqref{eq:vertical-graph-phase}, $n_+(D^2P_d)\le1$.
Using $\arctan t=t+O(t^3)$ as $t\to0$, we obtain for fixed $p\in\R^n$
\[
\begin{aligned}
0&=\epsilon^{2-d}\Theta(D^2\widehat u(\epsilon p))\\
&=\epsilon^{2-d}\tr D^2\widehat u(\epsilon p)+O(\epsilon^{2d-4})
  \longrightarrow\Delta P_d(p)\qquad\text{as }\epsilon\downarrow0.
\end{aligned}
\]
Since $\Delta P_d=0$ and $n_+(D^2P_d)\le1$, we apply
\cite[Corollary 2]{GW} to $-P_d$ and choose orthogonal coordinates
$p=(z,y)\in\R^2\times\R^{n-2}$ in which
\begin{equation}\label{eq:binary-hessian}
P_d=\operatorname{Re}\bigl(c_d(z_1+iz_2)^d\bigr),\qquad
\det D_z^2P_d=-d^2(d-1)^2|c_d|^2|z|^{2d-4},\qquad c_d\ne0.
\end{equation}

Since $P_d$ is the first nonzero Taylor term of $\widehat u$,
$\epsilon^{-d}\widehat u(\epsilon z,\epsilon y)=P_d(z)+O(\epsilon)$
is analytic at $\epsilon=0$. Equation \eqref{eq:binary-hessian} implies
$\det D_z^2P_d(z)\ne0$ for $z\ne0$. By the implicit function theorem,  
\[
D_z\bigl[\epsilon^{-d}\widehat u(\epsilon z,\epsilon y)-z\cdot q\bigr]
=D_zP_d(z)-q+O(\epsilon)=0
\]
has an analytic solution $z=z(\epsilon,q,y)$ near $(0,D_zP_d(z),y)$.
Set $w_\epsilon(q,y)\coloneqq\epsilon^{-d}\widehat u(\epsilon z,\epsilon y)-z\cdot q$
at $z=z(\epsilon,q,y)$. Differentiating at fixed $q$, we have
$D_yz=-(\widehat u_{zz})^{-1}\widehat u_{zy}$
and $D_yw_\epsilon=\epsilon^{1-d}D_y\widehat u(\epsilon z,\epsilon y)$.
Differentiating $D_yw_\epsilon$ at fixed $q$ and substituting
$D_yz=-(\widehat u_{zz})^{-1}\widehat u_{zy}$, we obtain
\begin{equation}\label{eq:W-quotient}
W(q,y)\coloneqq D_y^2w_\epsilon
=\epsilon^{2-d}\bigl(\widehat u_{yy}
-\widehat u_{yz}(\widehat u_{zz})^{-1}\widehat u_{zy}\bigr),
\end{equation}
with the derivatives of $\widehat u$ evaluated at $p=\epsilon(z,y)$.
By $\widehat u_{zy}+\widehat u_{zz}D_yz=0$, we have
\[
\begin{pmatrix}I_2&0\\(D_yz)^T&I_{n-2}\end{pmatrix}
D^2\widehat u(\epsilon z,\epsilon y)
\begin{pmatrix}I_2&D_yz\\0&I_{n-2}\end{pmatrix}
=\begin{pmatrix}\widehat u_{zz}&0\\0&\epsilon^{d-2}W\end{pmatrix}.
\]
Since $\epsilon^{2-d}\widehat u_{zz}(\epsilon z,\epsilon y)\to D_z^2P_d(z)$
as $\epsilon\downarrow0$, by \eqref{eq:binary-hessian},
$\widehat u_{zz}$ has one positive and one negative eigenvalue for small
$\epsilon>0$. Thus
$n_+(W)=n_+(D^2\widehat u(\epsilon z,\epsilon y))-1=0$ by
\eqref{eq:vertical-graph-phase}, so $W\le0$. Moreover,
\[
\operatorname{rank}D_px(\epsilon(z,y))
=\operatorname{rank}D^2\widehat u(\epsilon z,\epsilon y)
=2+\operatorname{rank}W(q,y).
\]
If $W=0$ on some open set for fixed $\epsilon>0$, then
$\operatorname{rank}d(\pi_x|_\Lambda)=2<n$ on some open subset of $\Lambda$,
contradicting $\operatorname{rank}d(\pi_x|_\Lambda)_{z_0}=n$ by analyticity.

For $z\ne0$, let $h_\epsilon(z,y)\coloneqq\tr W(q,y)$ at
$q=\epsilon^{1-d}D_z\widehat u(\epsilon z,\epsilon y)$.
By \eqref{eq:W-quotient} and \eqref{eq:binary-hessian}, $h_\epsilon(z,y)$
is analytic in $\epsilon$ at $0$ and satisfies
\[
h_\epsilon(rz,ry)=r^{d-2}h_{r\epsilon}(z,y)
\qquad\forall r>0,\ \epsilon\text{ near }0.
\]
At $\epsilon=0$, $w_0=P_d(z)-z\cdot q$ and $q=D_zP_d(z)$, so
$D_yw_0=0$ and $h_0=\Delta_yw_0=0$. If $\tr W=0$ on an open set for
a fixed $\epsilon>0$, then $W=0$ on the same open set since $W\le0$.
Then $\operatorname{rank}D_px=2$ on the open set, contrary to
$\operatorname{rank}d(\pi_x|_\Lambda)_{z_0}=n$ by analyticity.
Thus the expansion of $h_\epsilon$ has a first nonzero coefficient
$h_m$ with $m\ge1$. Using $h_\epsilon\le0$ for $\epsilon>0$ and comparing
the $\epsilon^m$ terms in $h_\epsilon(rz,ry)=r^{d-2}h_{r\epsilon}(z,y)$, we obtain
\begin{equation}\label{eq:normalized-passive-block}
h_\epsilon(z,y)=\epsilon^m\bigl(h_m(z,y)+O(\epsilon)\bigr),
\qquad h_m(z,y)\le0,\qquad
h_m(rz,ry)=r^{m+d-2}h_m(z,y).
\end{equation}
In the $(q,y)$ coordinates, $D_y^2w_\epsilon=O(\epsilon^m)$ by
$|W|\le-\tr W$ and \eqref{eq:normalized-passive-block}.

For $w=w(q,y)$, set
\begin{equation}\label{eq:normalized-mixed-determinant}
\mathscr F_\epsilon(D^2w)\coloneqq
\epsilon^{2-d}\operatorname{Im}\det\begin{pmatrix}
-D_q^2w+i\epsilon^{d-2} I_2&-D_{qy}^2w\\
i\epsilon^{d-2}D_{yq}^2w&I_{n-2}+i\epsilon^{d-2}D_y^2w
\end{pmatrix}.
\end{equation}
At $\epsilon^{d-2}=0$, the determinant in \eqref{eq:normalized-mixed-determinant}
equals $\det(-D_q^2w)\in\R$. Hence the imaginary part of the determinant
is divisible by $\epsilon^{d-2}$, so $\mathscr F_\epsilon$ is polynomial in
$\epsilon$ and the entries of $D^2w$.
By $z=-D_qw_\epsilon$, 
$\epsilon^{1-d}D\widehat u(\epsilon z,\epsilon y)=(q,D_yw_\epsilon)$, and  \eqref{eq:vertical-graph-phase}, we have 
\[
\begin{aligned}
\epsilon^{d-2}\mathscr F_\epsilon(D^2w_\epsilon)
&=\operatorname{Im}\det D_{(q,y)}
   \bigl[(z,y)+i\epsilon^{-1}D\widehat u(\epsilon z,\epsilon y)\bigr]\\
&=\det D_{(q,y)}(z,y)\,
   \operatorname{Im}\det(I+iD^2\widehat u(\epsilon z,\epsilon y))
=0.
\end{aligned}
\] 
Next, let $w_j$ be the coefficient of $\epsilon^j$ in the Taylor expansion
of $w_\epsilon$ at fixed $(q,y)$ and let
$\bar w_\epsilon\coloneqq\sum_{j=0}^{m-1}\epsilon^j w_j$.
The estimate $D_y^2w_\epsilon=O(\epsilon^m)$ implies
$D_y^2w_j=0$ for $0\le j<m$, so $D_y^2\bar w_\epsilon=0$.
For a square matrix $A$, $\adj(A)$ is the transpose of the cofactor
matrix of $A$. When $D_y^2w=0$, we add $D_{qy}^2w$ times the last
$n-2$ rows to the first two in \eqref{eq:normalized-mixed-determinant}
and obtain
\[
\begin{aligned}
\mathscr F_\epsilon(D^2w)
&=\epsilon^{2-d}\operatorname{Im}\det\bigl(
   -D_q^2w+i\epsilon^{d-2}(I_2+D_{qy}^2w\,D_{yq}^2w)\bigr)\\
&=-\tr\bigl(\adj(D_q^2w)(I_2+D_{qy}^2w\,D_{yq}^2w)\bigr).
\end{aligned}
\]
Since  $D_y^2\bar w_\epsilon=0$,  $D_y(D_{qy}^2\bar w_\epsilon)=0$ and
$\Delta_y D_q^2\bar w_\epsilon=0$.
As $\adj$ is linear on $2\times2$ matrices,
$\Delta_y\mathscr F_\epsilon(D^2\bar w_\epsilon)=0$.
Since $w_\epsilon-\bar w_\epsilon=\epsilon^m w_m+O(\epsilon^{m+1})$ and
$D_y(D^2w_0)=0$, expanding $\mathscr F_\epsilon$ at
$D^2\bar w_\epsilon$ yields
\[
\begin{aligned}
0&=\lim_{\epsilon\downarrow0}\epsilon^{-m}\Delta_y
 \bigl[\mathscr F_\epsilon(D^2w_\epsilon)-\mathscr F_\epsilon(D^2\bar w_\epsilon)\bigr] 
=D\mathscr F_0(D^2w_0)[D^2\Delta_y w_m].
\end{aligned}
\]

When $D_{qy}^2w=0$, the coefficient of $i\epsilon^{d-2}$ in the
determinant in \eqref{eq:normalized-mixed-determinant} is
$(\det D_q^2w)\Delta_y w-\Delta_q w$. The determinant in
\eqref{eq:normalized-mixed-determinant} is unchanged when $D_{qy}^2w$
and $D_{yq}^2w$ both change sign. Hence the linearization at
$D_{qy}^2w_0=0$ has zero coefficients for mixed derivatives. Therefore
\[
0=D\mathscr F_0(D^2w_0)[D^2\Delta_y w_m]
=(\det D_q^2w_0)\Delta_y^2 w_m-\Delta_q\Delta_y w_m.
\]
By \eqref{eq:normalized-passive-block} and the expansion of $w_\epsilon$,
$h_m(z,y)=(\Delta_y w_m)(D_zP_d(z),y)$.
Using $\Delta_zD_zP_d=0$,
$(D_z^2P_d)^2=-(\det D_z^2P_d)I_2$ and
$D_q^2w_0=-(D_z^2P_d)^{-1}$, we obtain
\[
\begin{aligned}
\Delta_{z,y}h_m(z,y)
&=(\Delta_y^2 w_m)(q,y)
 -(\det D_z^2P_d(z))(\Delta_q\Delta_y w_m)(q,y)\\
&=(\det D_z^2P_d(z))D\mathscr F_0(D^2w_0)[D^2\Delta_y w_m](q,y)=0,
\end{aligned}
\]
where $q=D_zP_d(z)$.

Since $h_m\leq 0$ and $\log|z|$ is subharmonic with $\log|z|=-\infty$ on
$\{z=0\}$, the harmonic function $h_m$ extends subharmonically across
$\{z=0\}$ by \cite[Theorem 6.1]{HLRemovable}. By
\eqref{eq:normalized-passive-block}, the subharmonic extension of $h_m$
has value $0$ at $(z,y)=(0,0)$, since
\[
0\ge\limsup_{(z,y)\to0,\ z\ne0}h_m(z,y)
\ge\lim_{r\downarrow0}r^{m+d-2}h_m((1,0),0)=0.
\]
The subharmonic extension of $h_m$ has maximum $0$ at $(z,y)=(0,0)$,
so the strong maximum principle implies $h_m=0$, which contradicts  the choice of $h_m$.
\end{proof}

\section*{On the use of AI}
The author worked with a large language model in developing the construction and proofs, starting from research notes and drafts prepared by the author. The model also assisted with intermediate calculations, technical checks, and the exposition. The author independently verified all mathematical arguments and takes full responsibility for the results.

\end{document}